\documentclass[letterpaper]{article}
\usepackage{authblk}
\usepackage{url}
\usepackage[margin=1in,footskip=0.25in]{geometry}
\usepackage[labelfont=bf]{caption}
\usepackage{graphicx}
\usepackage{amsmath,amssymb,amsfonts,amsmath}
\usepackage{graphicx}
\usepackage{subfigure}
\usepackage{caption}
\usepackage{tikz}
\usepackage{float}
\usepackage[dvipdf]{epsfig}
\usepackage{float}
\usepackage{lipsum}
\usepackage{wrapfig}
\usepackage{cite}
\usepackage{color}
\usepackage{easylist}
\usepackage{lineno}
\usepackage{booktabs}
\usepackage{hyperref}
\usepackage{cleveref}
\usepackage{multirow}
\usepackage{hhline}
\usepackage{soul}
\usepackage{blindtext}
\usepackage{threeparttable}
\usepackage[bbgreekl]{mathbbol}

\def \L2{L^2}

\newcommand{\eqn}[1]{
	\begin{eqnarray}
	#1
	\end{eqnarray}
}

\newcommand{\tb}[1]{
	{\bf #1}
}

\newcommand{\norm}[1]{
	\left\lVert #1\right\rVert
}

\newcommand{\ip}[1]{
	\left\langle #1 \right\rangle
}

\newcommand{\px}{
	\partial_{x}
}
\newcommand{\py}{
	\partial_{y}
}

\newcommand{\pt}{
	\partial_{\theta}
}

\newcommand{\pph}{
	\partial_{\phi}
}

\newcommand{\pti}{
	\partial_{t}
}

\newcommand{\di}{
	\text{div}
}

\newcommand{\bs}[1]{
	\boldsymbol #1
}

\title{The Projection Extension Method: A Spectrally Accurate Technique for Complex Domains}

\author[1]{Saad Qadeer} 
\author[2]{Ehssan Nazockdast} 
\author[3]{Boyce E.~Griffith} 
\affil[1]{Advanced Computing, Mathematics and Data Division, Pacific Northwest National Laboratory, Richland, WA, 99354, United States of America}
\affil[2]{Department of Applied Physical Sciences, University of North Carolina, Chapel Hill, NC, 27599, United States of America}
\affil[3]{Departments of Mathematics, Applied Physical Sciences, and Biomedical Engineering, University of North Carolina, Chapel Hill, NC, USA}
\affil[ ]{{\tt saad.qadeer@pnnl.gov}, {\tt ehssan@email.unc.edu}, and {\tt boyceg@email.unc.edu}}

\begin{document}

\maketitle

\begin{abstract}
An essential ingredient of a spectral method is the choice of suitable bases for test and trial spaces. On complex domains, these bases are harder to devise, necessitating the use of domain partitioning techniques such as the spectral element method. In this study, we introduce the Projection Extension (PE) method, an approach that yields spectrally accurate solutions to various problems on complex geometries without requiring domain decomposition. This technique builds on the insights used by extension methodologies such as the immersed boundary smooth extension and Smooth Forcing Extension (SFE) methods that are designed to improve the order of accuracy of the immersed boundary method. In particular, it couples an accurate extension procedure, that functions on arbitrary domains regardless of connectedness or regularity, with a least squares minimization of the boundary conditions. The resulting procedure is stable under iterative application and straightforward to generalize to higher dimensions. Moreover, it rapidly and robustly yields exponentially convergent solutions to a number of challenging test problems, including elliptic, parabolic, Newtonian fluid flow, and viscoelastic problems.
\end{abstract}



\section{Introduction}\label{secintro}

Spectral methods are a group of highly accurate numerical techniques for solving partial differential equations (PDEs) \cite{boyd2001chebyshev,gottlieb1977numerical}. Broadly speaking, these methods discretize the weak formulation of the problem by drawing test and trial functions from finite-dimensional spaces. In the classical description, these functions are supported on the entire domain, with the boundary conditions either enforced on all the trial functions or imposed as additional constraints. Provided the bases are chosen with care, this leads to fast and robust solvers that yield provably high-order accurate solutions \cite{bernardi1997spectral}; in particular, for infinitely smooth exact solutions, the errors decay exponentially with the dimension of the approximation space. Examples of such techniques include Galerkin, collocation, and tau methods, among others \cite{canuto2007spectral}.

However, the choice of suitable test and trial bases for arbitrary domains is far from obvious, with the result that the applicability of spectral methods, in the basic formulation, is somewhat limited. Domain decomposition methods resolve this problem by forgoing the use of global test and trial functions, proceeding instead by partitioning the domain into smaller parts and applying spectral techniques locally. Different strategies for handling communication across sub-domain boundaries have been proposed, leading to the development of patching \cite{canuto2012spectral} and spectral element methods \cite{karniadakis2013spectral,patera1984spectral}. While these techniques enable the deployment of spectral methods on more general geometries and also allow parallel computation, they rely on the decomposition of the domain into high-quality elements. In cases where such partitions are difficult to generate in the first place, the effectiveness of these approaches is severely curtailed.

In this paper, we introduce the Projection Extension (PE) method, an approach for developing spectrally accurate methods on complex domains. The fundamental point of divergence from the conventional framework of spectral methods is that we eschew the use of variational formulations and exact application of boundary constraints. Instead, we extend the problem to a larger computational domain and reformulate it as an error minimization problem, where the contributions to the error are drawn from the PDE on the original domain and the boundary conditions. Searching for the minimizer in a finite-dimensional space then yields the discretized problem. As the numerical results demonstrate, the resulting procedure converges to the true solution exponentially as the dimension of the underlying space is increased. As this strategy is valid for arbitrary domains, in particular placing no constraints on the connectedness or boundary regularity, it represents a significant advance in the reach of spectral methods.

Our technique can also be seen as building on existing fixed Cartesian grid methods. These techniques generally operate by embedding the physical domain in a computational domain, appropriately extending the problem, and solving it on a fixed mesh. Among others, instances of this approach include the immersed boundary (IB) \cite{peskin1972flow,peskin1977numerical}, the immersed interface (II) \cite{leveque1994immersed,li2006immersed}, and the immersed boundary smooth extension (IBSE) \cite{stein2016immersed,stein2017immersed} methods. The IB method, for example, includes regularized delta-functions in the extended PDE to act as Lagrange multipliers for enforcing the boundary conditions. The II method builds on the IB method by accounting for jumps in the derivatives of the extended solution across the physical boundary, typically gaining an additional order of convergence. The IBSE method proceeds by positing an arbitrarily smooth extension to the solution, which is then used to generate the corresponding inhomogeneities on the computational domain and recover the solution itself. By varying the regularity of the extension, one can, in principle, obtain any algebraic rate of convergence. High convergence rates are desirable, for instance, in viscoelastic fluid flow problems, where the low order convergence provided by the IB method is unable to capture the stress tensor values accurately \cite{spagnolie2015complex}.

Recently, similar ideas have been employed in the design of some highly promising methods. The Fourier Continuation Alternate Direction (FC-AD) Implicit method pioneered by Bruno and Lyon \cite{bruno2010high,lyon2010high}, relies on the ADI procedure to reduce an evolution equation to a sequence of one-dimensional elliptic problems, which are extended by a highly accurate Fourier Continuation routine to the appropriate computational domains. The smooth selection embedding method \cite{agress2019smooth,agress2021smooth} attempts to solve the extension problem by formulating it as a Sobolev norm optimization problem. Another important contribution is the partition of unity extension approach developed in the context of boundary integral methods and applied to heat and fluid flow problems \cite{klinteberg2019fast,fryklund2019integral,fryklund2018partition}. Most pertinent to the current work is the Smooth Forcing Extension (SFE) method \cite{qadeer2021smooth}, which searches for extensions to the inhomogeneous terms in a finite-dimensional space by imposing regularity constraints at the nodes of the discretized physical boundaries. The rates of convergence are in perfect agreement with theory for one-dimensional problems and comfortably exceed the expected rates for two-dimensional domains. However, the imposition of smoothness across the physical boundaries by matching the normal derivatives poses some difficulties. Primarily, the appropriate choice of normal vectors for domains with corners is not obvious, barring some simple cases. In addition, for time-dependent problems, the SFE technique requires a truncation of the forcing term before computing the normal derivatives at the boundary to ensure stability. Despite yielding convergent solutions, this procedure involves differentiating discontinuous functions and hence is ill-defined.

In the context of extension methods, the PE method can be viewed as a consequence of a more sophisticated extension procedure. At its heart lies the insight that the PDE and boundary conditions do not need to be satisfied exactly, and that a least squares error minimization is a worthwhile alternative to ensure the smoothness of the extension. As a result, the applicability of this approach is not limited by boundary regularity. Moreover, the PE technique naturally evaluates the forcing term only on the physical domain and, because it eschews explicit regularity constraints, is mathematically well-defined. As we demonstrate, these properties easily resolve the ``mountain-in-fog'' conundrum \cite{boyd2002comparison}, and enable the method to be deployed on a diverse range of geometries, including discs, stars, multiply connected domains with corners, curved surfaces, and channels with obstacles. In addition, its flexibility and stability under iterative application can be harnessed to solve a vast array of problems, including elliptic and parabolic equations as well as Newtonian and viscoelastic fluid flow models.

\section{Mathematical Formulation}\label{mathform}

In this section, we outline our approach for a model problem. As discussed in the previous section, our technique can be arrived at in two equivalent ways. The insights gleaned from both approaches are helpful in extending the scope of applications as the ideas developed here are extended to systems of other types later in the paper.

Let $\Omega \subset \mathbb{R}^d$ be an arbitrary bounded domain with boundary $\partial \Omega$, and $\mathcal{L}$ be a linear differential operator. Consider the problem 
\eqn{
	\left\{ \begin{matrix}
		\mathcal{L} u = f, &  \text{on } \Omega,  \\
		u = g, &   \text{on } \partial \Omega.  \\
	\end{matrix} \right.
	\label{model}
}

One methodology for solving this problem requires embedding $\Omega$ in a simpler computational domain $\Pi$, which allows the use of fast solvers. The extension problem is concerned with extending the forcing $f$ to $\Pi$, so that the solution $u_\text{e}$ to the extended problem is smooth and its restriction to $\Omega$ is the same as that of (\ref{model}). 

Let $\{\phi_j\}_{j\in \mathcal{J}}$ be a family of linearly independent functions on $\Pi$. We write the extended forcing as $f_\text{e} = \sum_{j\in \mathcal{J}} c_j \phi_j$. To find the unknown coefficients $\tb{c} = \left(c_j\right)_{j\in \mathcal{J}}$, we impose two conditions:
\begin{itemize}
	\item[(i)] the solution $u_\text{e}$ to $\mathcal{L}u_\text{e} = f_\text{e}$ should obey $u_\text{e} = g$ at the boundary $\partial \Omega$;
	\item[(ii)] the extended forcing should agree with the given forcing on $\Omega$: $f_\text{e}|_{\Omega} = f$. 
\end{itemize} 

To impose (i), we first note that the solution is given by
\eqn{
	u_\text{e} = \mathcal{L}^{-1} f_\text{e} = \mathcal{L}^{-1} \sum_{j\in \mathcal{J}} c_j \phi_j =   \sum_{j\in \mathcal{J}} c_j \mathcal{L}^{-1}\phi_j, \label{ue1} 
}  
assuming (for now) that $\mathcal{L}$ is invertible. Denoting by $S^*$ the interpolation operator to the boundary $\partial \Omega$, we require $S^*u_\text{e} = g$.  Setting $\theta_j := S^*\mathcal{L}^{-1} \phi_j$ for $j \in \mathcal{J}$, this condition translates to
\eqn{
	\sum_{j\in \mathcal{J}} c_j\theta_j = g. \label{ue2}
}

To apply both conditions, we define the objective function
\eqn{
	G(\tb{c}) = \norm{\sum_{j\in \mathcal{J}} c_j\theta_j - g}^2_{L^2(\partial\Omega)} + \norm{\sum_{j\in \mathcal{J}} c_j\phi_j - f}^2_{L^2(\Omega)} , \label{Gdefn} 
}
and find the $\tb{c}$ that minimizes $G(\tb{c})$. This can be accomplished by setting $\partial G/\partial c_l = 0$ for $l \in \mathcal{J}$; this yields
\eqn{
	\ip{\theta_l,\sum_{j\in \mathcal{J}} c_j\theta_j - g}_{\partial \Omega} + \ip{\phi_l,\sum_{j\in \mathcal{J}} c_j\phi_j - f}_{\Omega} = 0. \label{Gdiff}
}

Upon rearrangement, (\ref{Gdiff}) becomes
\eqn{
	\sum_{j\in \mathcal{J}}\left(\ip{\theta_l,\theta_j}_{\partial\Omega} + \ip{\phi_l,\phi_j}_{\Omega}\right)c_j = \ip{\theta_l,g}_{\partial\Omega} + \ip{\phi_l,f}_{\Omega} . \label{sqsys}
}

These equations form a self-adjoint linear system that can be used to solve for the unknown coefficients $\tb{c}$. Provided that the dimension $|\mathcal{J}|$ of the lower-dimensional subspace is sufficiently large, the resulting solution $u_\text{e} = \sum_{j\in \mathcal{J}} c_j\phi_j$ should agree with the actual solution of (\ref{model}) on $\Omega$. 

\subsection{Non-invertible $\mathcal{L}$}\label{NinvL}
In the case that $\mathcal{L}$ fails to be invertible on $\Pi$, we modify the formulation slightly. We assume that $\mathcal{L}$ is self-adjoint on $\Pi$, and let $\{\psi_k\}_{k\in\mathcal{K}}$ be an orthonormal basis for the null space of $\mathcal{L}$. First note that if $\mathcal{L}u_\text{e} = f_\text{e}$, then for all $k \in \mathcal{K}$, 
\eqn{
	\ip{\psi_k,f_\text{e}}_\Pi = \ip{\psi_k,\mathcal{L}u_\text{e}}_\Pi = \ip{\mathcal{L}\psi_k,u_\text{e}}_\Pi = 0, \label{fepsik}
}
using the self-adjoint property of $\mathcal{L}$. Next, define 
\eqn{
	d_k = \ip{\psi_k,u_\text{e}}_\Pi, \label{dkdef}
}
and set 
\eqn{
	u_0 = u_\text{e} - \sum_{d \in \mathcal{K}} d_k \psi_k. \label{u0def}
} 

In particular, note that for all $m \in \mathcal{K}$,
\eqn{
	\ip{\psi_m,u_0}_\Pi = \ip{\psi_m , u_\text{e} - \sum_{d \in \mathcal{K}} d_k \psi_k}_\Pi = \ip{\psi_m,u_\text{e}}_\Pi - \sum_{k \in \mathcal{K}} d_k\ip{\psi_m,\psi_k}_\Pi = 0. \label{u0eff}
}

Observe that $u_0$ is simply the orthogonal projection of $u_\text{e}$ onto $\mathcal{R}$, the range of $\mathcal{L}$. Set $\mathcal{A} = \mathcal{L}|_{\mathcal{R}}$ so that
\eqn{
	u_0 = \mathcal{A}^{-1}f_\text{e}. \label{Adef}
} 

To apply the additional constraints given by (\ref{fepsik}) on the extended forcing $f_\text{e}$, and to find the unknowns $\tb{d} = (d_k)_{k \in \mathcal{K}}$ to reconstruct $u_\text{e}$, we modify the objective function to
\eqn{
	G(\tb{c},\tb{d}) = \norm{\sum_{j\in \mathcal{J}} c_j\theta_j + \sum_{k\in \mathcal{K}} d_k\gamma_k - g}^2_{L^2(\partial\Omega)} + \norm{\sum_{j\in \mathcal{J}} c_j\phi_j - f}^2_{L^2(\Omega)}  + \sum_{k \in \mathcal{K}} \left|\sum_{j \in \mathcal{J}} c_j\ip{\psi_k,\phi_j}_\Pi\right|^2, \label{G2defn} 
}
where $\theta_j := S^*\mathcal{A}^{-1}\phi_j$ and $\gamma_k := S^*\psi_k$. Setting $\partial G/\partial c_l = 0$ for $l \in \mathcal{J}$ yields
\eqn{
	\sum_{j\in \mathcal{J}}\left(\ip{\theta_l,\theta_j}_{\partial\Omega} + \ip{\phi_l,\phi_j}_{\Omega} + \sum_{k \in \mathcal{K}} \ip{\psi_k,\phi_l}_\Pi\ip{\psi_k,\phi_j}_\Pi\right)c_j &+& \sum_{k \in \mathcal{K}} \ip{\theta_l,\gamma_k}_{\partial\Omega}d_k =\nonumber\\
	&& \ip{\theta_l,g}_{\partial\Omega} + \ip{\phi_l,f}_{\Omega} , \label{sqsys2a}
}
while $\partial G/\partial d_m = 0$ for $m \in \mathcal{K}$ gives
\eqn{
	\sum_{j \in \mathcal{J}} \ip{\gamma_m,\theta_j}_{\partial \Omega} c_j + \sum_{k \in \mathcal{K}} \ip{\gamma_m,\gamma_k}_{\partial \Omega} d_k = \ip{\gamma_m,g}_{\partial\Omega}. \label{sqsys2b}
}

Observe that these equations together again form a self-adjoint linear system that can be used to solve for the unknowns $\tb{c}$ and $\tb{d}$ and obtain $u_\text{e}$ according to (\ref{u0def}) and (\ref{Adef}).

\subsection{An Alternative Formulation}\label{AltForm}

Note that the preceding formulation does not, in principle, impose any smoothness constraints on the basis functions $\{\phi_j\}_{j \in \mathcal{J}}$ and, as a result, is valid for piecewise continuous functions as well. In practice, however, orthogonal polynomials or trigonometric functions are superior choices because of their rapid convergence properties. The smoothness of these functional families can then be leveraged to provide a greatly simplified alternative formulation. 

We expand the solution $u_\text{e}$ in terms of the family of smooth linearly independent functions $\{\phi_j\}_{j \in \mathcal{J}}$ as $u_\text{e} = \sum_{j \in \mathcal{J}} \hat{c}_j \phi_j$. Imposing the agreement of the corresponding forcing and boundary conditions requires us to minimize the objective function
\eqn{
	H(\tb{\hat c}) = \norm{\sum_{j\in \mathcal{J}} \hat{c}_j S^*\phi_j - g}^2_{L^2(\partial\Omega)} + \norm{\sum_{j\in \mathcal{J}} \hat{c}_j \mathcal{L}\phi_j - f}^2_{L^2(\Omega)}. \label{Hdefn} 
}

Setting $\upsilon_j := S^*\phi_j$ and requiring $\partial H/\partial \hat{c}_l = 0$ for $l \in \mathcal{J}$ yields
\eqn{
	&& \ip{\upsilon_l,\sum_{j\in \mathcal{J}} \hat{c}_j\upsilon_j - g}_{\partial \Omega} + \ip{\mathcal{L}\phi_l,\sum_{j\in \mathcal{J}} \hat{c}_j \mathcal{L}\phi_j - f}_{\Omega} = 0 \nonumber\\
	&\Rightarrow &   \sum_{j\in \mathcal{J}} \left(\ip{\upsilon_l,\upsilon_j }_{\partial \Omega} +  \ip{\mathcal{L}\phi_l,\mathcal{L}\phi_j}_{\Omega}\right)\hat{c}_j = \ip{\upsilon_l,g}_{\partial \Omega} + \ip{\mathcal{L}\phi_l,f}_{\Omega}.   \label{Hdiff}
}

Observe that this is also a self-adjoint linear system of exactly the same form as (\ref{sqsys}). This approach essentially differs from the previous one by not requiring the inversion of $\mathcal{L}$. An obvious advantage this confers is that, in case $\mathcal{L}$ fails to be invertible, we do not require the specialized treatment outlined in Subsection \ref{NinvL}. Another subtle benefit is that it affords us the freedom to choose a basis $\{\phi_j\}_{j \in \mathcal{J}}$ without having to deal with computing $\mathcal{L}^{-1}\phi_j$ for all $j$. For example, in the case where the basis functions are eigenfunctions of $\mathcal{L}$, such operations are trivial to execute. For more general bases, however, these terms comprise an additional layer of computation. Evaluating $\mathcal{L}\phi_j$, on the other hand, is a much more tractable task, assuming $\{\phi_j\}_{j \in \mathcal{J}}$ are sufficiently smooth, yielding significant savings in computational time and effort. We build on these ideas in Subsection \ref{AltStokes} to develop solvers for Stokes equations in more challenging settings, which are then used in Section \ref{ViscoFluids} to great effect.

\subsection{Numerical Implementation}\label{NumInt}
We briefly discuss the details underlying the numerical implementation of our technique. In the discussion that follows, we use equations (\ref{sqsys2a}) and (\ref{sqsys2b}) since they form a more general system than (\ref{sqsys}) or (\ref{Hdiff}). We also assume that the null space of $\mathcal{L}$ has finite dimension $K$; clearly, this assumption is not required if we use the formulation described in Subsection \ref{AltForm}.

The appropriate choice of the computational domain $\Pi$ depends on the form of the problem and the physical domain $\Omega$. Unless otherwise specified, we take $\Pi$ to be the $d$-dimensional torus $\mathbb{T}^d$. The obvious choice for the extension functions $\{\phi_j\}$ then is the Fourier basis $\{e^{i\tb{j}\cdot \tb{x}}\}_{\tb{j} \in \mathbb{Z}^d}$. This forms a basis for $C(\Pi)$, which in turn is dense in $L^2(\Pi)$. Moreover, this basis diagonalizes several ubiquitous differential operators, such as the Laplacian and Helmholtz operators, so computing $\mathcal{L}^{-1}\phi_j$ is particularly straightforward.       

In practice, we fix the highest frequency $N_\text{e}$ and consider the extension function $\{e^{i\tb{j}\cdot \tb{x}}\}_{\tb{j} \in \mathcal{J}(N_\text{e})}$, where $\mathcal{J}(N_\text{e}) = \{\tb{j} \in \mathbb{Z}^d: \max_{1 \leq i \leq d} |j_i| \leq N_\text{e}\}$. The computational domain is discretized by $N$ equidistant grid-points along each dimension $\{\tb{x}_{\tb{k}}\}_{\tb{k} \in \Gamma(N)}$, where $\Gamma(N) = \{\tb{k} \in \mathbb{N}^d: \max_{1 \leq i \leq d} |k_i| \leq N\}$. In addition, we define $\Omega_N$ to be the indices of the grid points that lie in $\Omega$ and denote the grid-spacing by $h = 2\pi/N$. Finally, the physical boundary $\partial \Omega$ is discretized by equidistant points $\{\tb{s}_i\}_{1 \leq i \leq n_\text{b}}$ with spacing $\Delta s$ along the arc length. The respective inner products can then be approximated by
\eqn{
	\ip{\alpha_1,\alpha_2}_{\partial \Omega} &=& \int_{\partial \Omega} \overline{\alpha_1(\tb{s})} \alpha_2(\tb{s}) \ ds \approx (\Delta s) \sum_{i = 1}^{n_\text{b}} \overline{\alpha_1(\tb{s}_i)} \alpha_2(\tb{s}_i),  \label{POip}\\
	\ip{\beta_1,\beta_2}_{\Pi} &=& \int_{\Pi} \overline{\beta_1(\tb{x})} \beta_2(\tb{x}) \ d\tb{x} \approx h^d \sum_{\tb{k} \in \Gamma(N)} \overline{\beta_1(\tb{x}_\tb{k})} \beta_2(\tb{x}_\tb{k}),  \label{Cip} \\
	\ip{\delta_1,\delta_2}_{\Omega} &=& \int_{\Pi} \chi_{\Omega}(\tb{x})\overline{\delta_1(\tb{x})} \delta_2(\tb{x}) \ d\tb{x} \approx h^d \sum_{\tb{k} \in \Omega_N}\overline{\delta_1(\tb{x}_\tb{k})} \delta_2(\tb{x}_\tb{k}).  \label{Oip}
}

Define the matrices $P$, $Q$, $R$ and $S$ by
\eqn{
	P_{i\tb{j}} = \theta_\tb{j}(\tb{s}_i)(\Delta s)^{1/2}, \quad Q_{\tb{k}\tb{j}} = \phi_{\tb{j}}(\tb{x}_{\tb{k}}) h^{d/2}, \quad R_{k\tb{j}} = \ip{\psi_k,\phi_{\tb{j}}}_\Pi, \quad S_{ik} = \gamma_k(\tb{s}_i)(\Delta s)^{1/2}, \label{pqrsmats}
}
for $1 \leq i \leq n_\text{b}$, $\tb{j} \in \mathcal{J}(N_\text{e})$, $\tb{k} \in \Omega_N$ and $1 \leq k \leq K$ (an appropriate ordering is applied to the multi-indices when we use them as row or column indices). Similarly, define the vectors $\tb{g}$ and $\tb{f}$, with entries
\eqn{
	g_i = g(\tb{s}_i)(\Delta s)^{1/2}, \quad f_{\tb{k}} = f(\tb{x}_{\tb{k}})h^{d/2}, \label{gfs}
}
for $1 \leq i \leq n_\text{b}$ and $\tb{k} \in \Omega_N$. Setting
\eqn{
	A = \begin{pmatrix}
		P & S \\ Q & 0_{|\Omega_N| \times K} \\ R & 0_{K \times K}
	\end{pmatrix}, \quad \tb{r} = \begin{pmatrix}
		\tb{g} \\ \tb{f} \\ 0_{K \times 1} \\
	\end{pmatrix}, \quad \tb{z} = \begin{pmatrix}
		\tb{c} \\ \tb{d}
	\end{pmatrix} \label{Arzdef}
}
then allows us to write the system (\ref{sqsys2a},\ref{sqsys2b}) as
\eqn{
	A^*A\tb{z} = A^*\tb{r}, \label{normeqs}
}
as the matrix product ends up computing the approximations to the inner products, as given by (\ref{POip}) and (\ref{Oip}). Observe that $A$ is a matrix of size $(n_\text{b} + |\Omega_N| + K) \times (|\mathcal{J}(N_\text{e})| + K)$. In general, we need $n_\text{b}$ and $N$ to be large so the inner product approximations are computed accurately. As a result, (\ref{normeqs}) can be identified as the normal equations for the overdetermined problem $A\tb{z} = \tb{r}$. This can be solved efficiently in the least squares sense by using the $QR$ decomposition of $A$. This approach has the benefit of not squaring the condition number, as (\ref{normeqs}) does, while also allowing us to solve this problem a large number of times where $\tb{r}$ may vary over iterations, as is the case, for example, for time-dependent problems.

\section{Numerical Results for Elliptic and Parabolic Problems}\label{SimpTests}
We begin by applying the PE method to the one-dimensional Poisson equation
\eqn{
	\left\{
	\begin{matrix}
		u''(x) = 1/x , \quad x \in (2,5), \\
		u(2) = 1, \ u(5) = -1. 
	\end{matrix} \right. \label{Poi1D}
}

We use equations (\ref{sqsys2a}) and (\ref{sqsys2b}) for this example to illustrate the effectiveness of both the extension and solution methodologies. For functions defined on $C = \mathbb{T}$, the Laplacian operator $\partial_x^2$ has a one-dimensional null-space composed of constant functions. We set $N = 2^{10}$ and allow $N_\text{e}$ to vary; note that the number of grid-points simply controls the accuracy of the numerical integrations in (\ref{Cip}) and (\ref{Oip}). The convergence plots for the corresponding extensions $f_\text{e}^{[N_\text{e}]}$ and solutions $u_\text{e}^{[N_\text{e}]}$ are shown in Figure \ref{Poi1Dpl}. Observe that both converge geometrically to the true functions with increasing values of $N_\text{e}$. In particular, the rapid convergence of the forcing demonstrates that searching for a minimizer in a finite-dimensional space, as in (\ref{Gdefn}), is a powerful method for extending arbitrary functions to larger computational domains. Note also that the specific forcing term possesses a singularity in the extension region. The fact that our algorithm still manages to produce a highly accurate extension shows that it circumvents the ``mountain-in-fog'' problem \cite{boyd2002comparison}. This is essentially a consequence of only using the values of the function on the physical domain, as in the vector $\tb{f}$ in (\ref{Arzdef}). This compares very favorably with the approach employed in the SFE method that required knowledge of the derivatives of the function at the boundary and implicitly made use of the full functional form.   

\def \sclam{0.21}
\def \sclbm{0.14}

\begin{figure}[tbph]
	\centering
	\subfigure[Convergence of $f_\text{e}^{[N_\text{e}]}$]
	{\includegraphics[scale=\sclbm]{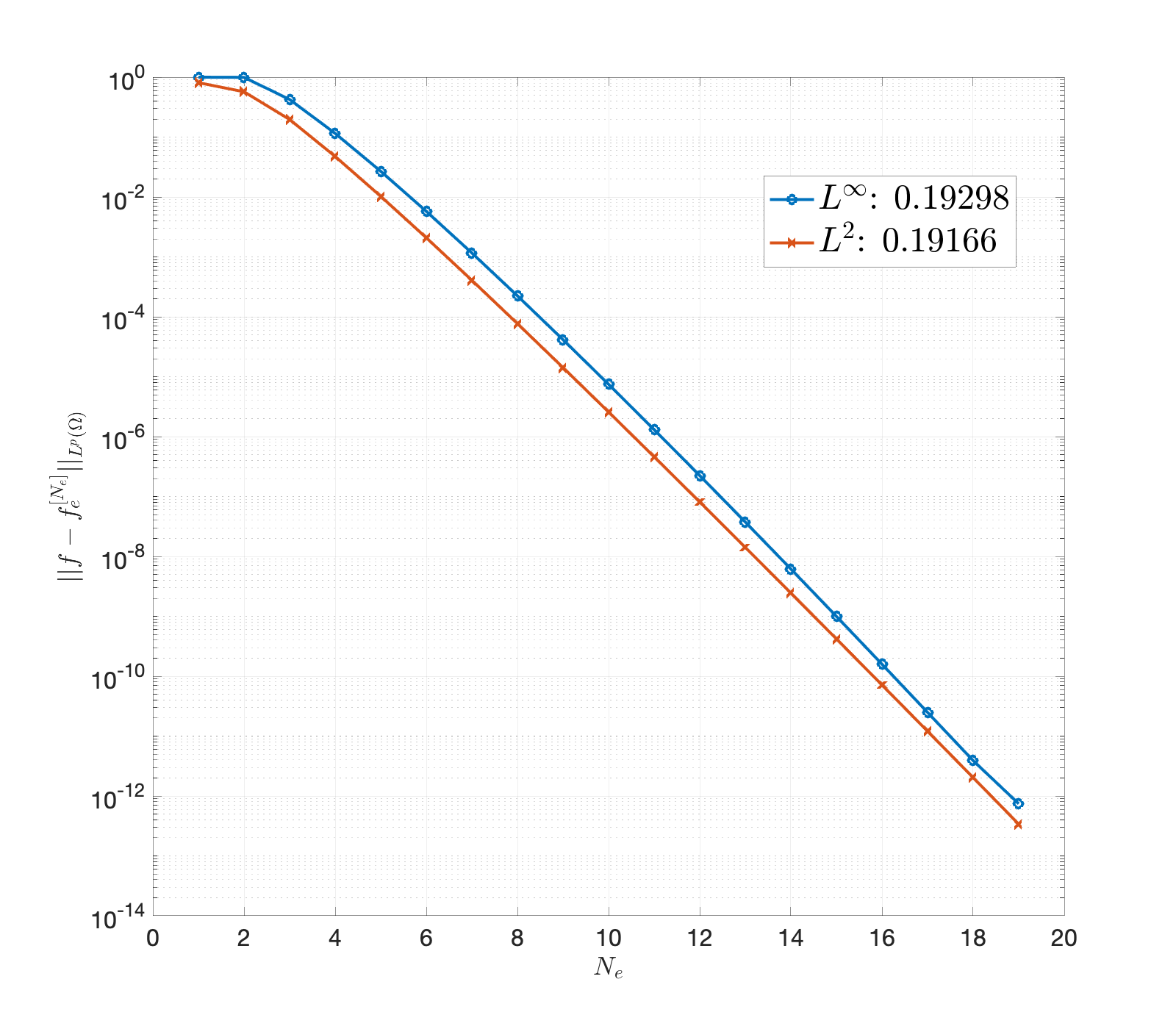}
		\label{Poi1Df}
	}
	\subfigure[Convergence of $u_\text{e}^{[N_\text{e}]}$]
	{\includegraphics[scale=\sclbm]{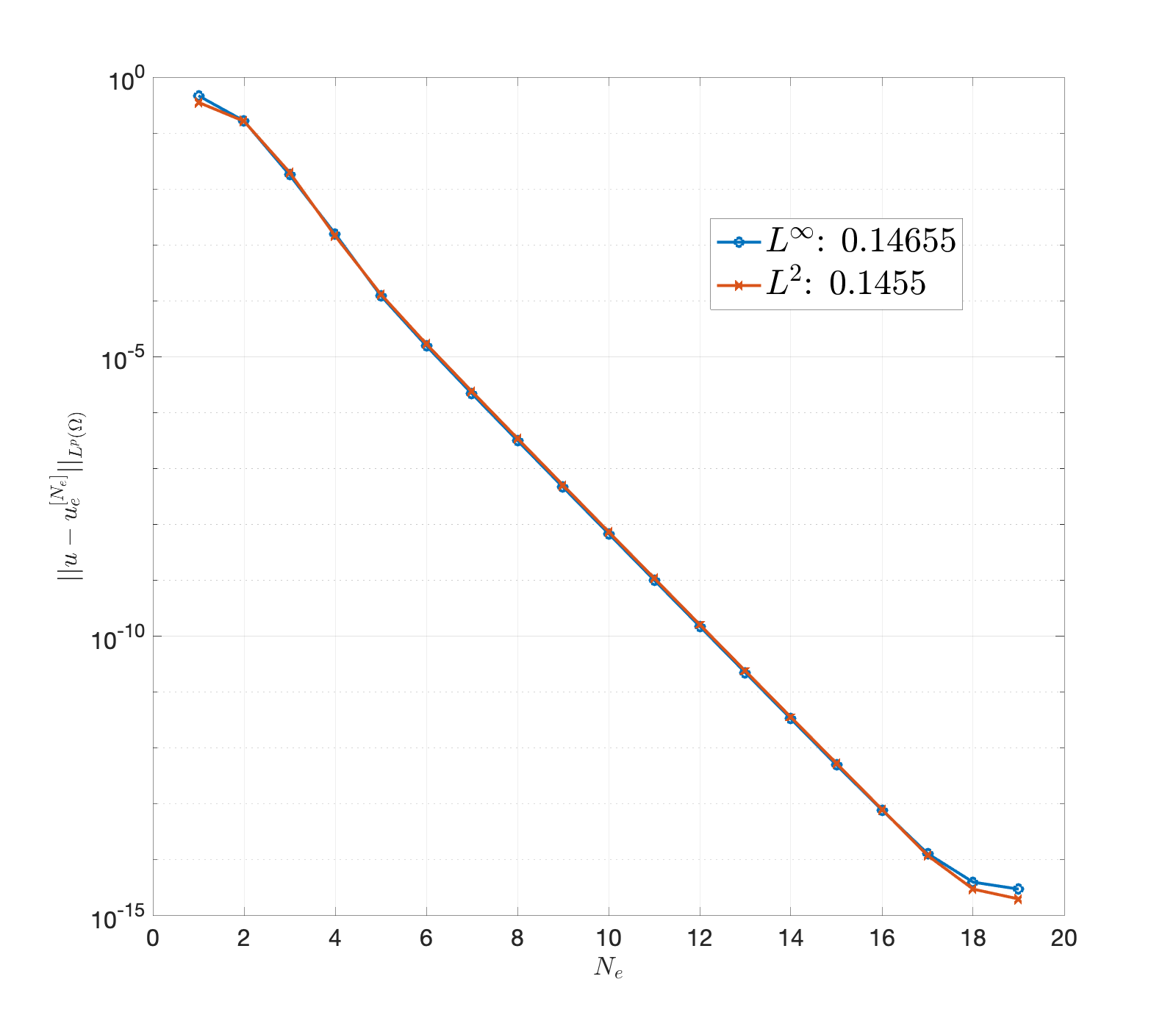}
		\label{Poi1Du}
	}
	\caption{Convergence plots for the extended forcing and solutions in $L^{\infty}(\Omega)$ and $L^2(\Omega)$ norms, for the one-dimensional Poisson problem (\ref{Poi1D}). Observe that both converge exponentially at $O(a^{N_\text{e}})$; the corresponding best-fit values of $a$ are shown for the different norms.}\label{Poi1Dpl}
\end{figure}

Next, we adapt our method to solve the heat equation. Consider the more general time-dependent problem 
\eqn{
	\left\{ \begin{matrix}
		u_t - \mathcal{L} u = f(t,x), &  \text{for }x \in \Omega,  \\
		u(0,x) = u_0(x), &  \text{for }x \in \Omega,  \\
		u(t,x) = g(t,x),  & \text{for }x \in  \partial \Omega, \ t > 0.   \\
	\end{matrix} \right.
	\label{tdepeqn}	
}

We employ the following iteration scheme, obtained from the four-step Backward Differentiation Formula (BDF-4), to discretize the time derivative:
\eqn{
	\left(\mathbb{I} - \frac{12\Delta t}{25} \mathcal{L}\right)u^{n+1} = \frac{12\Delta t f^{n+1} + 48u^{n} - 36u^{n-1} + 16u^{n-2} - 3u^{n-3}}{25},	\label{bdf4}
}

Setting
\eqn{
	\hat{ \mathcal{L}} = \mathbb{I} - \frac{12\Delta t}{25} \mathcal{L}, \qquad F^{n+1} = \frac{12\Delta t f^{n+1} + 48u^{n} - 36u^{n-1} + 16u^{n-2} - 3u^{n-3}}{25}	\label{Lhat}
}
allows us to write (\ref{bdf4}) as
\eqn{
	\hat{\mathcal{{L}}}u^{n+1} = F^{n+1} \label{bdf41}	
}
with corresponding boundary conditions $u^{n+1}(x) = g((n+1)\Delta t,x)$, for $x \in \partial \Omega$. Thus, this scheme lends itself naturally to the approach described earlier. The matrices corresponding to equations (\ref{qrform}) need to be built just once (for a specified $\Delta t$) for the entirety of a simulation. 

\begin{figure}[tbph]
	\centering
	\subfigure[Convergence of $u_\text{e}^{[N_\text{e}]}$]
	{\includegraphics[scale=\sclbm]{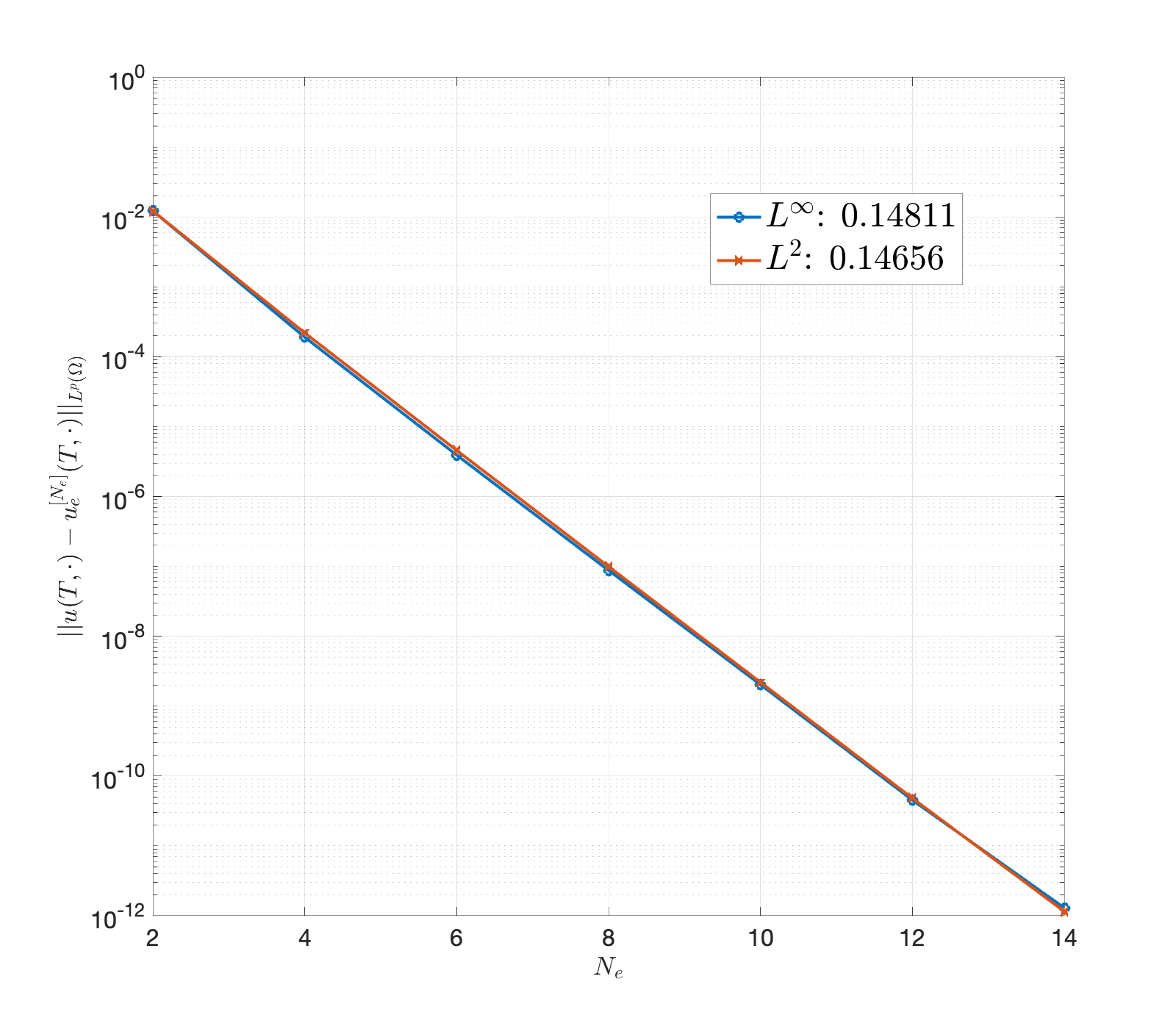}
		\label{heat1Dpl}
	}
	\subfigure[Evolution of errors with $t$ for $N_\text{e} = 14$]
	{\includegraphics[scale=\sclbm]{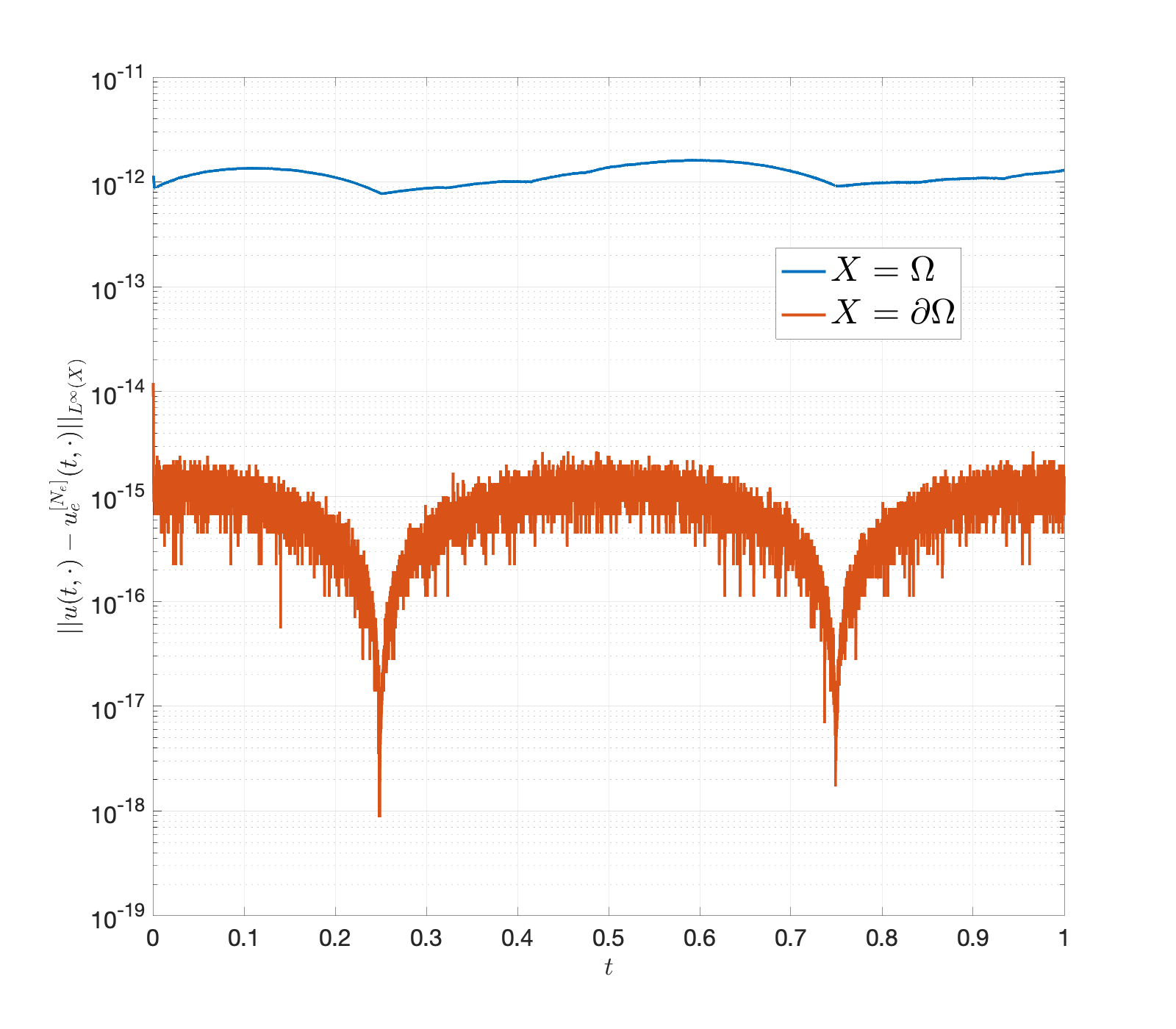}
		\label{heat1DT}
	}
	\caption{Results for the one-dimensional heat equation solved up to $T = 1$. (a) Exponential convergence of the form $O(a^{N_\text{e}})$ with the estimated values $a$ shown for the different norms. (b) Time discretization errors are negligible throughout, as seen here by tracking the $L^{\infty}(\Omega)$ error in the computed solution, as well as the deviation in the boundary values.}\label{heat1D}
\end{figure}

For the heat equation, we have $\mathcal{L} = \Delta$ so (\ref{bdf41}) reduces to a Helmholtz-type equation. On the one-dimensional torus $\mathbb{T}$, the operator $\hat{\mathcal{L}} = I - (12\Delta t/25)\Delta$ is invertible so we do not need to make provisions for the null-space. We test our technique by using the forcing and boundary conditions derived from the exact solution $u(t,x) = \ln(x)\cos(2\pi t)$ on $\Omega = (2,5)$ for $t \in [0,1]$. Note that our time-stepping scheme requires values of the solution for the first three steps to be provided. In the results shown here, we have used the exact solution for that; alternatively, we could have employed the backward Euler scheme for initializing the multistep method. 

Figure \ref{heat1Dpl} shows results with $N = 2^{8}$ and $\Delta t = 10^{-4}$. The high-order time-stepping routine we have used ensures that the time integration error is negligible, as shown in Figure \ref{heat1DT}, allowing for a comprehensive test of the accuracy and stability properties of our algorithm. The errors can again seen to decay exponentially with $N_\text{e}$, even when the exact solution and forcing possess singularities in the extension region. More importantly, this demonstrates that iterative application of our technique is stable without the introduction of any specific modifications, unlike the SFE method \cite{qadeer2021smooth}.

\begin{figure}[tbph]
	\centering
	\subfigure[The star-shaped domain $\Omega_1$]
	{\includegraphics[scale=\sclbm]{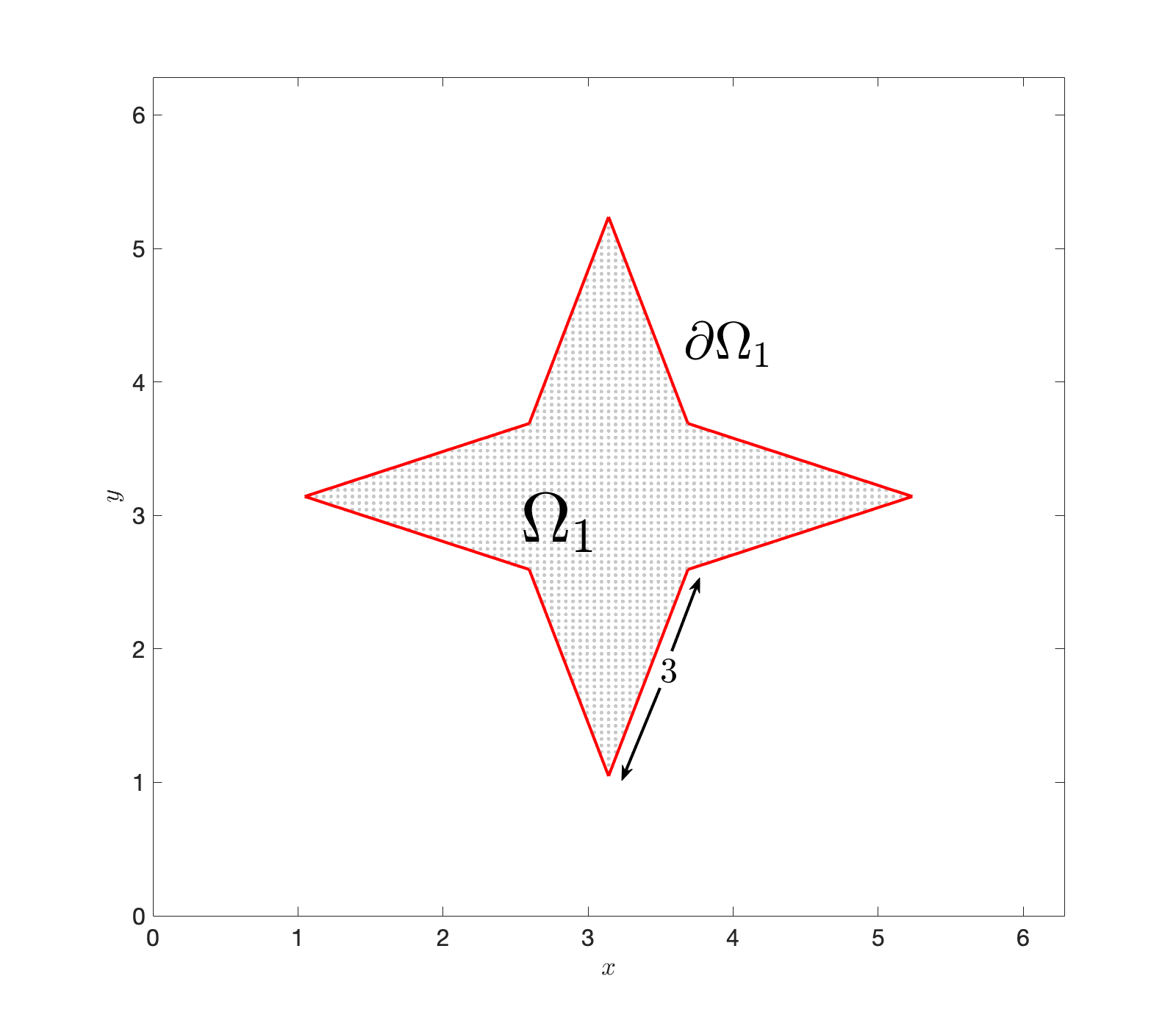}
		\label{starshape}
	}
	\subfigure[The extended solution for $N_\text{e} = 30$]
	{\includegraphics[scale=\sclbm]{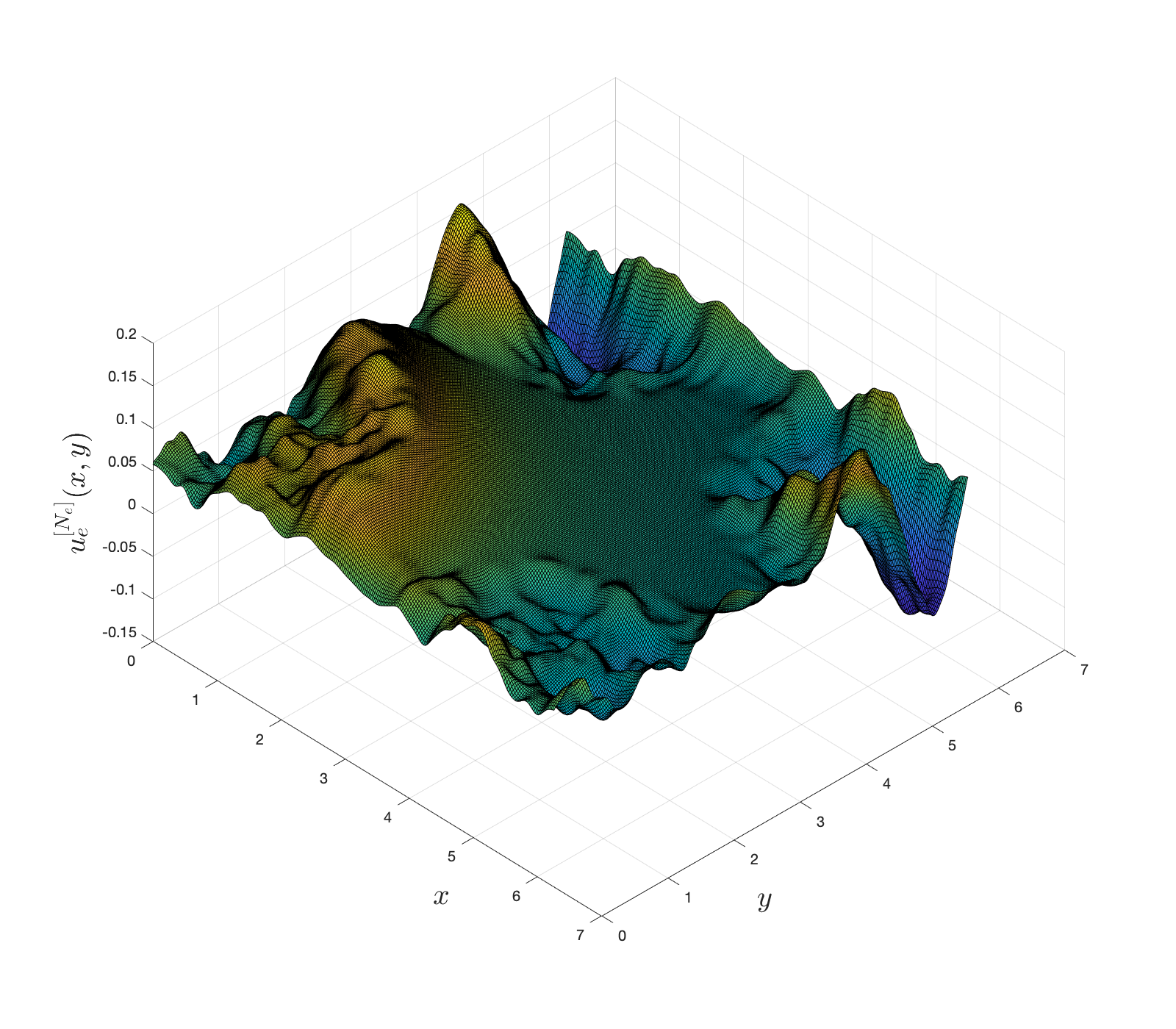}
		\label{Poi2Dudiam}
	}
	\subfigure[Convergence for problem (\ref{Poi2D})]
	{\includegraphics[scale=\sclbm]{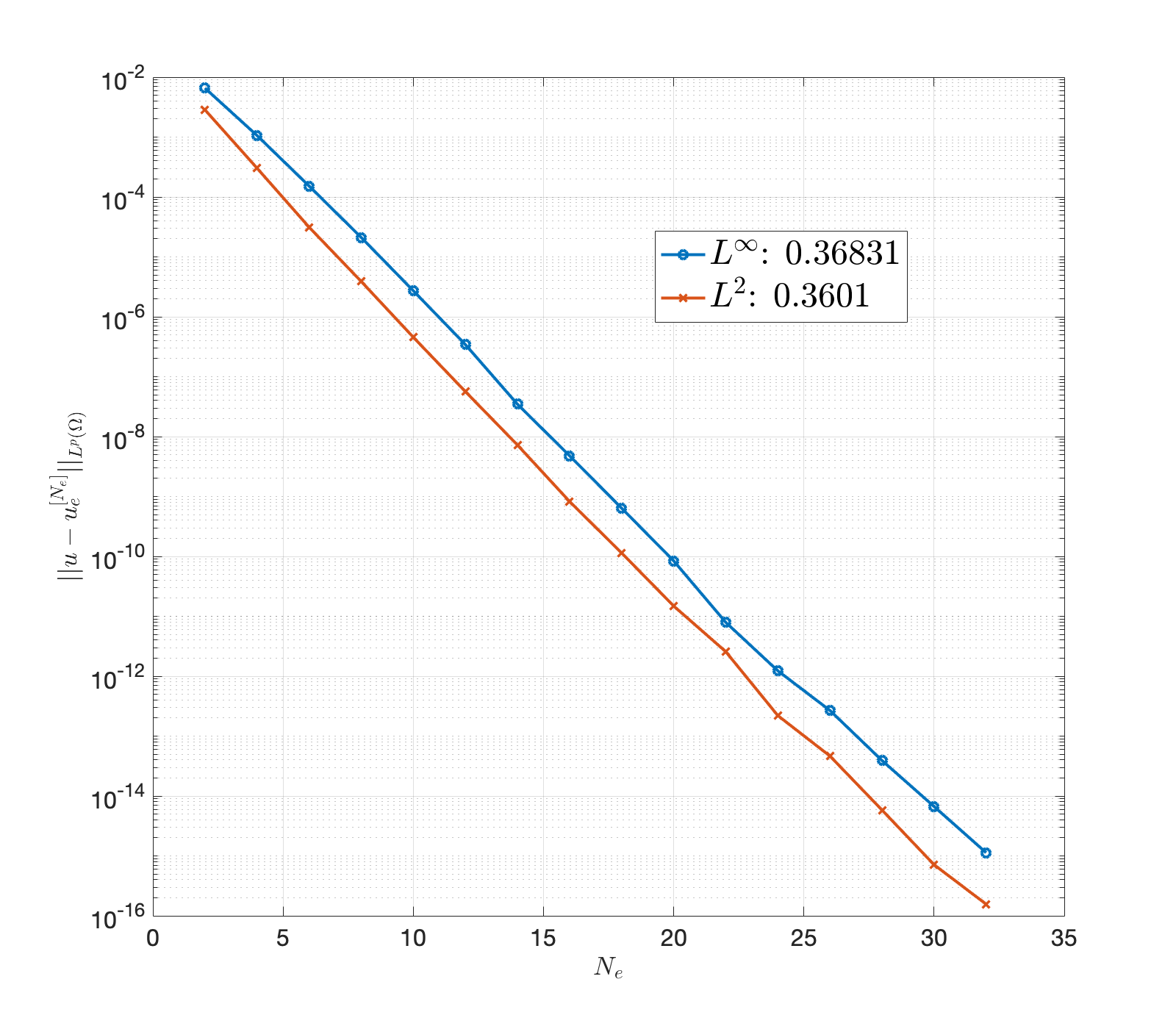}
		\label{Poi2Ddiam}
	}
	\subfigure[Convergence for problem (\ref{Poi2DN})]
	{\includegraphics[scale=\sclbm]{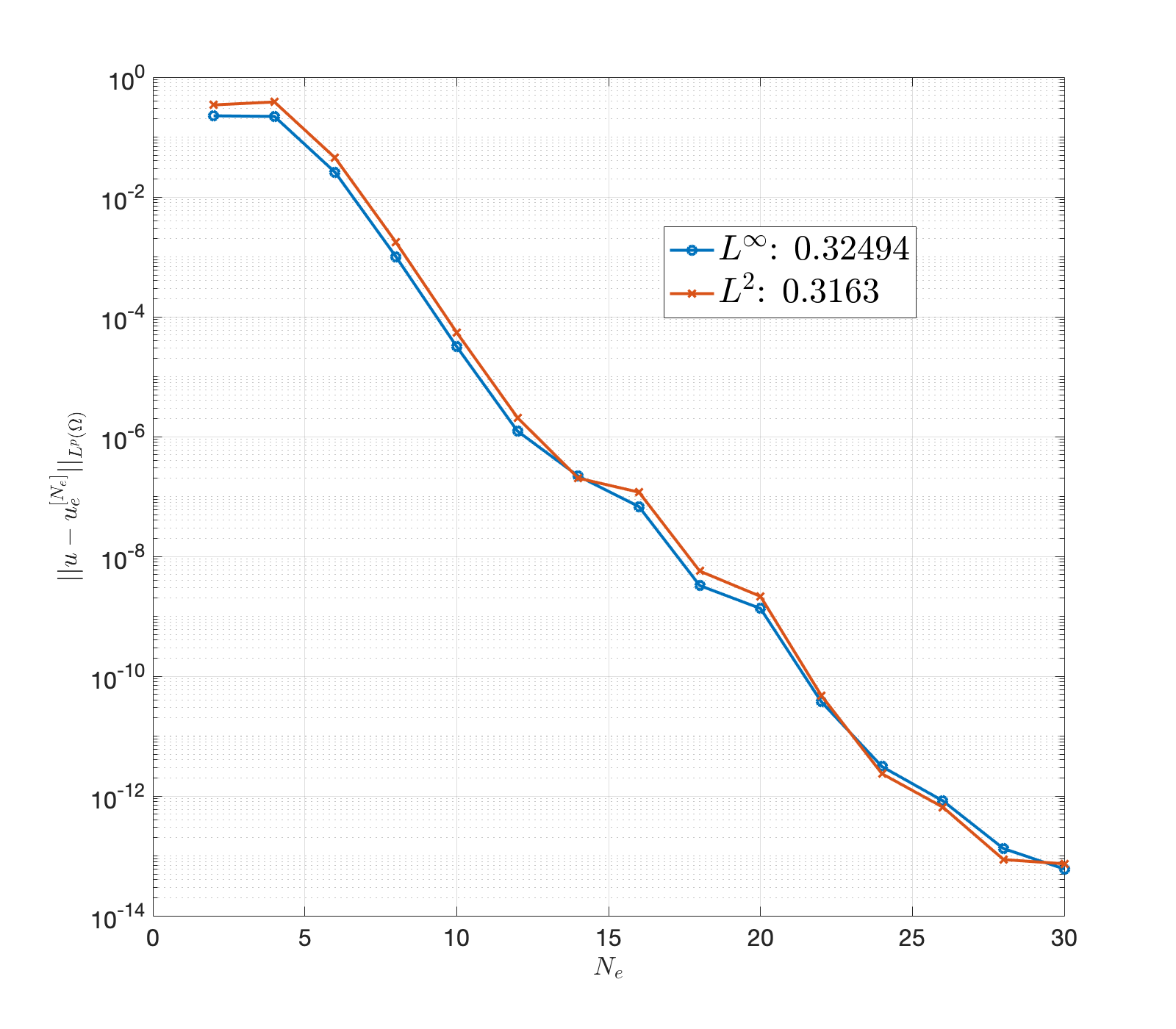}
		\label{Poi2DNeum}
	}
	\caption{Results for the two-dimensional Poisson problems using $N = 2^8$. (a) The domain used in problem (\ref{Poi2D}). (b) The extended solution computed by the algorithm for (\ref{Poi2D}) using $N_\text{e} = 32$. Observe in particular the smooth transition across the boundary $\partial \Omega_1$. (c) Exponential convergence of the form $O(a^{N_\text{e}})$ for the Poisson problem (\ref{Poi2D}) on the domain shown in (a); the estimated values $a$ are shown for the different norms. (d) Exponential convergence for the mixed boundary condition problem (\ref{Poi2DN}).}\label{Poi2Dpl}
\end{figure}

We now move on to two-dimensional test problems. Let $\Omega_1$ be the star-shaped domain shown in Figure \ref{starshape} and consider the Poisson problem  
\eqn{
	\left\{\begin{matrix}
		-\Delta u = -4/(x^2+y^2)^2, & \text{on } \Omega_1, \\
		u = 1/(x^2+y^2), & \text{on } \partial \Omega_1. 	
	\end{matrix}\right.	\label{Poi2D}
} 

This problem has the exact solution $u(x,y) = 1/(x^2+y^2)$. We discretize the boundary $\partial \Omega_1$ by placing $n_\text{b}$ equidistant points along the boundary; as the boundary is a closed loop, the corresponding boundary integral computation (\ref{POip}) via the trapezoid rule is spectrally accurate. We have generally found that choosing $n_\text{b}$ so that the resulting node spacing $\Delta s$ equals the grid spacing $h$ leads to a fairly accurate and robust solver. Figure \ref{Poi2Dudiam} displays an instance of the extended solution computed by the algorithm and Figure \ref{Poi2Ddiam} shows that the convergence is again exponential in $N_\text{e}$ in both norms. Observe that the solution surface is smooth with the transition over the physical boundary barely perceptible. This indicates the effectiveness of our extension routine and an insight into its edge over other extension methods.

Next, let $\Omega_2 = B_2(\pi,\pi)$, i.e., the circular disc of radius 2 centered at $(\pi,\pi)$, and consider
\eqn{
	-\Delta u = -\frac{2}{xy}\left(\frac{1}{x^2} + \frac{1}{y^2}\right), \quad \text{on } \Omega_2, \label{Poi2DN}
} 
with Neumann conditions applied on the upper semicircle of the boundary and Dirichlet on the lower. These constraints are drawn from the exact solution $u(x,y) = 1/(xy)$. The Neumann conditions can be applied by appropriately modifying the objective function (\ref{G2defn}). The convergence results are shown in Figure \ref{Poi2DNeum}.  

Finally, we apply our methodology to the two-dimensional heat equation
\eqn{
	\left\{\begin{matrix}
		u_t-\Delta u = -2\pi\ln(x^2+y^2)\sin(2\pi t), & \text{for } x,y \in \Omega_2, t > 0, \\
		u(0,x,y) = \ln(x^2+y^2), & \text{for } x,y \in  \Omega_2, \\
		u(t,x,y) = \ln(x^2+y^2)\cos(2\pi t), & \text{for } x,y \in  \partial\Omega_2, t > 0. 	\\
	\end{matrix}\right.	\label{Heat2D}
} 

\begin{figure}[tbph]
	\centering
	\subfigure[Convergence of $u_\text{e}^{[N_\text{e}]}$]
	{\includegraphics[scale=\sclbm]{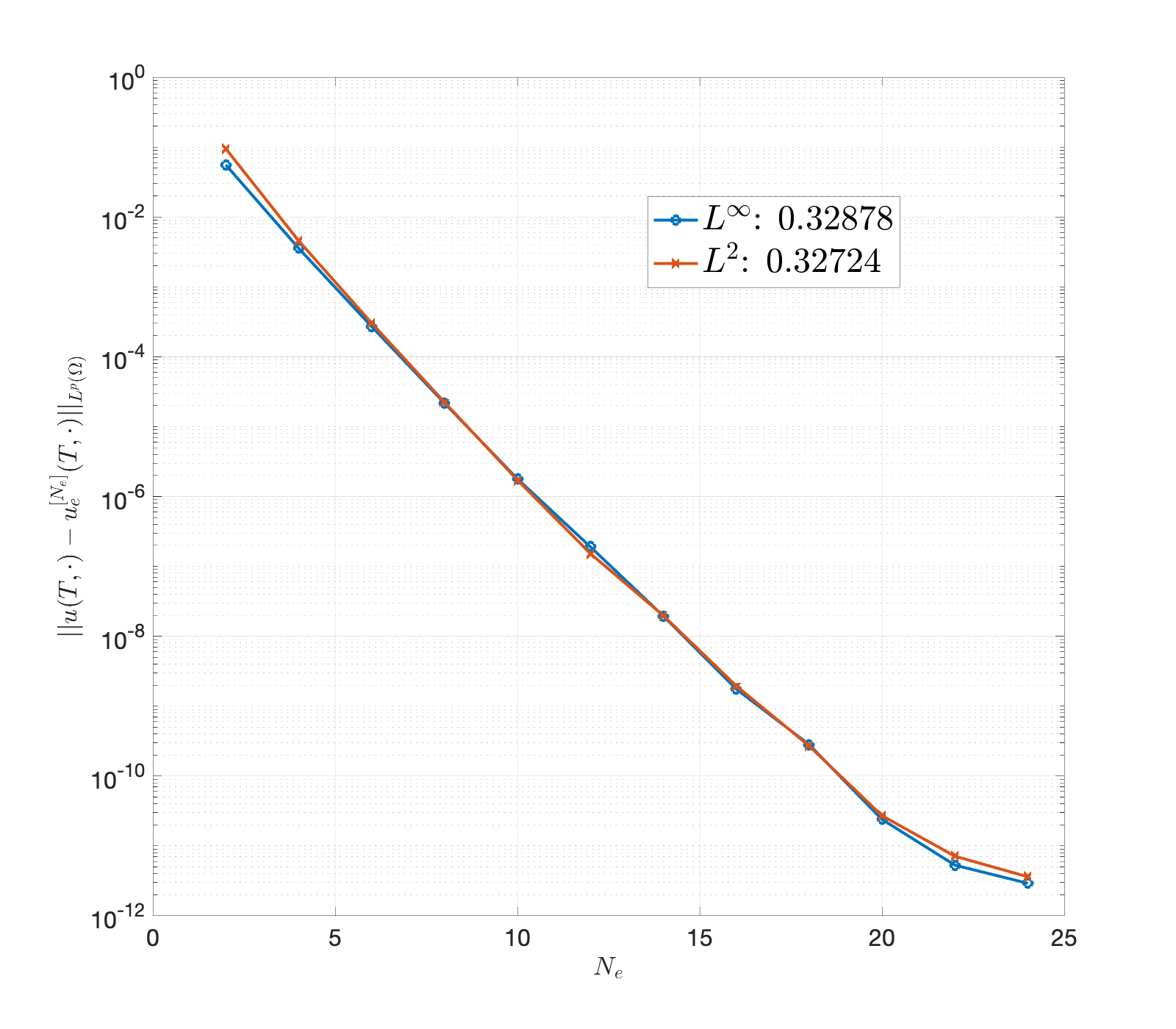}
		\label{heat2Dpl}
	}
	\subfigure[Evolution of errors with $t$ for $N_\text{e} = 32$]
	{\includegraphics[scale=\sclbm]{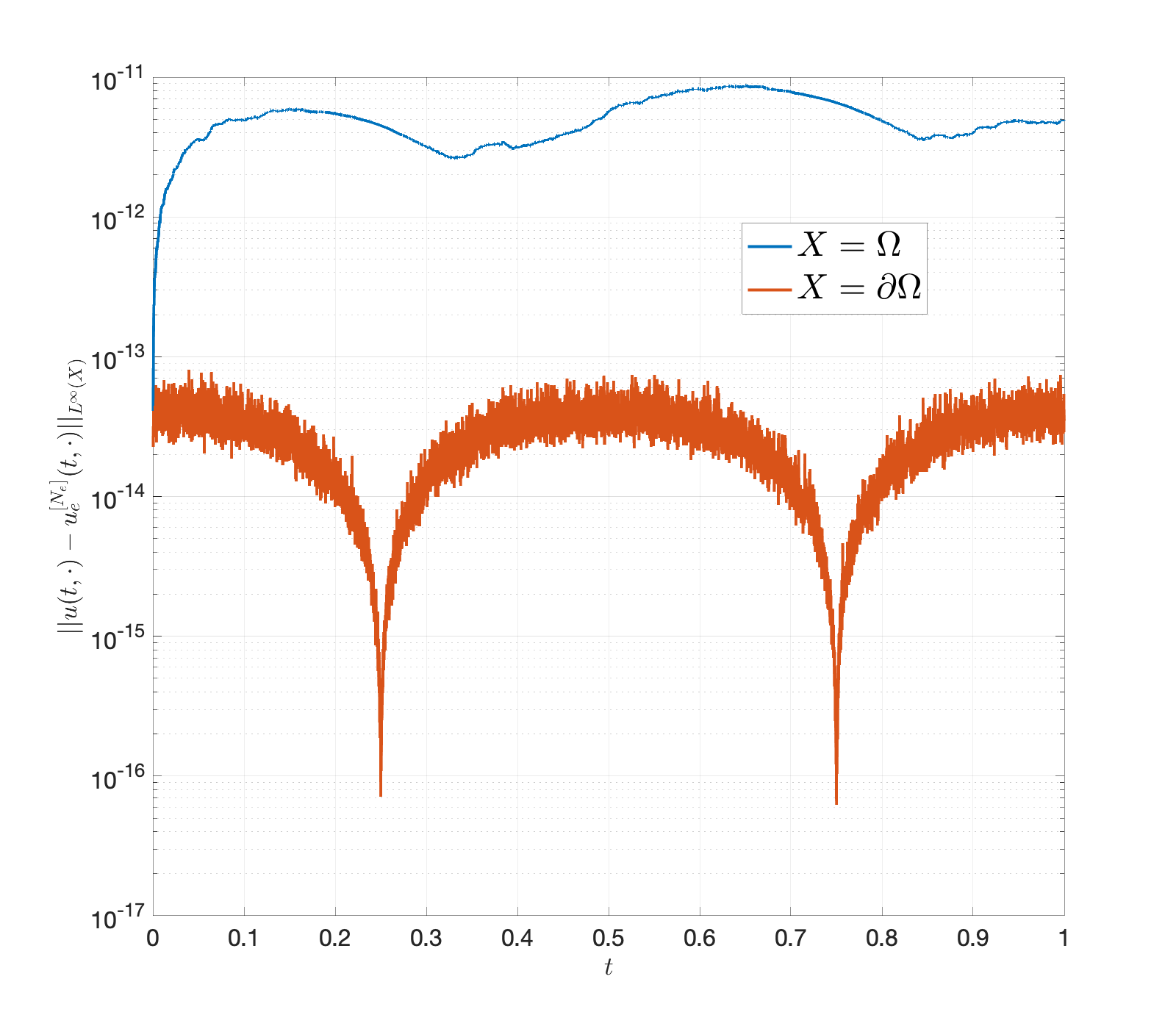}
		\label{heat2DT}
	}
	\caption{Results for the two-dimensional heat equation solved up to $T = 1$. (a) Exponential convergence of the form $O(a^{N_\text{e}})$ with the estimated values $a$ shown for the different norms. (b) Time discretization errors are negligible throughout, as seen here by tracking the $L^{\infty}(\Omega)$ error in the computed solution, as well as the deviation in the boundary values.}\label{heat2D}
\end{figure}

This problem has the exact solution $u(t,x,y) = \ln(x^2+y^2)\cos(2\pi t)$. We solve it up to $T = 1$, with time-step $\Delta t = 10^{-4}$ for the BDF-4 scheme as before, and choose $N = 2^8$ grid-points in each dimension. The resulting convergence can be again seen to be geometric with $N_\text{e}$ in Figure \ref{heat2Dpl}.

\section{Newtonian Fluid Problems}\label{FluidProbs}
In this section, we provide the details of our method applied to some of the more common Newtonian fluid models. 

\subsection{The General Stokes Problem}\label{stokessub}
Let $\mathcal{L}$ and $\Omega$ be as in the Section \ref{mathform}. We begin by considering the general Stokes system
\eqn{
	\left\{ \begin{matrix}
		\mathcal{L} \tb{u} + \nabla p = \tb{f}, &  \text{on } \Omega,  \\
		\nabla \cdot \tb{u} = 0, &  \text{on } \Omega,  \\
		\tb{u} = \tb{g}, &   \text{on } \partial \Omega. \\ 
	\end{matrix} \right. 
	\label{stokes}
} 

Following the spirit of our earlier exposition, we seek an extended forcing $\tb{f}_\text{e}$ over a computationally simpler domain $\Pi$. Note that the invertibility of the Stokes system on $\Pi$ depends on the choice of $\mathcal{L}$. As in Subsection \ref{NinvL}, we assume that $\mathcal{L}$ is self-adjoint on $\Pi$, set $\mathcal{R} = \text{range}(\mathcal{L})$, and let $\mathcal{A} = \mathcal{L}|_\mathcal{R}$. In addition, denote by $\{\Psi^{(k)}\}_{k \in \mathcal{K}}$ a divergence-free orthonormal basis for the null space of $\mathcal{L}$ on $\Pi$. If $\tb{u}_\text{e}$ is the solution to the extended problem on $\Pi$, then we must have, for all $k \in \mathcal{K}$,
\eqn{
	\ip{\Psi^{(k)},\tb{f}_\text{e}}_\Pi = \ip{\Psi^{(k)},\mathcal{L}\tb{u}_\text{e} + \nabla p_\text{e}}_\Pi = \ip{\mathcal{L}\Psi^{(k)},\tb{u}_\text{e}}_\Pi - \ip{\nabla \cdot \Psi^{(k)} , p_\text{e}}_\Pi  = 0. \label{fepsikst}
}

In addition, we set $d_k = \ip{\Psi^{(k)},\tb{u}_\text{e}}_\Pi$ and define
\eqn{
	\tb{u}_0 = \tb{u}_\text{e} - \sum_{k \in \mathcal{K}}d_k\Psi^{(k)}. \label{u0stdef}
}

This obeys the ``range-restricted'' Stokes system:
\eqn{
	\left\{ \begin{matrix}
		\mathcal{A} \tb{u}_0 + \nabla p_0 = \tb{f}_\text{e}, &  \text{on } \Pi,  \\
		\nabla \cdot \tb{u}_0 = 0, &  \text{on } \Pi.  \\
	\end{matrix} \right. 
	\label{stokesR}
}

Denote by $\tb{u}^{(l,j)}$ the solution to the range-restricted Stokes system on $\Pi$ with forcing $\Phi^{(l,j)}$ whose components ae $\Phi^{(l,j)}_i = \delta_{il}\phi_j$, for $1 \leq i,l \leq d$ and $j \in \mathcal{J}$, i.e., 
\eqn{
	\left\{ \begin{matrix}
		\mathcal{A} \tb{u}^{(l,j)} + \nabla p^{(l,j)} = \Phi^{(l,j)} &  \text{on } \Pi,  \\
		\nabla \cdot \tb{u}^{(l,j)}  = 0, &  \text{on } \Pi.  \\
	\end{matrix} \right. 
	\label{stokesforc}
}

Let $\Theta^{(l,j)} = S^*\tb{u}^{(l,j)} $ be the restriction of this solution to the boundary $\partial \Omega$, and similarly set $\gamma^{(k)} = S^*\Psi^{(k)}$. Next, expand the extended forcing as 
\eqn{
	\tb{f}_\text{e} = \sum_{j \in \mathcal{J}} \sum_{l = 1}^d c_j^{(l)}\Phi^{(l,j)} = \sum_{j \in \mathcal{J}} \begin{pmatrix} c_j^{(1)} \\ c_j^{(2)} \\ \vdots \\ c_j^{(d)}  \end{pmatrix} \phi_j. \label{forcext}
}

The boundary and forcing agreement conditions are then imposed by minimizing 
\eqn{
	G(\tb{c},\tb{d})  &=& \norm{\sum_{j\in \mathcal{J}} \sum_{l = 1}^d c^{(l)}_j\Theta^{(l,j)} + \sum_{k\in \mathcal{K}} d_k\gamma^{(k)} - \tb{g}}^2_{L^2(\partial\Omega)} + \nonumber\\
	&& \norm{\sum_{j\in \mathcal{J}} \sum_{l = 1}^d c^{(l)}_j\Phi^{(l,j)} - \tb{f}}^2_{L^2(\Omega)}  + \sum_{k \in \mathcal{K}} \left|\sum_{j \in \mathcal{J}} \sum_{l = 1}^d c_j^{(l)}\ip{\Psi^{(k)},\Phi^{(l,j)}}_\Pi\right|^2. \label{GStdefn}
}

Setting $\partial G/\partial c_{j_0}^{(l_0)} = 0$ then yields
\eqn{
	&& \sum_{j\in \mathcal{J}} \sum_{l = 1}^d \left(\ip{\Theta^{(l_0,j_0)},\Theta^{(l,j)}}_{\partial \Omega} + \ip{\Phi^{(l_0,j_0)},\Phi^{(l,j)}}_{\Omega} + \sum_{k \in \mathcal{K}} \ip{\Psi^{(k)},\Phi^{(l_0,j_0)}}_\Pi \ip{\Psi^{(k)},\Phi^{(l,j)}}_\Pi \right)c_j^{(l)} + \nonumber\\
	&& \hspace{1 cm}  \sum_{k \in \mathcal{K}}  \ip{\Theta^{(l_0,j_0)},\gamma^{(k)}}_{\partial \Omega} d_k \hspace{0.5 cm} = \hspace{0.5 cm}  \ip{\Theta^{(l_0,j_0)},\tb{g}}_{\partial \Omega} + \ip{\Phi^{(l_0,j_0)},\tb{f}}_{\Omega}, \label{sqsys3a}
}
while $\partial G/\partial d_{k_0} = 0$ gives
\eqn{
	\sum_{j\in \mathcal{J}} \sum_{l = 1}^d \ip{\gamma^{(k_0)} , \Theta^{(l,j)}}_{\partial \Omega}  c^{(l)}_j + \sum_{k\in \mathcal{K}} \ip{\gamma^{(k_0)} , \gamma^{(k)}}_{\partial \Omega}  d_k = \ip{\gamma^{(k_0)} , \tb{g}}_{\partial \Omega}. \label{sqsys3b}
}

Equations (\ref{sqsys3a}) and (\ref{sqsys3b}) again form a self-adjoint linear system. Upon truncating the basis $\{\phi_j\}_{j \in \mathcal{J}}$, it reduces to a Hermitian linear problem which can be further simplified by discretizing the inner products and solved by a QR decomposition, as detailed in Subsection \ref{NumInt}.

\begin{figure}[tbph]
	\centering
	\subfigure[]
	{\includegraphics[scale=\sclbm]{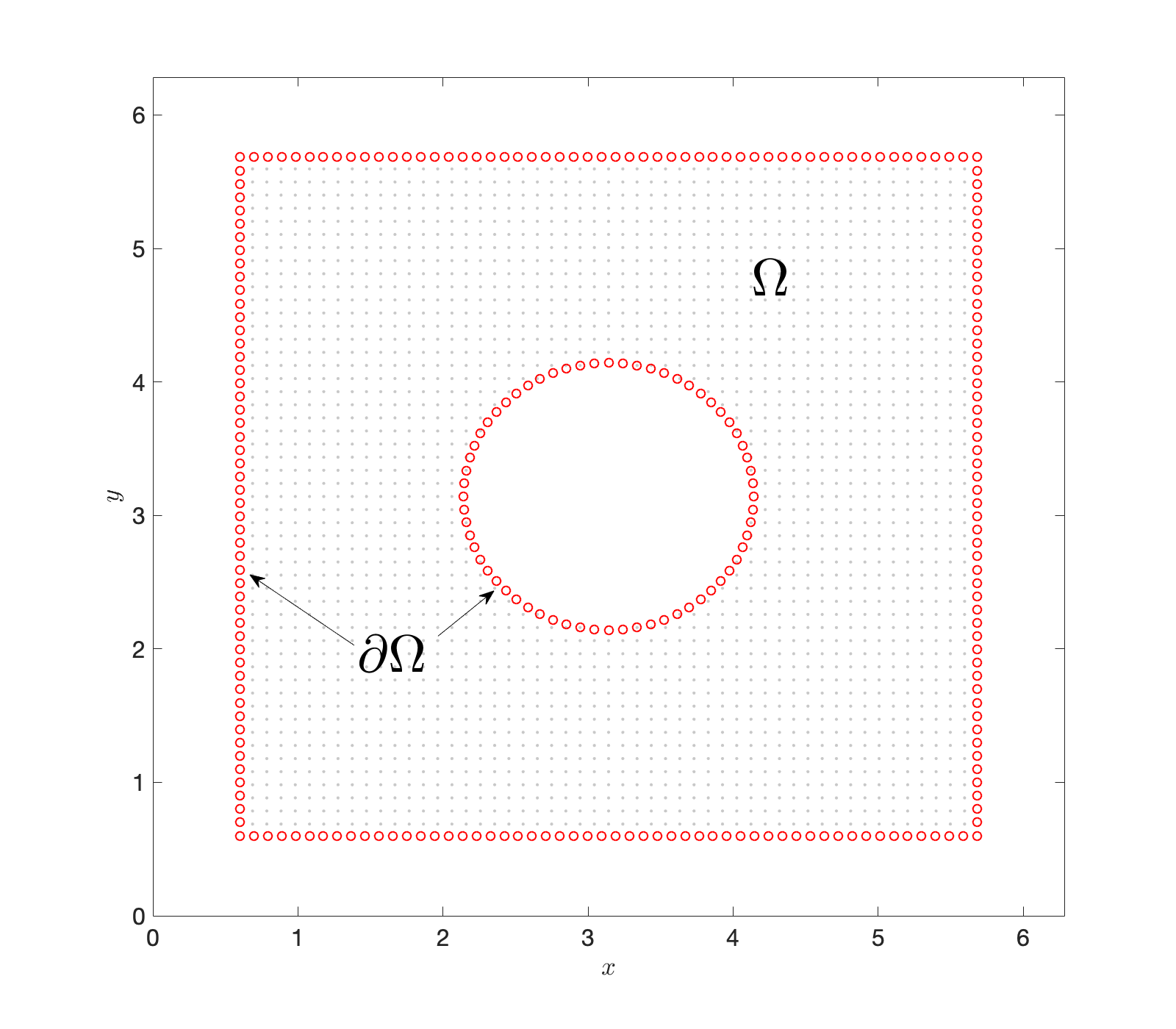}
		\label{rectdom}
	}
	\subfigure[]
	{\includegraphics[scale=\sclbm]{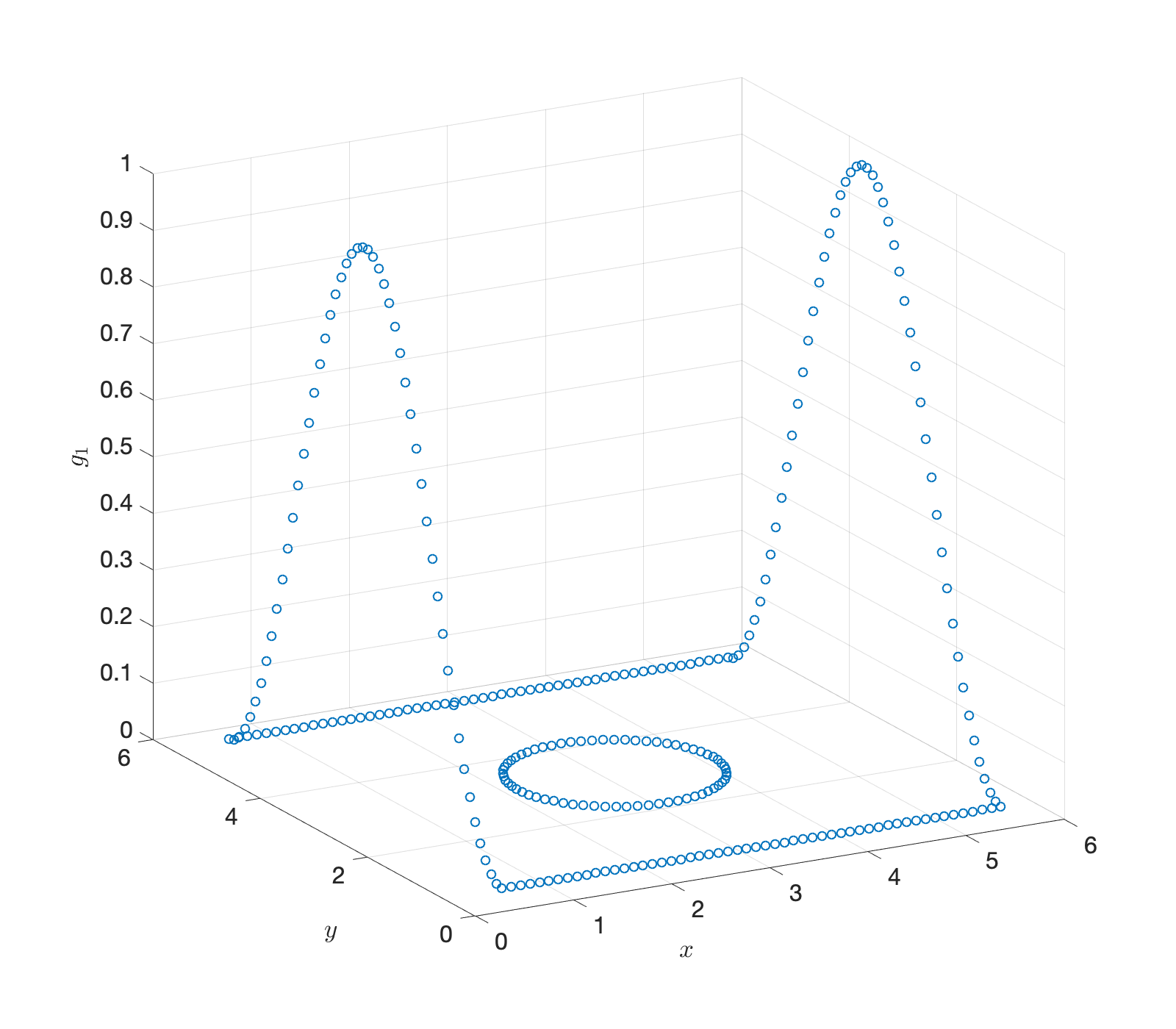}
		\label{flowBC}
	}
	\caption{(a) The domain $\Omega$ used for the fluid test problems. The walls are at $x = 0.6$, $2\pi-0.6$ and $y = 0.6$, $2\pi-0.6$, and the unit disc is centered at $(\pi,\pi)$. (b) The $x$-component of the boundary condition for test problem (\ref{stokesEx1}). An identical inflow and outflow condition is prescribed, and no-slip is enforced everywhere else.}\label{architecture}
\end{figure}

We test our method by solving three Stokes problems with $\mathcal{L} = -\Delta$ on the domain $\Omega = [0.6,2\pi-0.6]^2 - \overline{B_1(\pi,\pi)}$ (see Figure \ref{rectdom}). We again take $\Pi$ to be the two-dimensional torus $\mathbb{T}^2$ and use the Fourier basis as the extension functions. As our first example, we consider the problem obtained from the exact solution
\eqn{
	\tb{u}(x,y) =  \begin{pmatrix}
		e^{\sin(x)}\cos(y) \\
		-e^{\sin(x)}\sin(y)\cos(x)
	\end{pmatrix}, \qquad p(x,y) = e^{2\cos(x)}. \label{stokesEx0}
} 

\def \sclcm {0.1}
\begin{figure}[tbph]
	\centering
	\subfigure[Problem (\ref{stokesEx0})]
	{\includegraphics[scale=\sclcm]{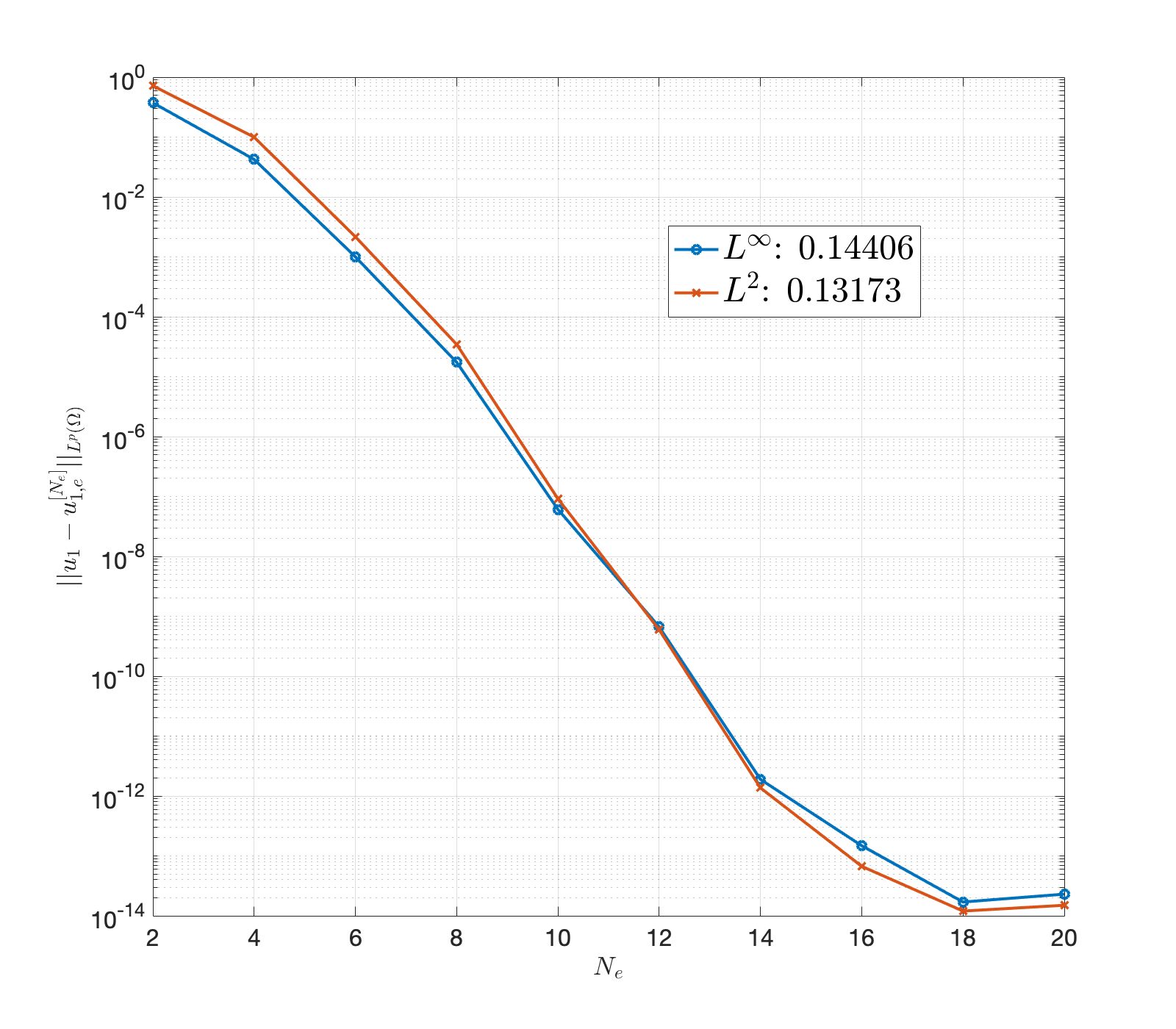}
		\includegraphics[scale=\sclcm]{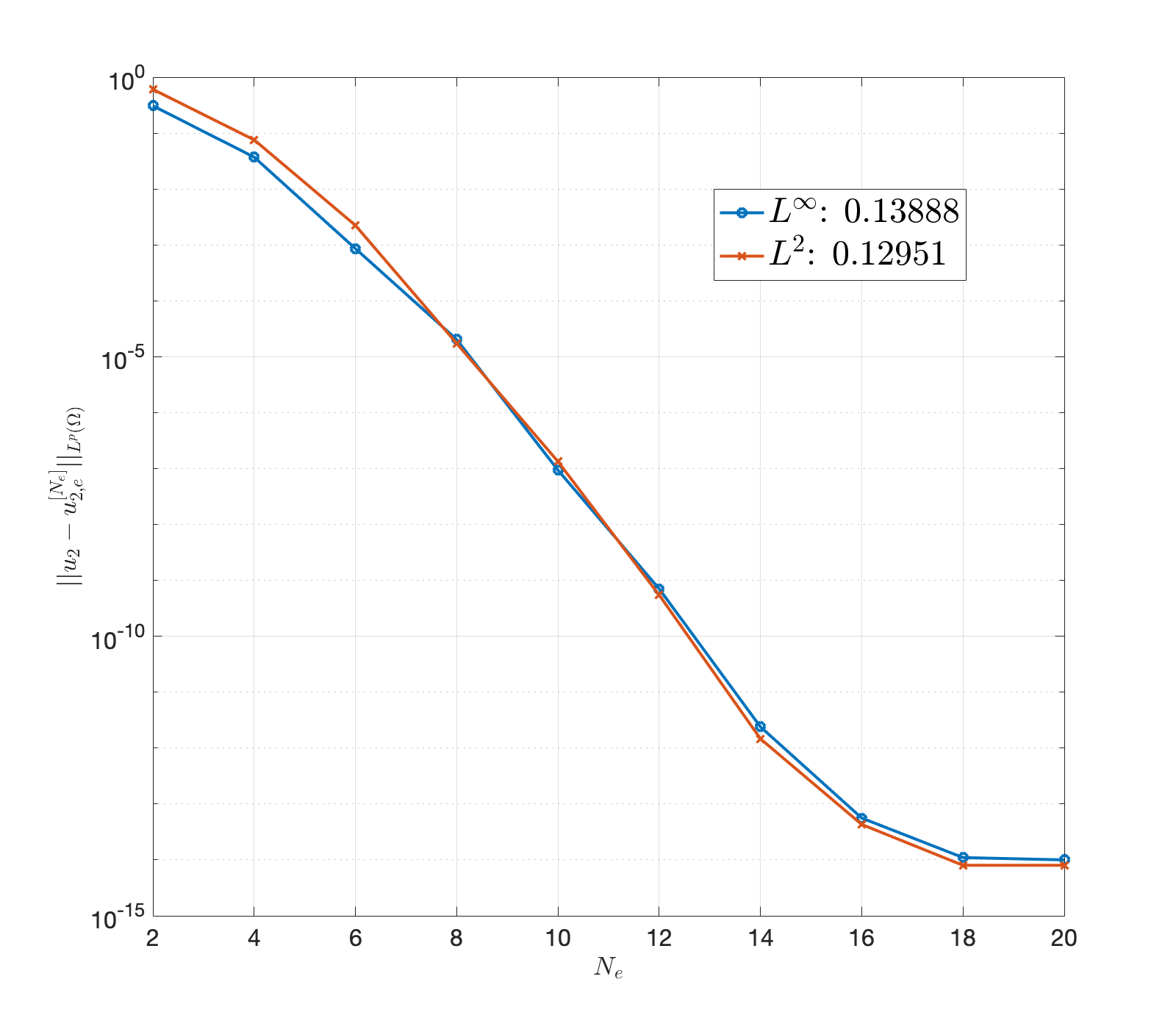}
		\includegraphics[scale=\sclcm]{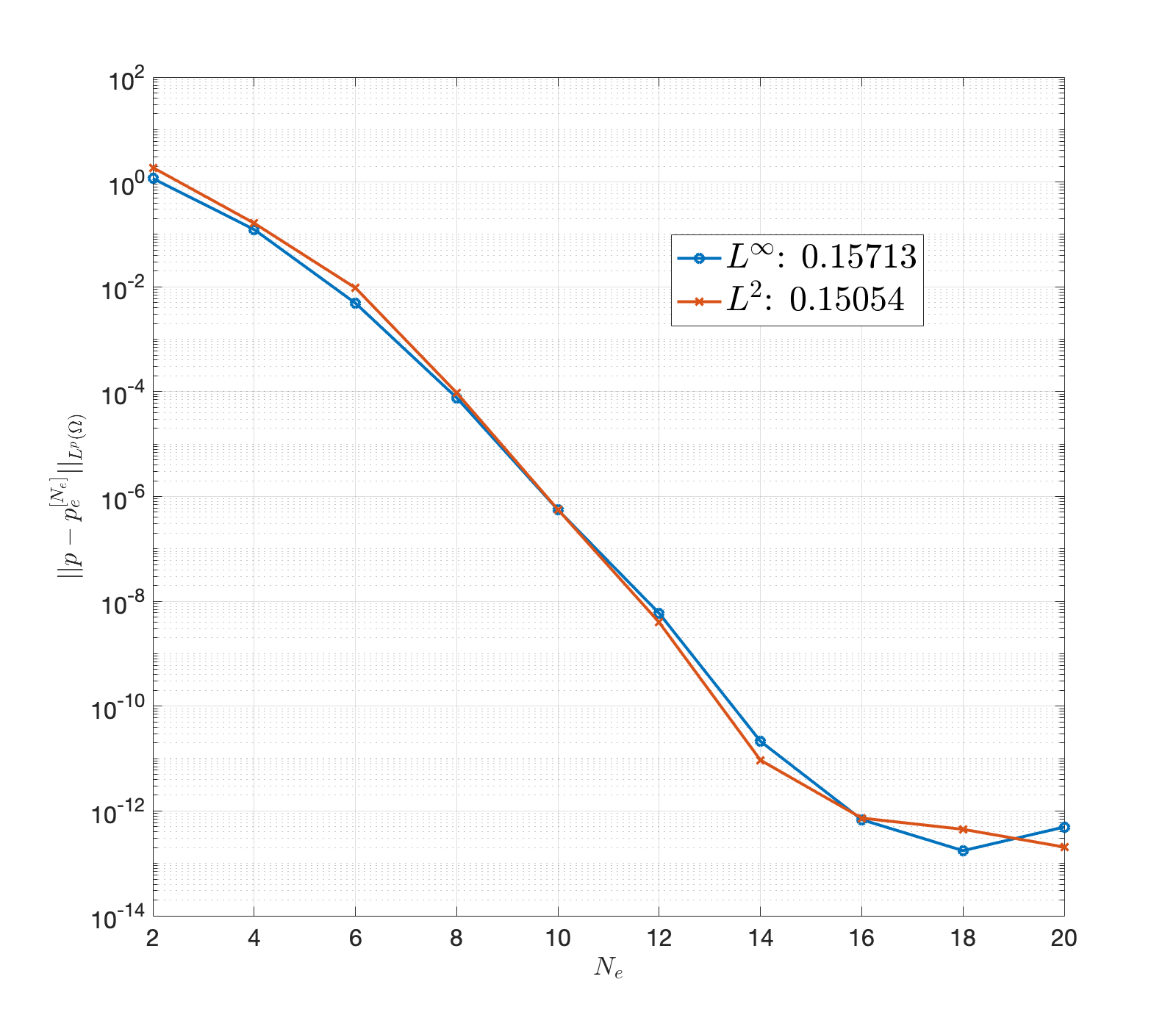}
		\label{Stex0}
	}
	\subfigure[Problem (\ref{stokesEx2})]
	{\includegraphics[scale=\sclcm]{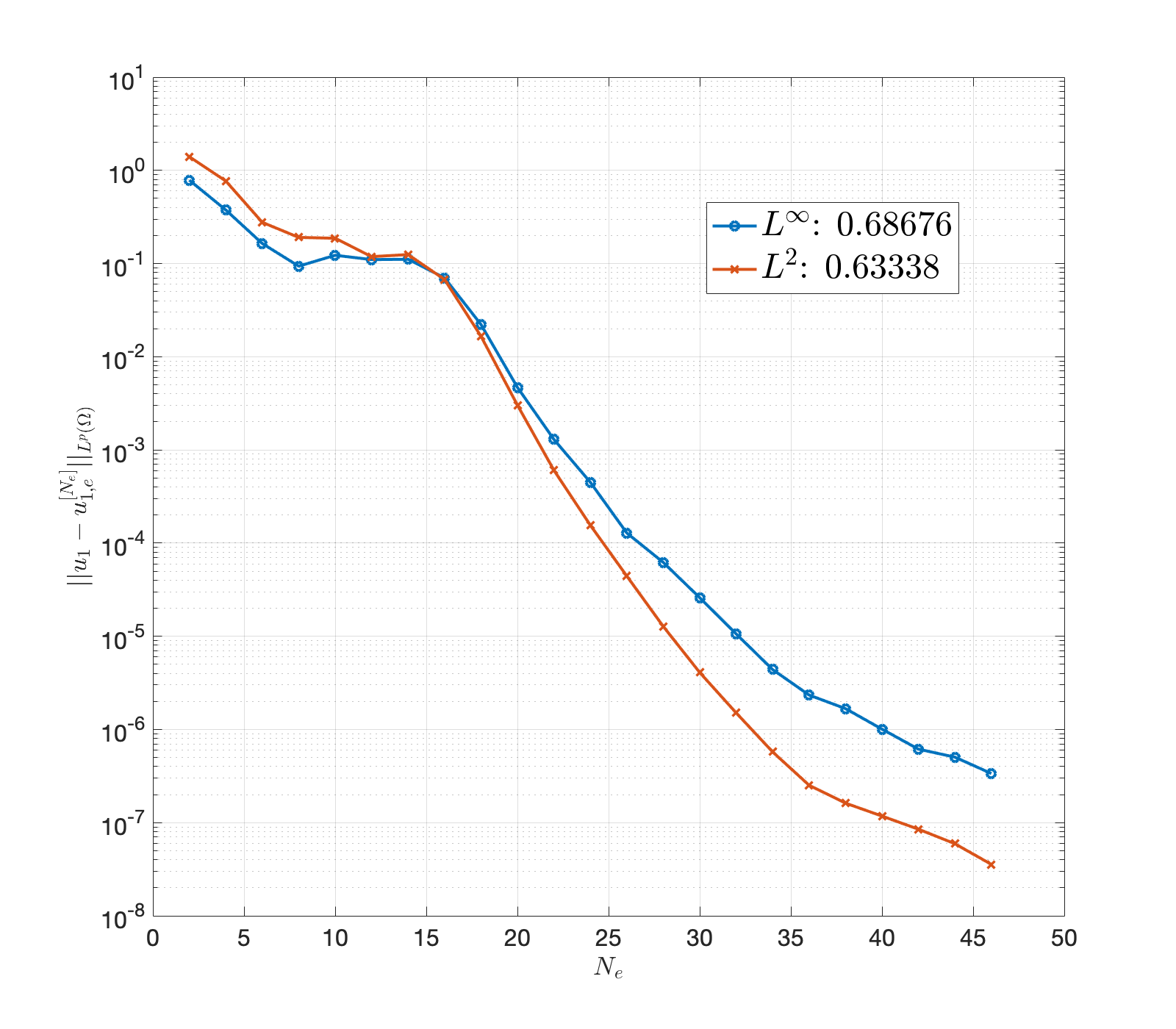}
		\includegraphics[scale=\sclcm]{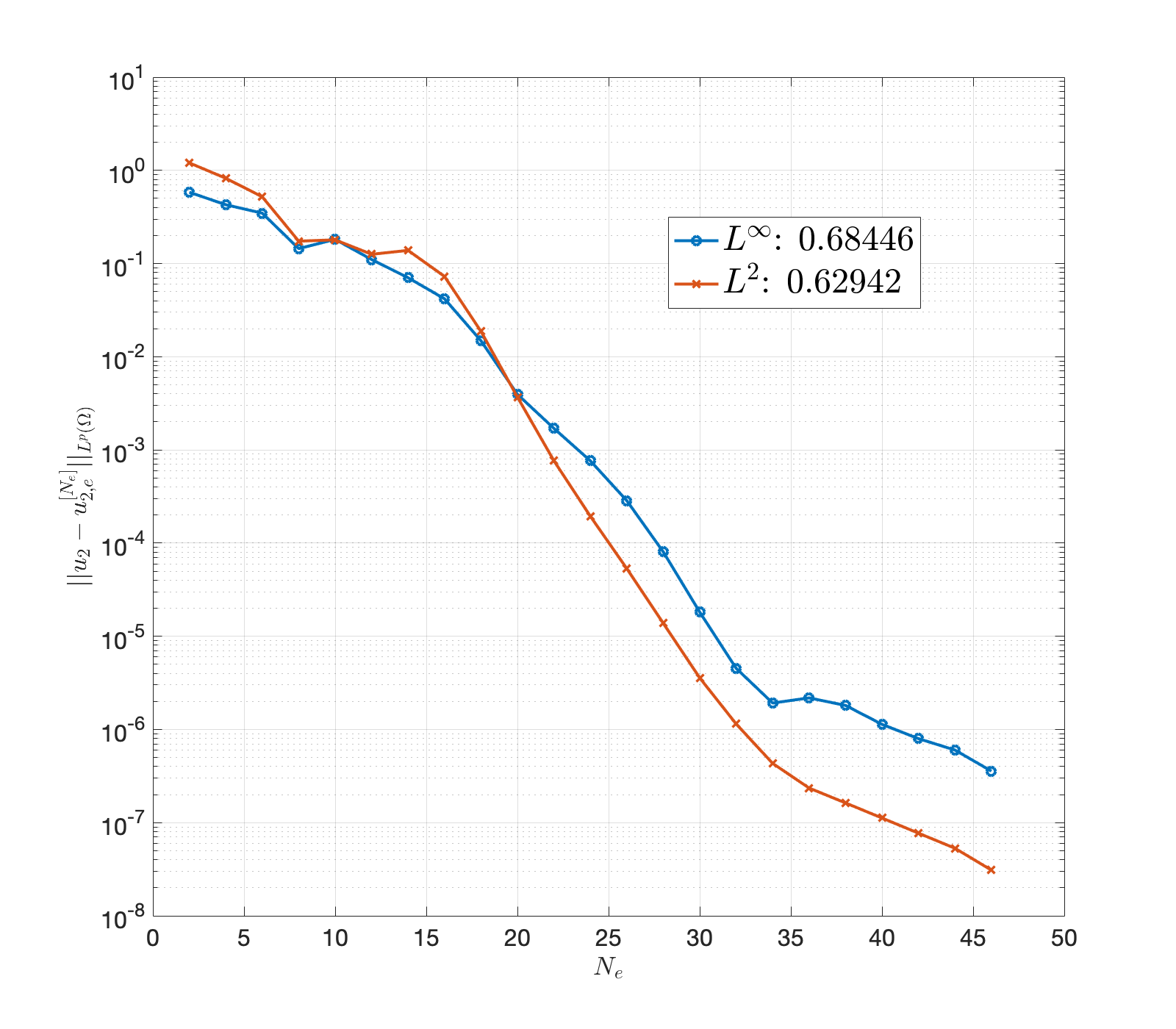}
		\includegraphics[scale=\sclcm]{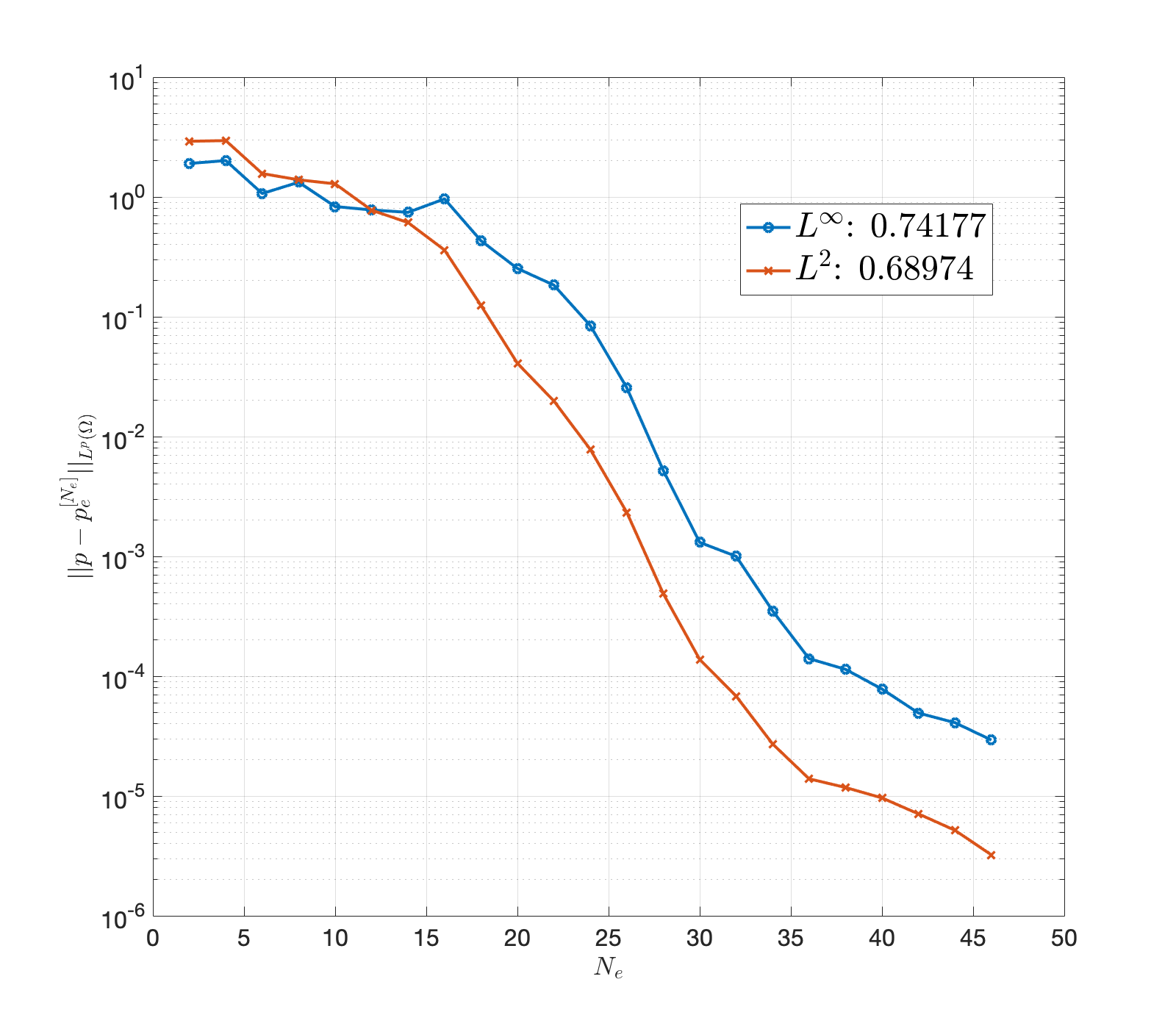}
		\label{Stex2}
	}
	\subfigure[Problem (\ref{stokesEx1})]
	{\includegraphics[scale=\sclcm]{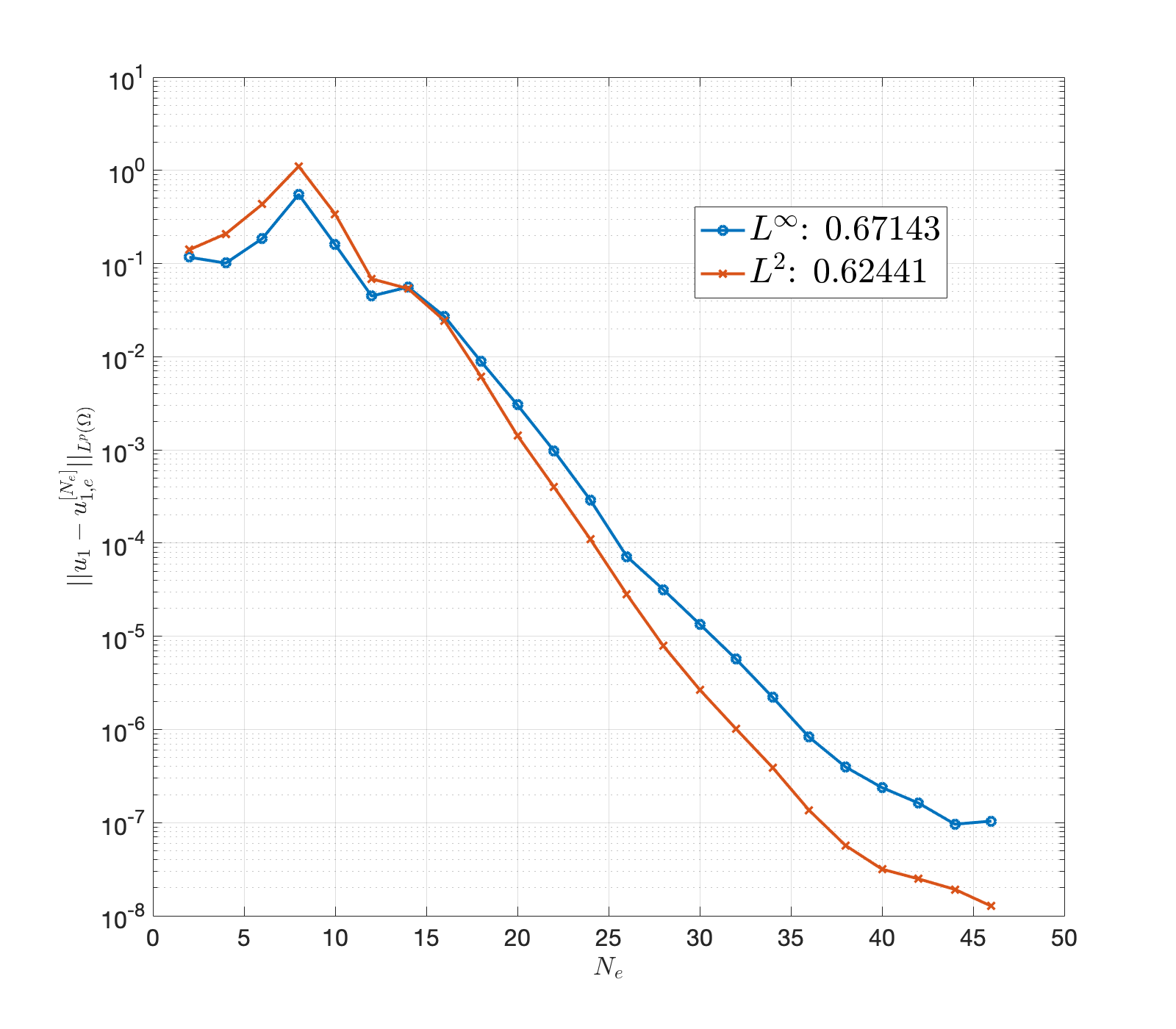}
		\includegraphics[scale=\sclcm]{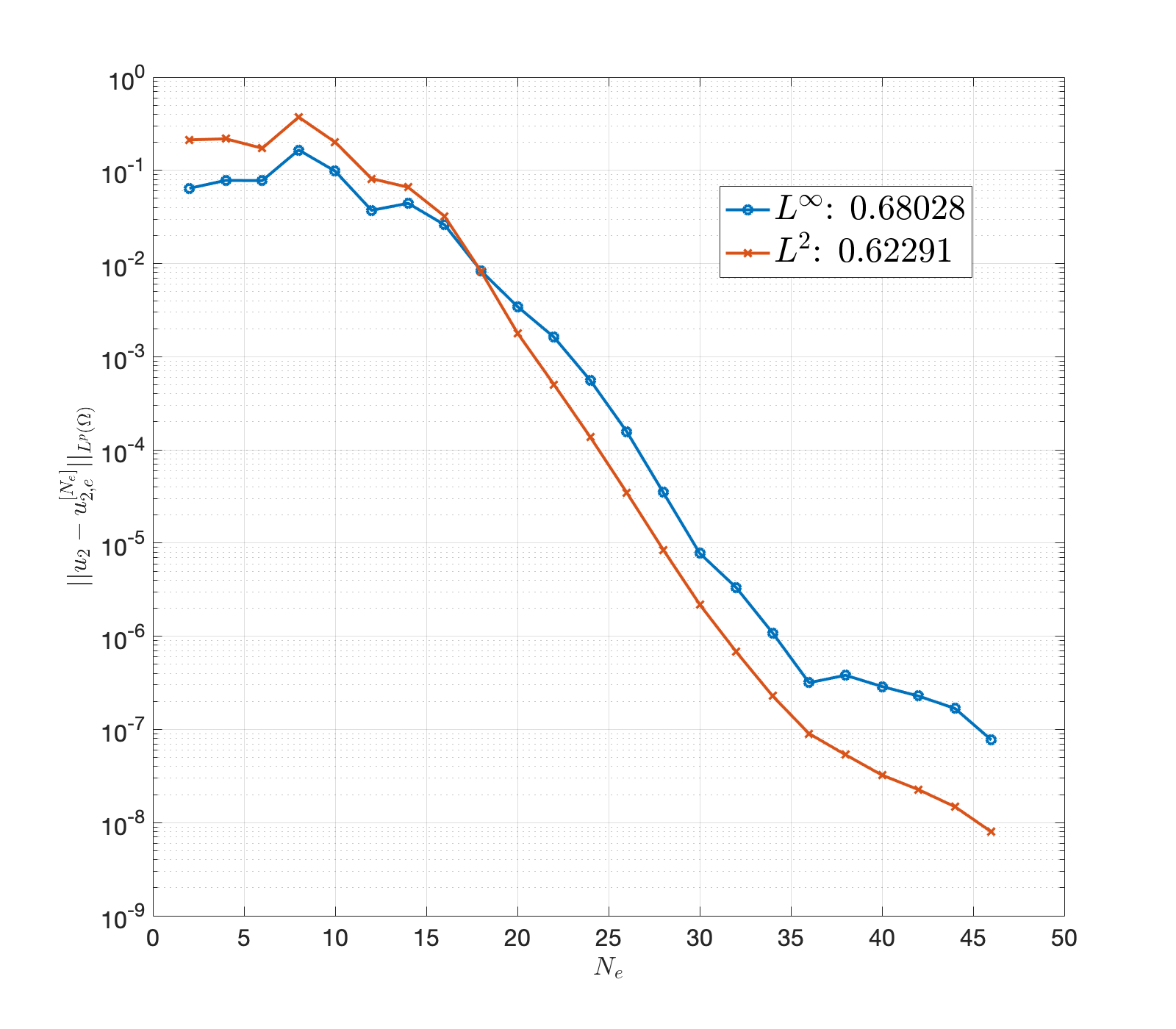}
		\includegraphics[scale=\sclcm]{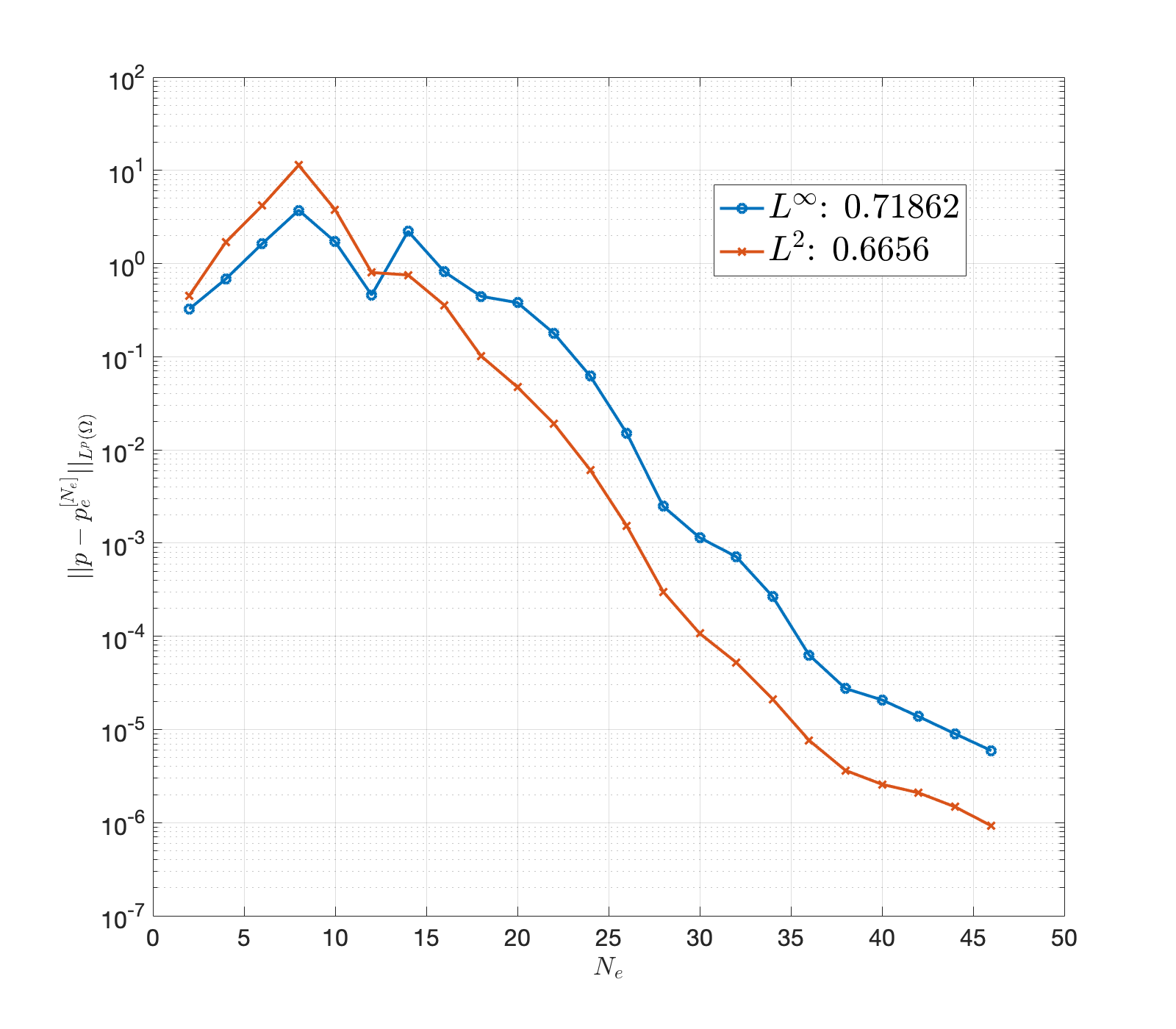}
		\label{Stex1}
	}
	\caption{Convergence plots for the solutions $u_1$, $u_2$ and $p$ to Stokes problems on the domain shown in Figure \ref{rectdom}. All unknowns appear to converge exponentially at $O(a^{N_\text{e}})$; the corresponding best-fit values of $a$ are shown for the different norms.
	}\label{stokesconv}
\end{figure}

Next, we apply our technique to 
\eqn{
	\left\{ \begin{matrix}
		-\Delta \tb{u} + \nabla p = \tb{f}
		, &  \text{on } \Omega,  \\
		\nabla \cdot \tb{u} = 0, &  \text{on } \Omega,  \\
		\tb{u} = \tb{0}, &   \text{on } \partial \Omega, \\ 
	\end{matrix} \right. 
	\label{stokesEx2}
}
with
\eqn{
	\tb{f}(x,y) =  \begin{pmatrix}
		e^{\sin(x)}\cos(y)(1+\sin(x)-\cos^2(x)) - 2\sin(x)e^{2\cos(x)} \\
		e^{\sin(x)}\sin(y)\cos(x)(\cos^2(x) - 3\sin(x) - 2)
	\end{pmatrix}, \label{stokesf}
}

This is the same forcing as in problem (\ref{stokesEx0}), only this time we combine it with homogeneous Dirichlet boundary conditions. Finally, we solve the problem with zero forcing and boundary conditions
\eqn{
	\tb{g}(x,y) = \left\{
	\begin{matrix}
		\begin{pmatrix} \sin^2\left(\frac{\pi (y-0.6)}{2\pi-1.2}\right) \\ 0
		\end{pmatrix} & x = 0.6, x = 2\pi-0.6, 0.6 \leq y \leq 2\pi-0.6, \\
		\tb{0} & \text{otherwise}
	\end{matrix} \right. \label{stokesEx1}
}

This choice is meant to replicate inflow and outflow with no-slip enforced everywhere else (see Figure \ref{flowBC}). 

The operator $\mathcal{L} = -\Delta$ on $\Pi$ has a two-dimensional null-space, spanned by $\left\{ \begin{pmatrix}
1 \\ 0
\end{pmatrix} , \begin{pmatrix}
0 \\ 1
\end{pmatrix} \right\}$, so we need to factor this into our system of equations. We set $N = 2^8$ and place boundary points roughly $\Delta s = 2\pi/N$ apart. Recall that our technique makes no assumptions about the connectedness of the domain or smoothness of its boundary; as a result, it does not require domain-specific adjustments. 

The corresponding convergence plots are shown in Figure \ref{stokesconv}. The solutions to problem (\ref{stokesEx0}) converge geometrically almost instantly and continue rapidly until machine precision is achieved (at around $N_\text{e} = 18$). The two problems with unknown solutions appear to require a minimum number of extension functions, after which the convergence is exponential in velocity and pressure values. This suggests that our technique is well-equipped to handle complicated geometries as well as arbitrary forcing and boundary conditions.

\subsection{The Alternative Formulation for Stokes Equation}\label{AltStokes}
While the approach developed out in Subsection \ref{stokessub} is fairly general and flexible, it is advantageous to explore the use of other computational domains and, by extension, other bases. These basis functions may not be as straightforward to invert as a Fourier basis but offer other advantages with respect to the problem and geometry. For example, a common benchmark problem is flow in a channel past an obstacle. One way to model this is to impose periodicity at the inlet-outlet walls and enforce the no-slip condition at the lateral walls, somewhat similar to Figure \ref{flowBC}. In this case, it is helpful to consider an alternate computational domain and basis so that the inherent periodicity and no-slip constraints are satisfied by default. This can be achieved, e.g.,  by combining the Fourier basis in the $x$-direction and a polynomial basis in the $y$-direction. Accordingly, this calls for a methodology built upon the alternative formulation in Subsection \ref{AltForm}. For simplicity, we only provide the details for a two-dimensional problem. 

We begin by partitioning the boundary $\partial \Omega = \partial \Omega_\text{I} \cup \partial \Omega_\text{A}$ into intrinsic and auxiliary parts respectively. We shall enforce the no-slip boundary conditions on $\partial \Omega_\text{I}$ by default, and those on $\partial\Omega_\text{A}$ by the minimization procedure. Denote the latter boundary condition by $\tb{g}_\text{A}$. In contrast with the earlier approach, where we only searched for an appropriate projection of the forcing, we now need to project the velocity and pressure terms. The velocity basis functions must be chosen so they obey both the intrinsic boundary conditions and the incompressibility constraint. We enforce these by choosing the basis functions to be the anti-symmetric derivatives of a family of appropriate scalar functions $\{\phi_j\}_{\mathcal{J}}$ on $\Pi$. More precisely, we set
\eqn{
	\Phi_j = \begin{pmatrix}
		\py \phi_j \\ -\px \phi_j
	\end{pmatrix}, \label{AsymDer}
}
where the $\{\phi_j\}_{\mathcal{J}}$ are chosen in the first place to ensure that the velocity basis functions $\{\Phi_j\}_{j \in \mathcal{J}}$ meet the intrinsic constraints. No such provisions are necessary for the pressure term expansion so we simply choose a complete basis $\{\eta_k\}_{k \in \mathcal{K}}$. Write
\eqn{
	\tb{u}_\text{e} = \sum_{j \in \mathcal{J}} c_j \Phi_j, \qquad p_\text{e} = \sum_{k \in \mathcal{K}} d_k \eta_k, \label{uepedef} 
}
and define $\Upsilon_j := S^*_\text{A} \Phi_j$ as the restriction of the velocity basis functions to the auxiliary boundary. We then seek to minimize
\eqn{
	G(\tb{c},\tb{d}) = \norm{\sum_{j \in \mathcal{J}} c_j \Upsilon_j - \tb{g}_\text{A}}^2_{L^2(\partial \Omega_\text{A})} + \norm{\sum_{j \in \mathcal{J}} c_j \mathcal{L} \Phi_j+ \sum_{k \in \mathcal{K}} d_k \nabla  \eta_k - \tb{f}}^2_{L^2(\Omega)}. \label{GSt2defn}
}

Setting $\partial G/\partial c_{j_0} = 0$ for all $j_0 \in \mathcal{J}$ yields
\eqn{
	\sum_{j \in \mathcal{J}} \left(\ip{\Upsilon_{j_0},\Upsilon_j}_{\partial \Omega_\text{A}}  + \ip{\mathcal{L} \Phi_{j_0}, \mathcal{L} \Phi_{j}}_{\Omega} \right)c_j + \sum_{k \in \mathcal{K}} \ip{\mathcal{L} \Phi_{j_0} , \nabla \eta_k}_{\Omega} d_k = \ip{\Upsilon_{j_0},\tb{g}_\text{A}}_{\partial \Omega_\text{A}} + \ip{\mathcal{L} \Phi_{j_0}, \tb{f}}_{\Omega}, \label{sqsys4a}  
}
while $\partial G/\partial d_{k_0} = 0$ for all $k_0 \in \mathcal{K}$ leads to
\eqn{
	\sum_{j \in \mathcal{J}} \ip{\nabla \eta_{k_0},\mathcal{L} \Phi_j}_{\Omega}c_j + \sum_{k \in \mathcal{K}} \ip{\nabla \eta_{k_0} , \nabla \eta_k}_{\Omega}d_k = \ip{\nabla \eta_{k_0}, \tb{f}}_{\Omega}. \label{sqsys4b}  
}

This self-adjoint problem can be discretized to yield a Hermitian system, which can in turn be solved by the techniques detailed earlier. 

To illustrate this approach, we solve the problem of Stokes flow in a channel past a circular obstacle. Let $\Omega = (\mathbb{T} \times (-2,2)) - \overline{B_1(\pi,0)}$, and consider
\eqn{
	\left\{ \begin{matrix}
		-\Delta \tb{u} + \nabla p = (\alpha,0), &  \text{on } \Omega,  \\
		\nabla \cdot \tb{u} = 0, &  \text{on } \Omega,  \\
		\tb{u} = \tb{0}, &   \text{for } y = \pm 2, \\ 
		\tb{u} = \tb{0}, &   \text{on } \partial B_1(\pi,0), \\
		\int_{-2}^2 u_1 \ dy = q, & \text{at } x =  0. \\  
	\end{matrix} \right. 
	\label{stokesEx3}
}

The last condition serves as a stand-in for explicit inlet-outlet boundary flow by assigning a certain flow rate $q$ to the horizontal component of the velocity. As the solution depends linearly on the value of $q$, we can set $q = 1$ in our tests without loss of generality. It also follows from the bulk incompressibility and no-slip boundary constraints that this condition could have been posed at any $x$ value. The constant forcing $(\alpha,0)$ acts as a Lagrange multiplier to enforce the flow-rate condition.  

We take $\Pi = \mathbb{T} \times (-2,2)$ to be the computational domain, the lateral walls $\mathbb{T} \times \{-2,2\}$ to be the intrinsic boundary $\partial \Omega_\text{I}$, and the obstacle wall $\partial B_1(\pi,0)$ as the auxiliary boundary $\partial \Omega_\text{A}$. An appropriate pressure basis can be built by combining trigonometric and Chebyshev polynomials on $[-2,2]$ in the form $\eta_{k_1,k_2}(x,y) = e^{ik_1x}T_{k_2}(y)$ for $k_1 \in \mathbb{Z}$ and $k_2 \geq 0$. However, as constant terms do not affect the pressure, we can discount $\eta_{0,0}$. 

The velocity basis can be built by assembling a complete scalar basis and following (\ref{AsymDer}), while ensuring the intrinsic boundary conditions are satisfied. We start by combining Fourier and polynomial bases as for the pressure basis above. Next, we require that all elements of the scalar basis possess roots of multiplicity two at $y = \pm 2$, to wit $\phi_{j_1,j_2}(x,y) = e^{ij_1x}T_{j_2}(y)(y^2-4)^2$ for $j_1 \in \mathbb{Z}$, $j_2 \geq 0$, so that both components of the corresponding $\Phi_{j_1,j_2}$ vanish at the lateral walls. This requirement however turns out to be too stringent if the scalar basis function is independent of $x$, in which case we only need its $y$-derivative to vanish at $y = \pm 2$. Note that $\{\phi_{0,j_2}\}_{j_2 \geq 0}$ already account for all the polynomials of degree at least four that have this property. Thus, we only need to include $\tilde \phi(x,y) = y(y^2 - 12)$, as no polynomial of lower degree possesses this property. Taken together, the collection $\{\tilde \Phi\} \cup \{\Phi_{j_1,j_2}\}_{j_1 \in \mathbb{Z}, j_2 \geq 0}$ forms a complete velocity basis.

The  flow rate condition in (\ref{stokesEx3}) is accounted for in the objective function (\ref{GSt2defn}) by including the additional term $\left|\sum_{j_1,j_2} c_{j_1,j_2} \int_{-2}^2 \left(\Phi_{j_1,j_2}\right)_1(0,y) \ dy + \tilde{c}\int_{-2}^2 \tilde \Phi_1(0,y) \ dy - q\right|^2$. Using (\ref{AsymDer}), this can be rewritten as
\eqn{
	&& \left|\sum_{j_1,j_2} c_{j_1,j_2}\int_{-2}^2 \py \phi_{j_1,j_2}(0,y) \ dy + \tilde{c}\int_{-2}^2 \py \tilde\phi(0,y) \ dy - q\right|^2 \nonumber\\
	&=& \left|\sum_{j_1,j_2 } c_{j_1,j_2}\left(\phi_{j_1,j_2}(0,2) - \phi_{j_1,j_2}(0,-2)\right) + \tilde{c}\left(\tilde\phi(0,2) - \tilde\phi(0,-2)\right) - q \right|^2. \label{FRcond}
}  

As $\phi_{j_1,j_2}$ vanish at $y = \pm 2$ for all $j_1, j_2$, they do not contribute towards (\ref{FRcond}) so this term is greatly simplified. The resulting objective function is
\eqn{
	G(\tb{c},\tb{d},\tilde c,\alpha) &=&  \norm{\sum_{(j_1,j_2) \in \mathcal{J}} c_{j_1,j_2} \Upsilon_{j_1,j_2} + \tilde c \tilde \Upsilon}^2_{L^2(\partial B_1(\pi,0))} + \left|\tilde{c}\left(\tilde\phi(0,2) - \tilde\phi(0,-2)\right) - q\right|^2 \nonumber\\
	&& \norm{\sum_{(j_1,j_2) \in \mathcal{J}} -c_{j_1,j_2} \Delta \Phi_{j_1,j_2} - \tilde{c} \Delta \tilde \Phi + \sum_{(k_1,k_2) \in \mathcal{K}} d_{k_1,k_2} \nabla  \eta_{k_1,k_2} - (\alpha,0)}^2_{L^2(\Omega)}, \label{GSt3defn}
}
where $\mathcal{J} = \{(j_1,j_2): j_1 \in \mathbb{Z} \text{ and } j_2 \geq 0\}$, $\mathcal{K} = \mathcal{J} \backslash \{(0,0)\}$, $\tb{c} = \left(c_{j_1,j_2}\right)_{(j_1,j_2) \in \mathcal{J}}$, and $\tb{d} = \left(d_{k_1,k_2}\right)_{(k_1,k_2) \in \mathcal{K}}$. In practice, we truncate these by setting a cut-off frequency: for a given $N_\text{e}$, set $\mathcal{J}(N_\text{e}) = \{(j_1,j_2) \in \mathcal{J}: |j_1| \leq N_\text{e} \text{ and } j_2 \leq N_\text{e}\}$ and similarly for $\mathcal{K}(N_\text{e})$.

\def \sclccm {0.09}

\begin{figure}[tbph]
	\centering
	\subfigure[Convergence of $u^{[N_\text{e}]}_{1,\text{e}}$]
	{\includegraphics[scale=\sclccm]{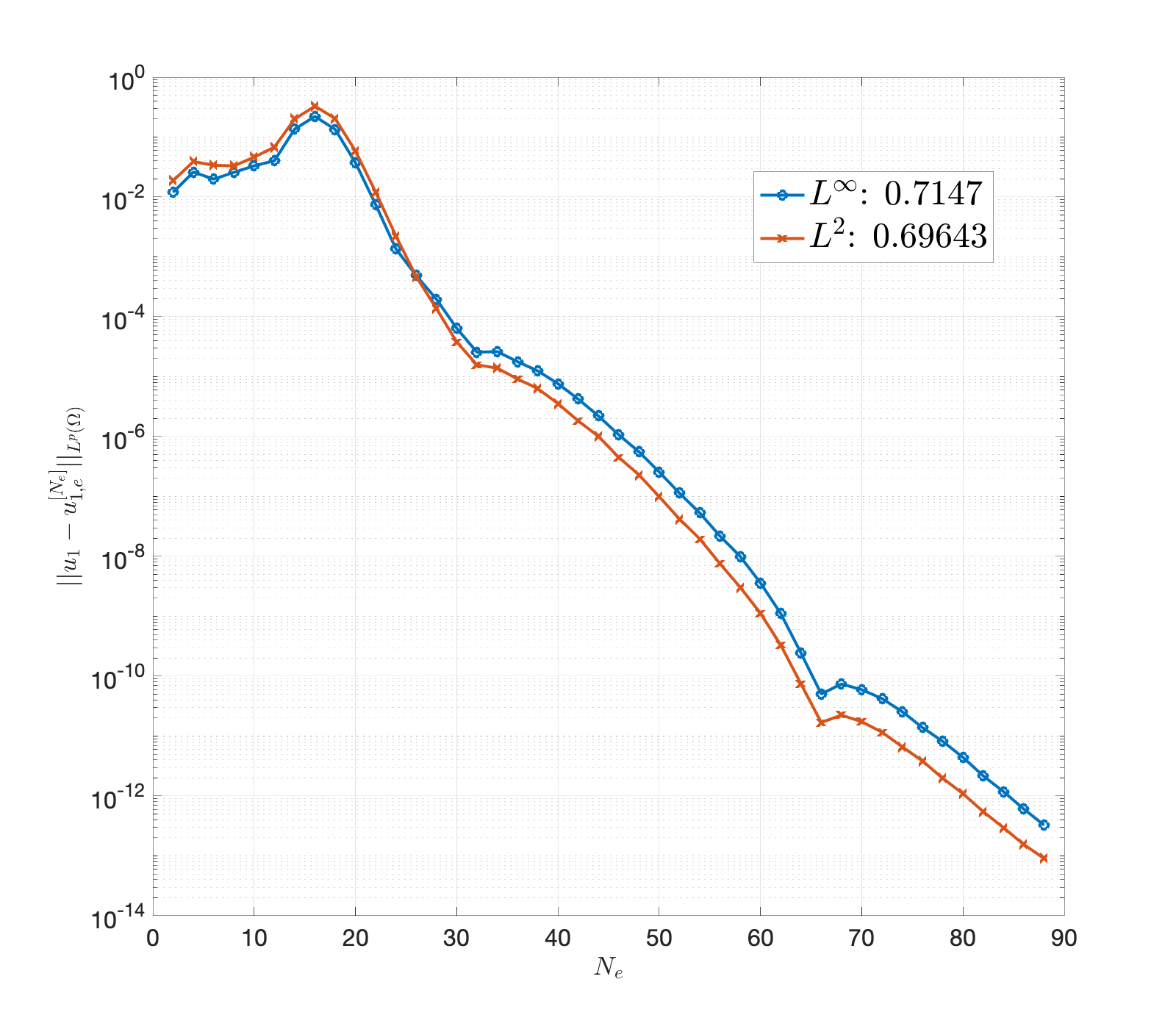}
		\label{Istokes3x}
	}
	\subfigure[Convergence of $u^{[N_\text{e}]}_{2,\text{e}}$]
	{\includegraphics[scale=\sclccm]{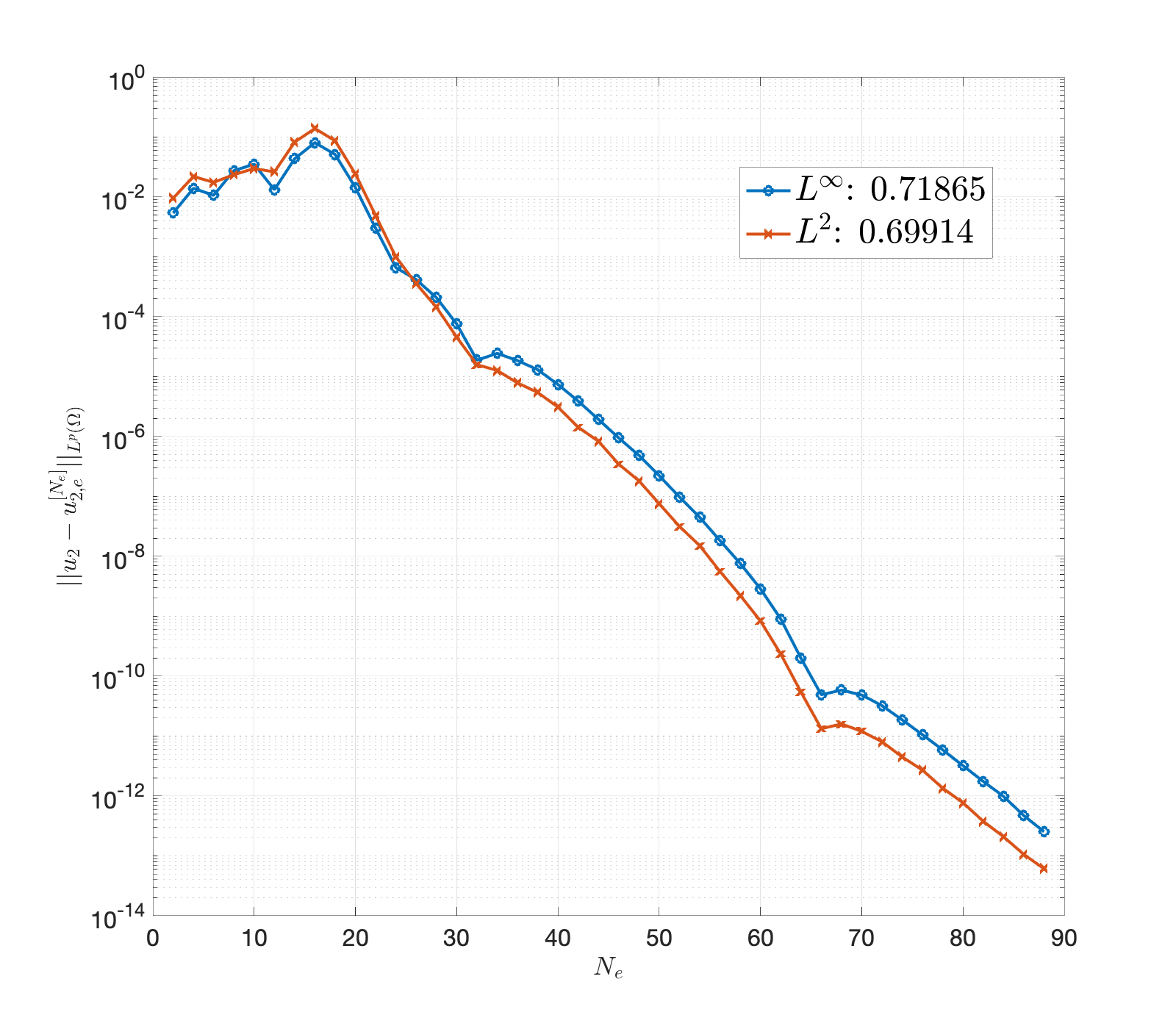}
		\label{Istokes3y}
	}
	\subfigure[Convergence of $p^{[N_\text{e}]}_\text{e}$]
	{\includegraphics[scale=\sclccm]{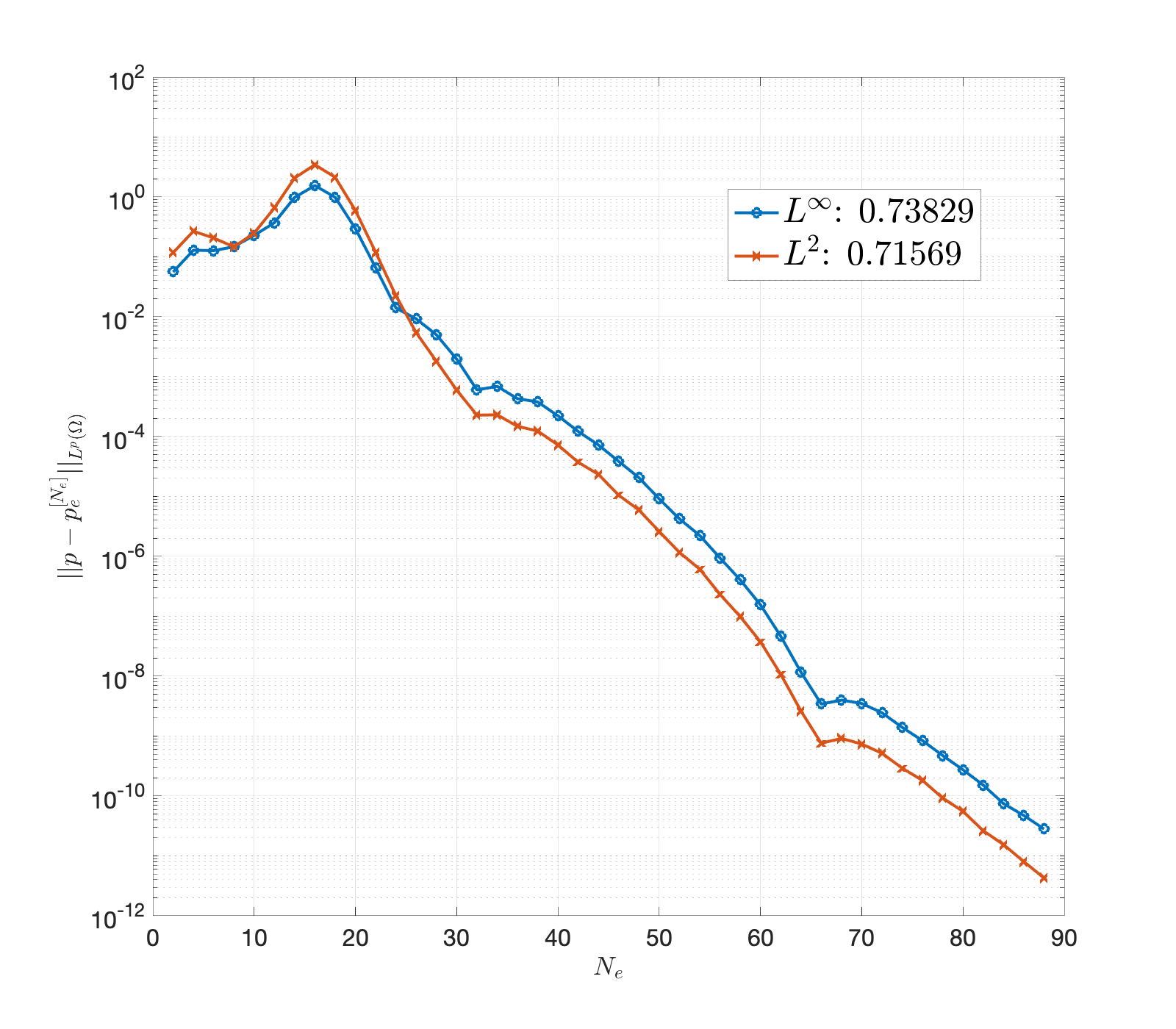}
		\label{Istokes3p}
	}
	\caption{Convergence results for velocity $u_1,u_2$ and pressure $p$ for the Stokes problem (\ref{stokesEx3}). For $N_\text{e} \gtrapprox 18$, the convergence is exponential and yields much smaller errors than the solvers used earlier.}\label{Istokes3}
\end{figure}

The changes to the technique also necessitate some modifications to the implementation. We keep the uniform grid in the $x$-direction, as it is suitable for integrating periodic functions, and couple it with a Gaussian quadrature scheme in the $y$-direction to accurately integrate the polynomials. More precisely, let $\{x_k\}_{k = 1}^{N_x}$ be equidistant points on $[0,2\pi)$ and let $\{y_l\}_{l = 1}^{N_y}$ be Gaussian quadrature nodes on $[-2,2]$ with associated weights $\{w_i\}_{i = 1}^{N_y}$. The inner product evaluation (\ref{Oip}) is replaced by
\eqn{
	\ip{\delta_1,\delta_2}_{\Omega} &=& \int_{\Pi} \chi_{\Omega}(\tb{x})\overline{\delta_1(\tb{x})} \delta_2(\tb{x}) \ d\tb{x} \approx h \sum_{k = 1}^{N_x} \sum_{l = 1}^{N_y} \chi_{\Omega}(x_k,y_l) \overline{\delta_1(x_k,y_l)} \delta_2(x_k,y_l)w_l.  \label{Oip2}
}

The resulting modifications to (\ref{pqrsmats}) and (\ref{sqsys4a}) are straightforward and are omitted here. Figure \ref{Istokes3} shows results for $N_x = 2^8$ and $N_y = 190$ as $N_\text{e}$ is varied. The errors are again computed by comparing successive solutions as the exact solution is unknown. As in problems (\ref{stokesEx2}) and (\ref{stokesEx1}), we see that a certain threshold of the number of extension functions is required beyond which solutions start converging. More strikingly, this approach allows us to attain much lower errors than the earlier examples of Stokes solvers, which is primarily a consequence of accounting for the boundary conditions in the design of the basis functions. This also has the benefit of significantly reducing the cost of solving the resulting linear systems as the number of boundary nodes is substantially reduced.

As another example of this approach, we solve the Stokes equation on a two-dimensional sphere. This problem shows up frequently while studying cell membranes, the dynamics of which can be modeled as viscous flows on two-dimensional curved domains \cite{henle2010hydrodynamics,saffman1975brownian}. While exact solutions can be found in simple cases, the problem becomes more complex if the surface has non-uniform curvature or if it contains obstacles. Our technique, however, is sufficiently flexible to enable us to tackle these challenges.   

Denote by $\Pi$ be the two-dimensional sphere of radius $\rho = 1$ parameterized by spherical coordinates, and define the physical domain
\eqn{
	\Omega = \{(\theta,\phi) \in \Pi: \sin(\theta)\sin(\phi) < 0.8\}. \label{spheredef}
} 

In Cartesian coordinates, this can be seen to be a spherical surface of unit radius missing the skullcap given by $y \geq 0.8$, which models an obstacle. The Stokes system is given by

\eqn{
	\left\{ \begin{matrix}
		(-\Delta-K) \tb{u} + \nabla p = \tb{f}, &  \text{on } \Omega,  \\
		\nabla \cdot \tb{u} = 0, &  \text{on } \Omega,  \\
		\tb{u} = \tb{g}, &   \text{on } \partial \Omega. \\
	\end{matrix} \right. 
	\label{stokesEx4}
}

Here, $K$ is the curvature of the surface which, on a sphere, takes the constant value $K = 1/\rho^2$. The boundary is simply the circle of radius $\sqrt{1-0.8^2} = 0.6$ centered at $(0,0.8,0)$, lying in the $y = 0.8$ plane.

We again construct appropriate velocity and pressure bases on $\Pi$, and use them to minimize the objective function (\ref{GSt2defn}). A helpful observation is that the family of spherical harmonics $\{Y_{l}^m\}_{l \geq 0, |m| \leq l}$ is an eigenbasis for the Laplacian $(-\Delta)$ on $\Pi$ with corresponding eigenvalues $l(l+1)$. Recall that these functions are given by
\eqn{
	Y_l^m(\theta,\phi) = e^{im\phi}P_l^m(\cos(\theta)), \label{spharmdef}
}
where $\{P_l^m\}$ are the associated Legendre polynomials. The completeness of the harmonics allows us to use them as the pressure basis (without the constant $Y_0^0$ element), and to generate the velocity basis by computing their anti-symmetric derivatives
\eqn{
	\Phi_l^m(\theta,\phi) = \pt Y_l^m(\theta,\phi)\hat {\bs{\phi}} -  \frac{1}{\sin(\theta)}\pph Y_l^m(\theta,\phi) \hat{\bs{\theta}}. \label{sphvelbasis}
} 

The eigenvalue property of this basis with respect to the Laplacian, and hence the operator $\mathcal{L} = (-\Delta-K)$, greatly simplifies the form of the objective function. We again set a cutoff frequency and only access velocity and pressure basis functions from $\mathcal{J}(N_\text{e}) = \{(l,m): 0\leq l \leq N_\text{e} \text{ and } |m|\leq l\}$ and $\mathcal{K}(N_\text{e}) = \mathcal{J}(N_\text{e})\backslash \{(0,0)\}$ respectively. For the inner product evaluation, we use a uniform grid $\{\phi_k\}_{k = 1}^{N_{\phi}}$ along the azimuthal angle to leverage the periodicity, with a high-order Gaussian quadrature scheme $\{(\theta_l,w_l)\}_{l = 1}^{N_{\theta}}$ along the polar angle. As a result, we have
\eqn{
	\ip{\delta_1,\delta_2}_{\Omega} &=& \int_0^{2\pi}\int_0^{\pi} \chi_{\Omega}(\theta,\phi) \overline{\delta_1(\theta,\phi)} \delta_2(\theta,\phi)  \sin(\theta) \ d\theta \ d\phi \nonumber\\
	&\approx& \left(\frac{2\pi}{N_{\phi}}\right) \sum_{k = 1}^{N_\phi} \sum_{l = 1}^{N_\theta} \chi_{\Omega}(\theta_l,\phi_k) \overline{\delta_1(\theta_l,\phi_k)} \delta_2(\theta_l,\phi_k) \sin(\theta_l) w_l.  \label{Oip3}
} 

Note that $\sin(\theta_l) > 0$ and $w_l > 0$ for all $l$; as a result, we can compute square-roots of $\sin(\theta_l)w_l$ and attach them to the basis function terms while setting up the matrices, as in (\ref{pqrsmats}).

\begin{figure}[tbph]
	\centering
	\subfigure[Convergence of $u^{[N_\text{e}]}_{\phi,\text{e}}$]
	{\includegraphics[scale=\sclccm]{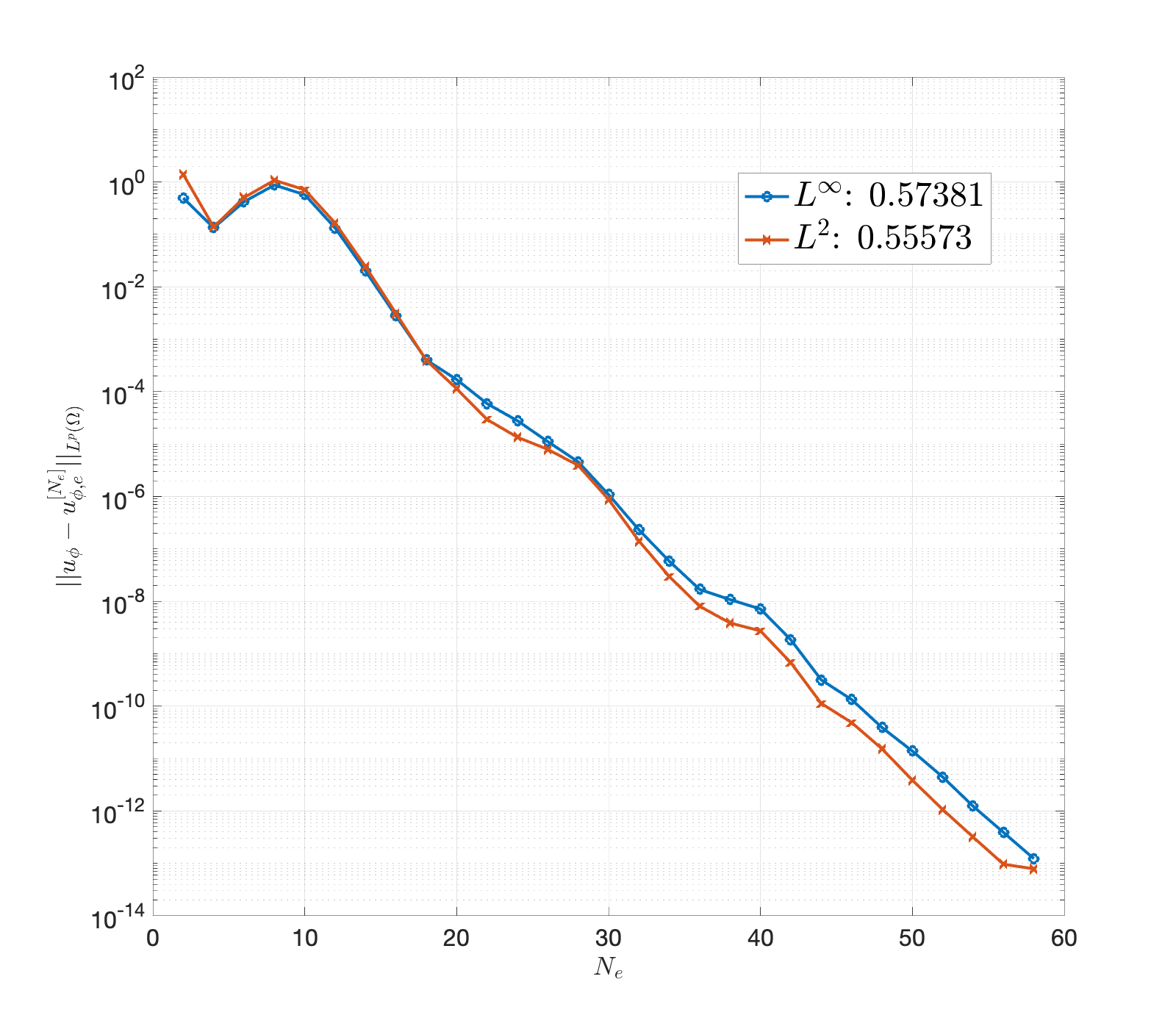}
		\label{Istokes4phi}
	}
	\subfigure[Convergence of $u^{[N_\text{e}]}_{\theta,\text{e}}$]
	{\includegraphics[scale=\sclccm]{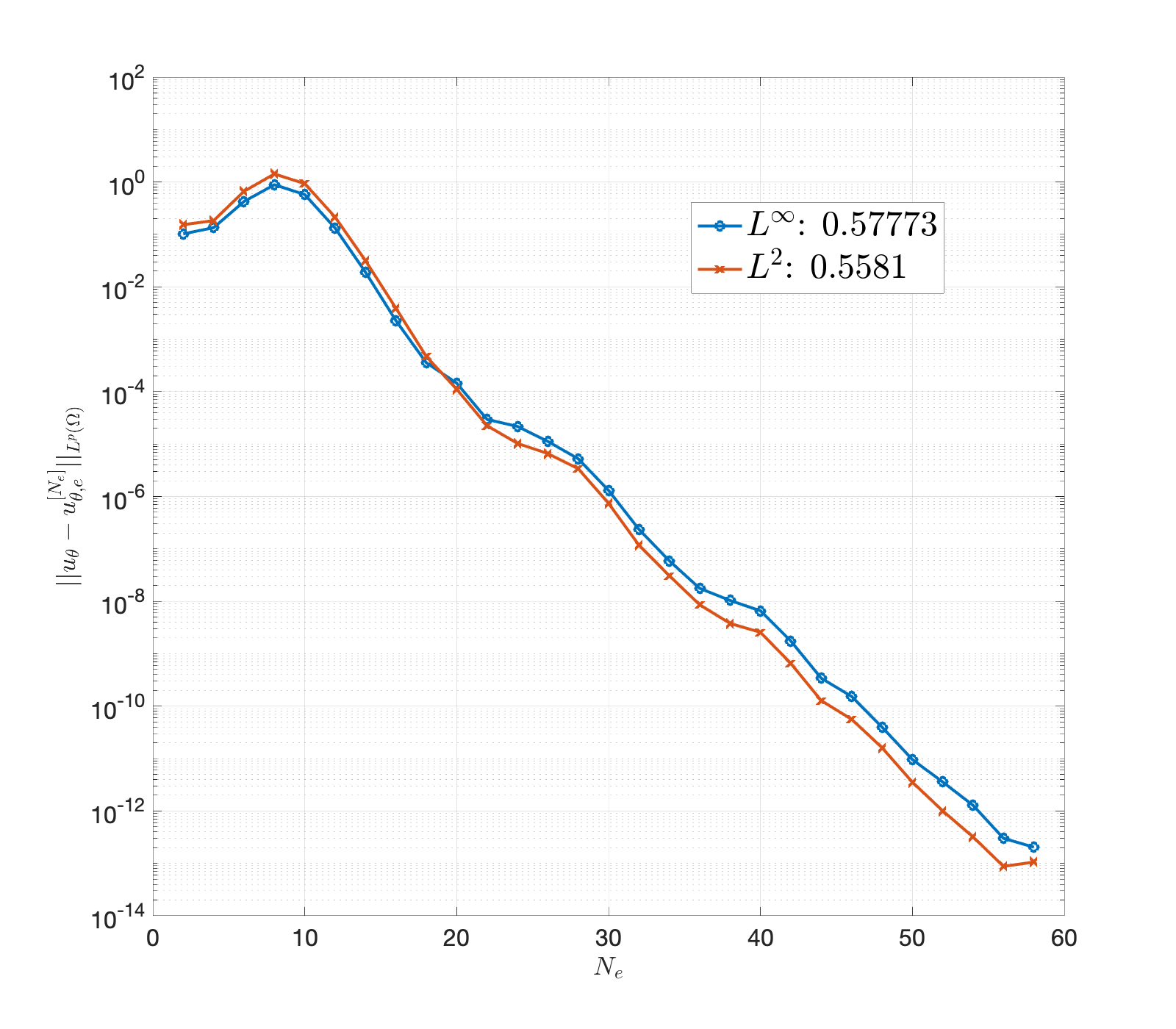}
		\label{Istokes4theta}
	}
	\subfigure[Convergence of $p^{[N_\text{e}]}_\text{e}$]
	{\includegraphics[scale=\sclccm]{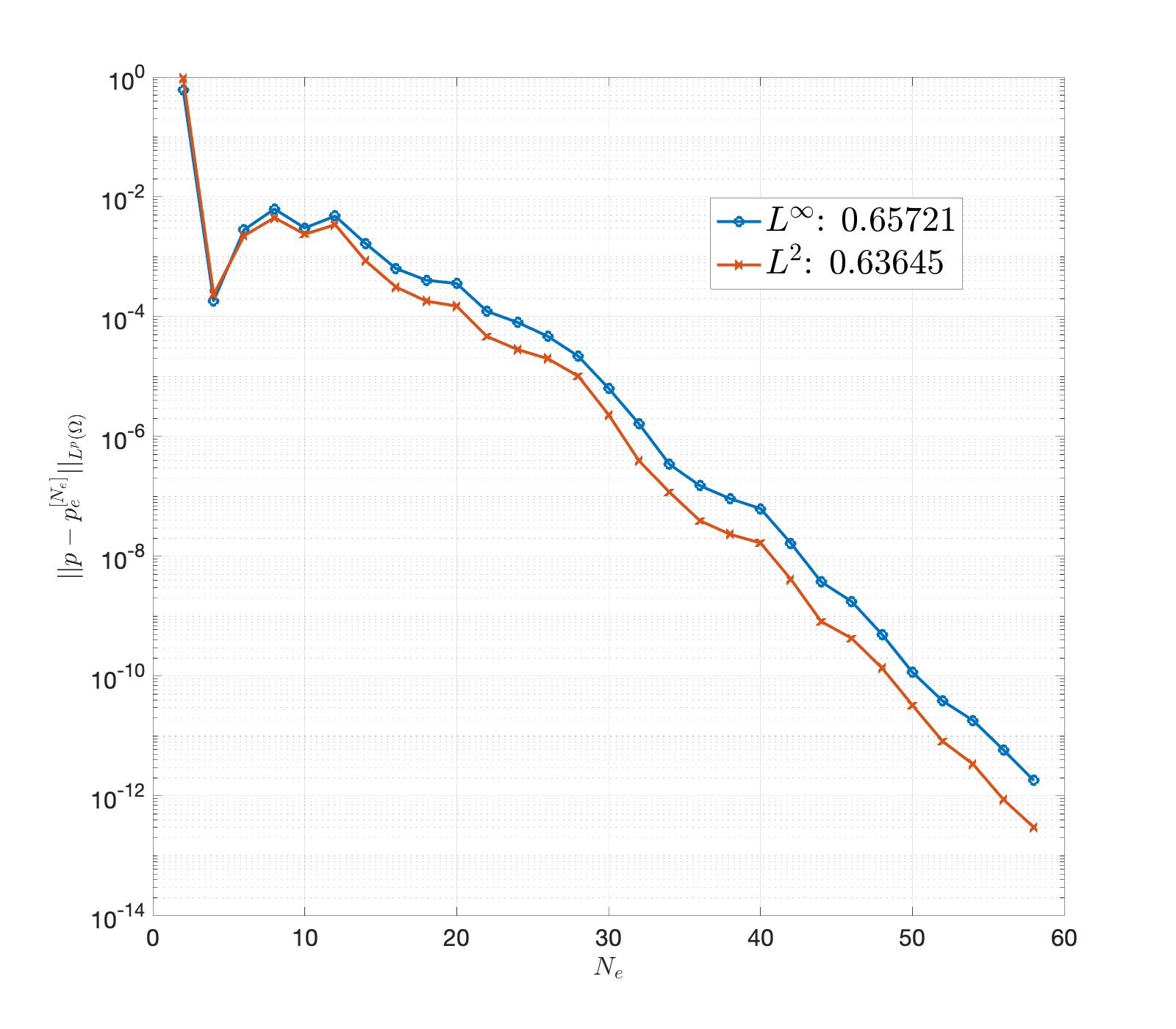}
		\label{Istokes4p}
	}
	\caption{Convergence results for velocity $u_\phi,u_\theta$ and pressure $p$ for Stokes flow (\ref{stokesEx4}) on a two-dimensional sphere. Beyond $N_\text{e} \gtrapprox 12$, the convergence is exponential.}\label{Istokes4}
\end{figure}

We solve (\ref{stokesEx4}) with $\tb{g} = \tb{0}$ and
\eqn{
	\tb{f}(\theta,\phi) = \left(\sin(\phi)\cos(\theta) + \cos(2\phi)\sin(\theta)\right)\hat{\bs{\phi}} + \left(-\cos(\phi) + \sin(2\phi)\sin(\theta)\cos(\theta)\right)\hat{\bs{\theta}}. \label{stokes4forc} 
}

As the exact solution is unknown, we again use successively refined solutions to compute the errors. The resulting convergence plots, shown in Figure \ref{Istokes4} for $N_{\phi} = 2^7$ and $N_{\theta} = 140$, demonstrate that all the quantities converge rapidly as $N_\text{e}$ is varied. We remark that this strategy can be easily generalized to surfaces with non-constant curvature, provided a computational domain with a complete basis can be found.

\subsection{The Navier--Stokes System}\label{navierstokessub}

We now briefly describe how the approaches developed for Stokes equation in Subsections \ref{stokessub} and \ref{AltStokes} can be extended to the Navier--Stokes problem. The key insight is that using a multi-step method for discretizing the time derivative yields a sequence of iterative Stokes equations, somewhat akin to our approach for solving the heat equation in Section \ref{SimpTests}. 

Consider the Navier--Stokes system
\eqn{
	\left\{ \begin{matrix}
		\tb{u}_t + (\tb{u} \cdot \nabla)\tb{u} = -\nabla p + \Delta \tb{u} + \tb{f}
		, &  \text{for } \tb{x} \in \Omega, t > 0,  \\
		\nabla \cdot \tb{u} = 0, &  \text{for } \tb{x} \in  \Omega, t > 0,  \\
		\tb{u} = \tb{g}, &   \text{for } \tb{x} \in \partial \Omega, t > 0, \\ 
		\tb{u}(0,\tb{x}) = \tb{u}_0(\tb{x}), & \text{for } \tb{x} \in \Omega. \\
	\end{matrix} \right. 
	\label{NSEx}
}

\begin{figure}[tbph]
	\centering
	\subfigure[Problem (\ref{NSEx0})]
	{\includegraphics[scale=\sclcm]{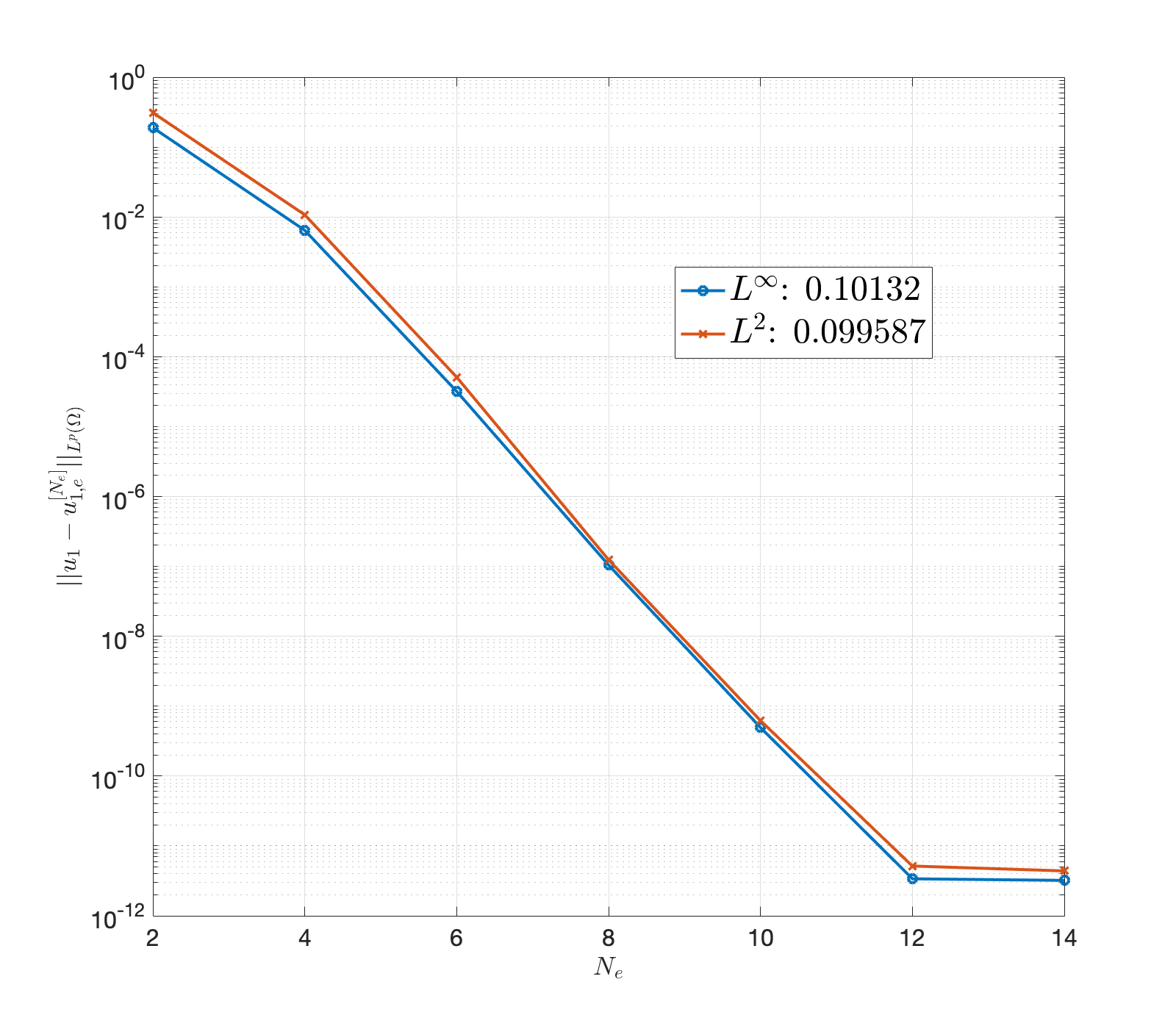}
		\includegraphics[scale=\sclcm]{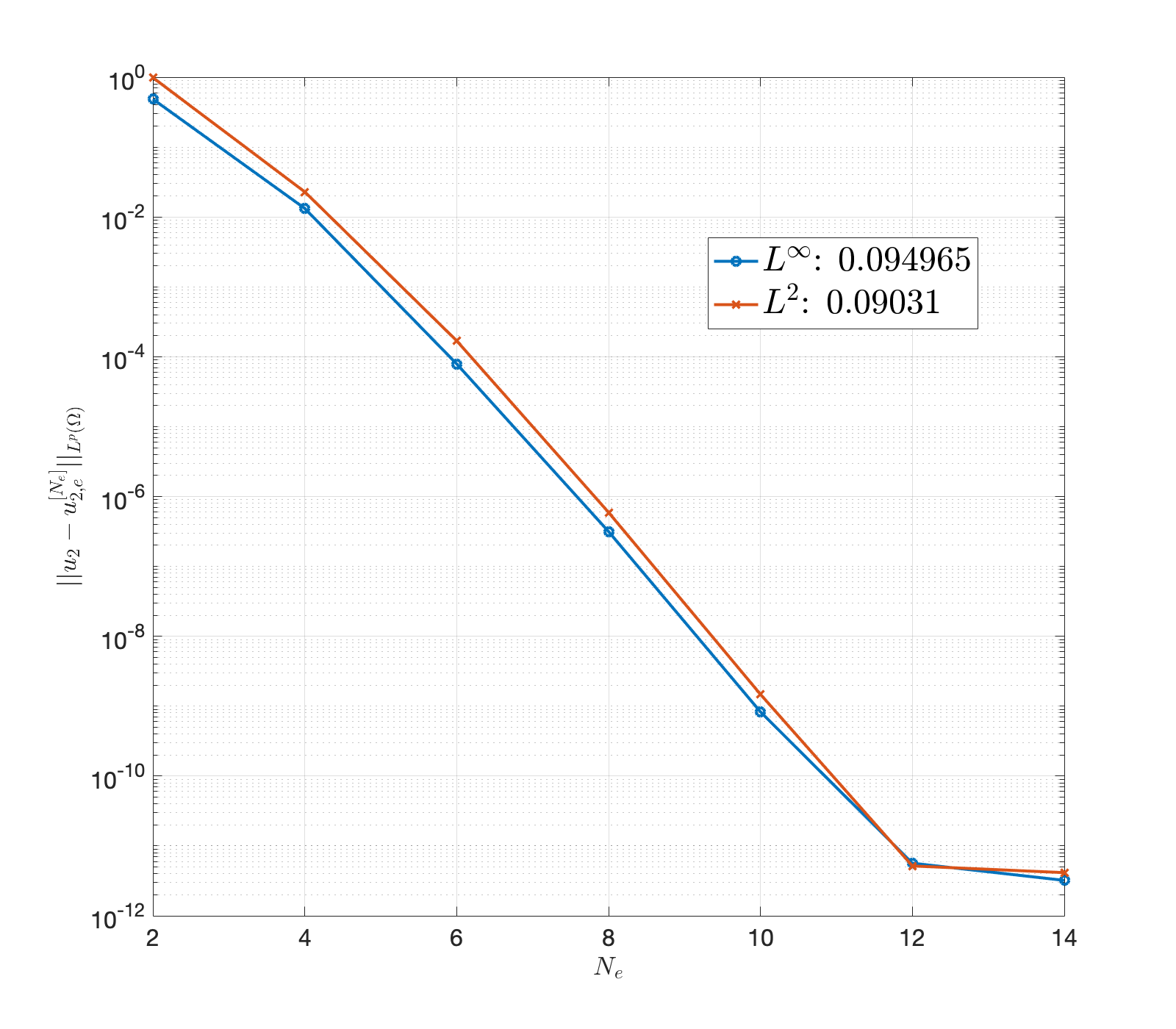}
		\includegraphics[scale=\sclcm]{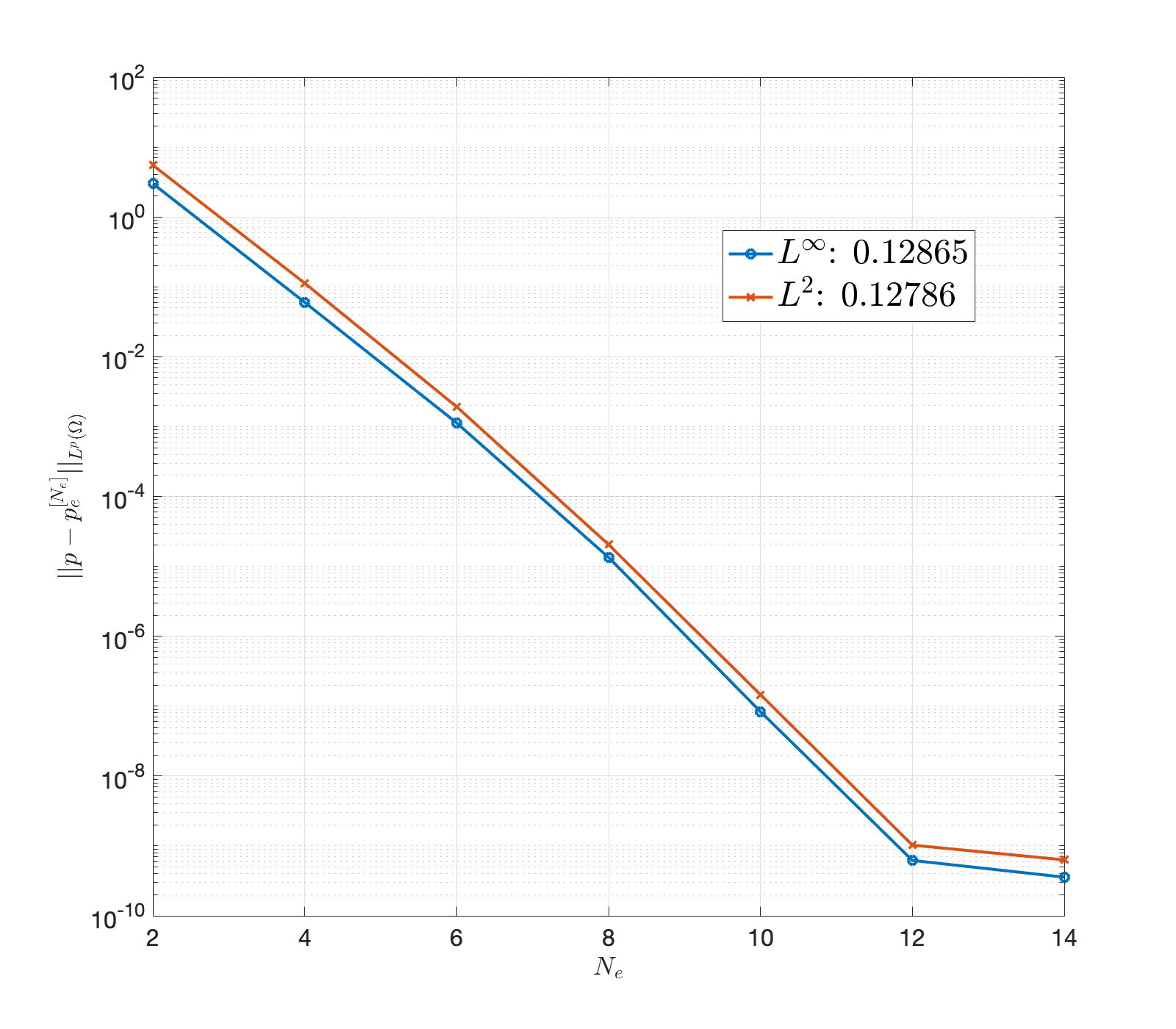}
		\label{NSex0}
	}
	\subfigure[Same forcing as (\ref{NSEx0}) with homogeneous Dirichlet boundary conditions]
	{\includegraphics[scale=\sclcm]{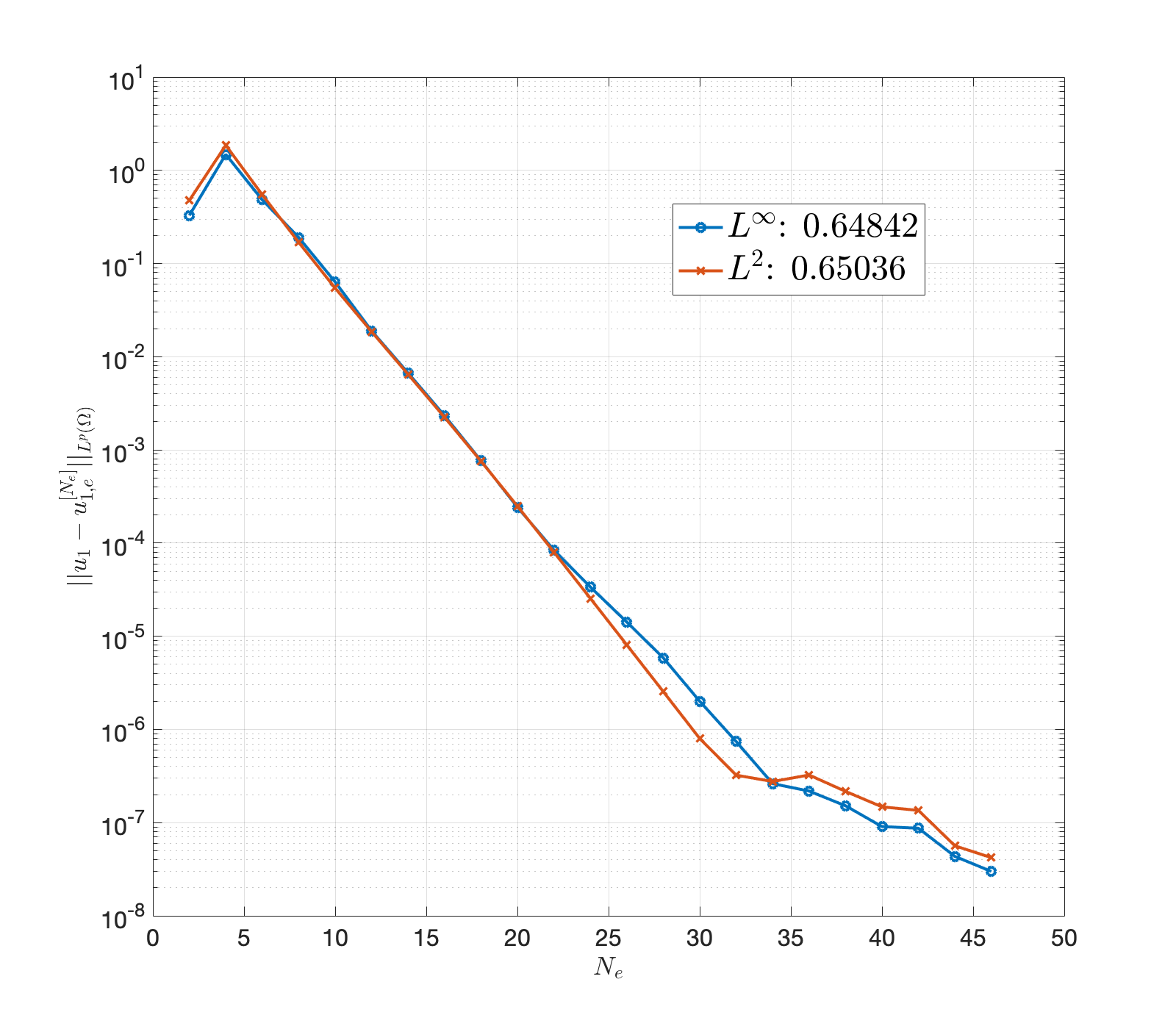}
		\includegraphics[scale=\sclcm]{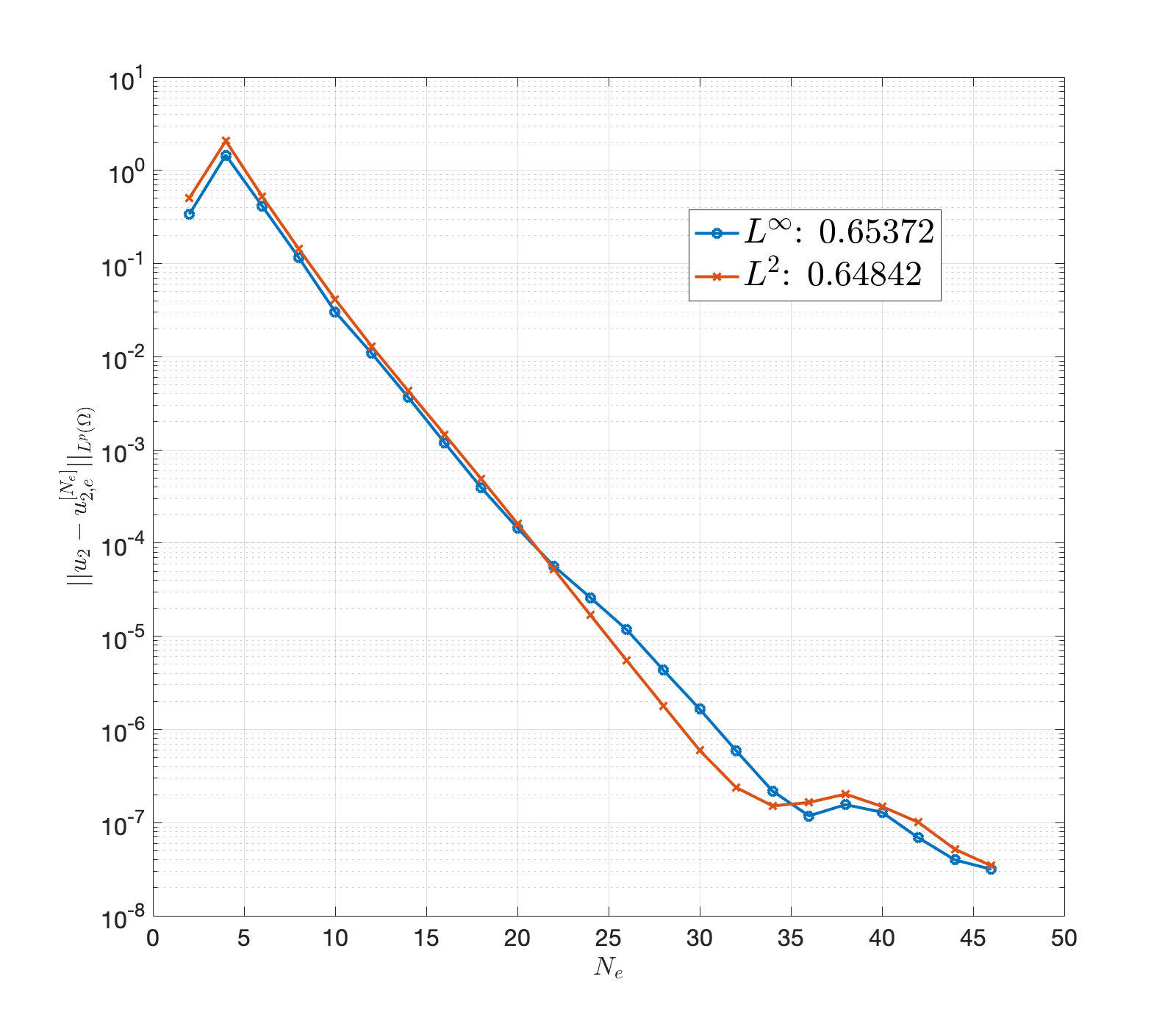}
		\includegraphics[scale=\sclcm]{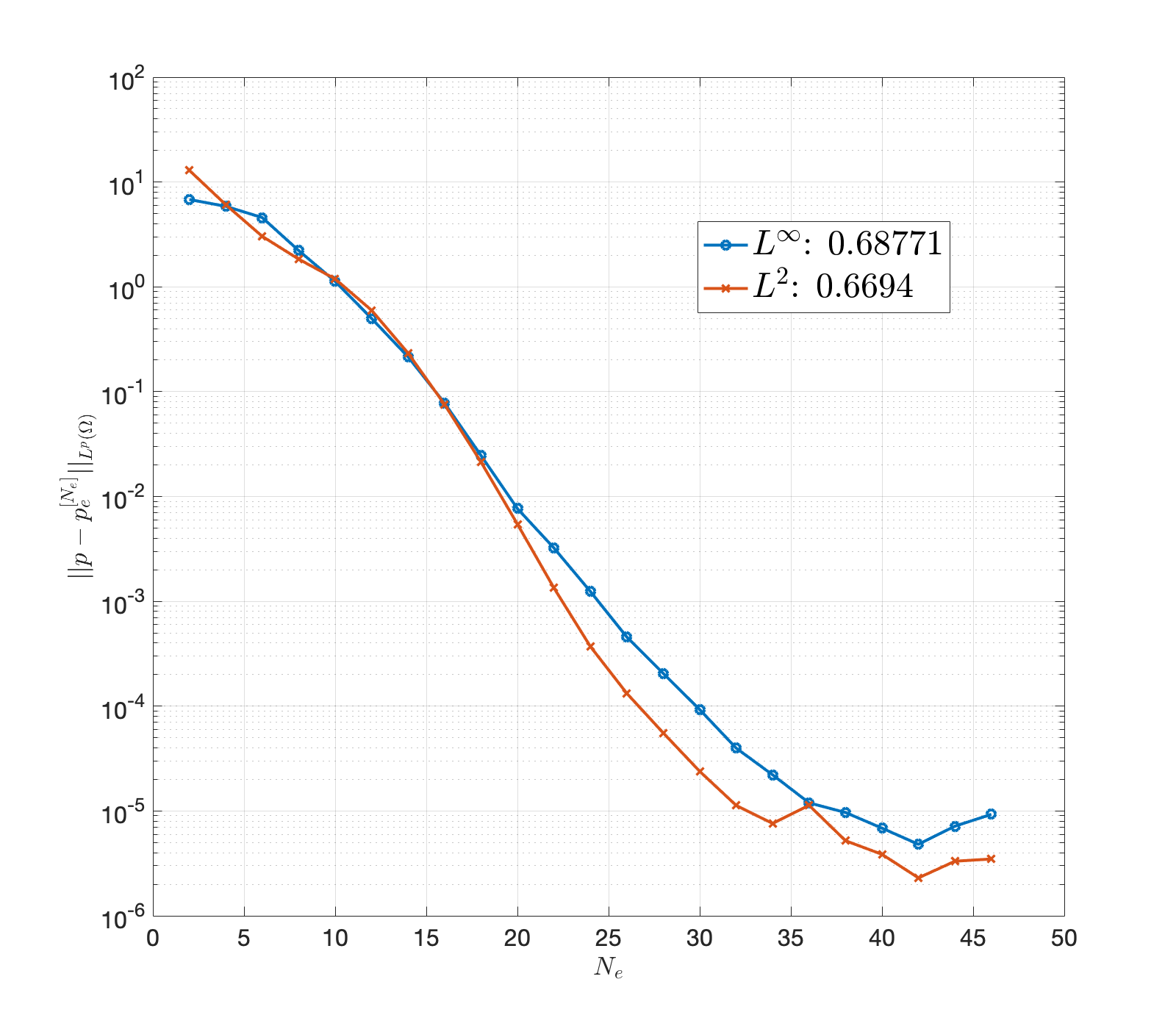}
		\label{NSex2}
	}
	\subfigure[Zero forcing with inflow-outflow boundary conditions (\ref{stokesEx1})]
	{\includegraphics[scale=\sclcm]{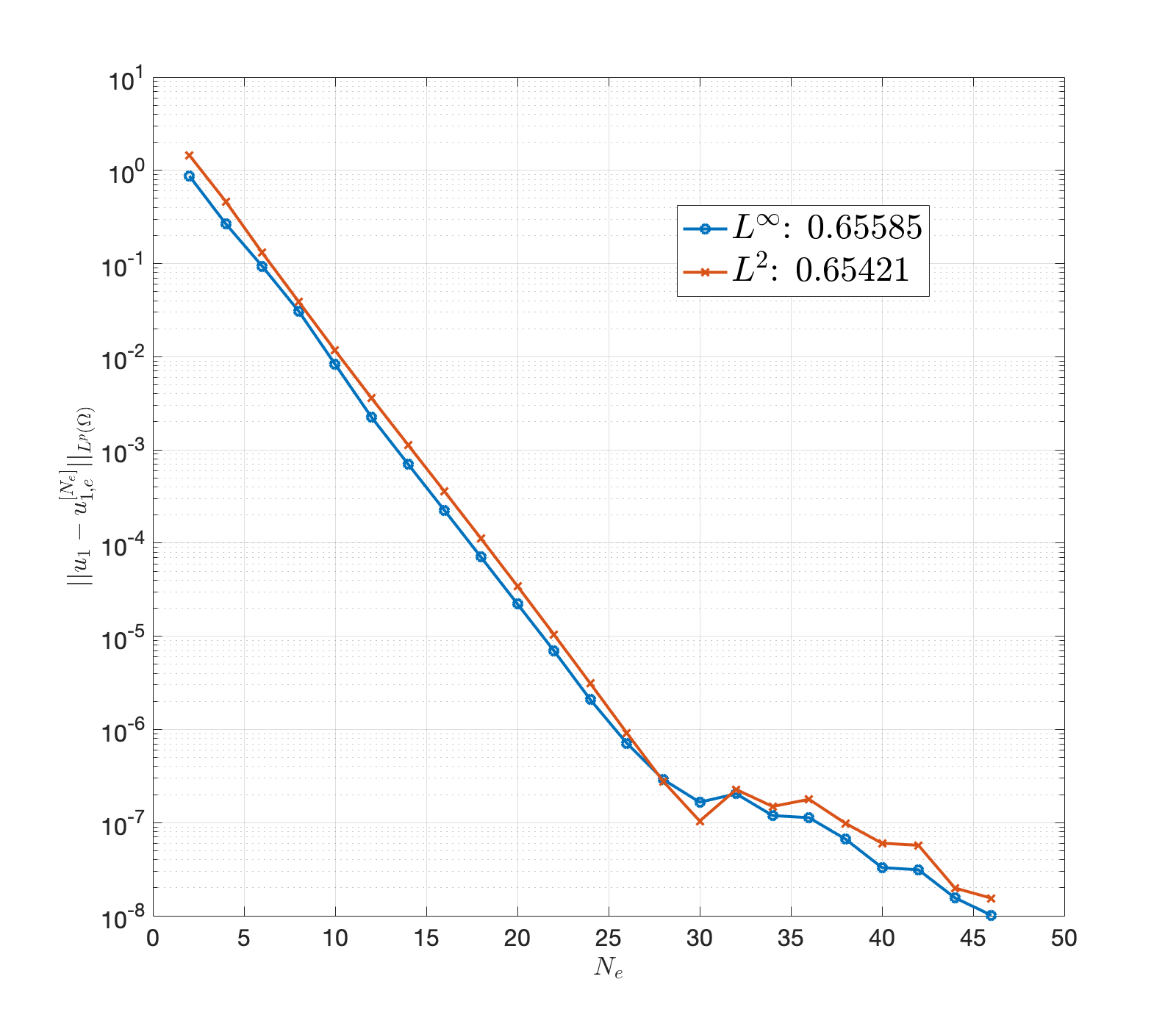}
		\includegraphics[scale=\sclcm]{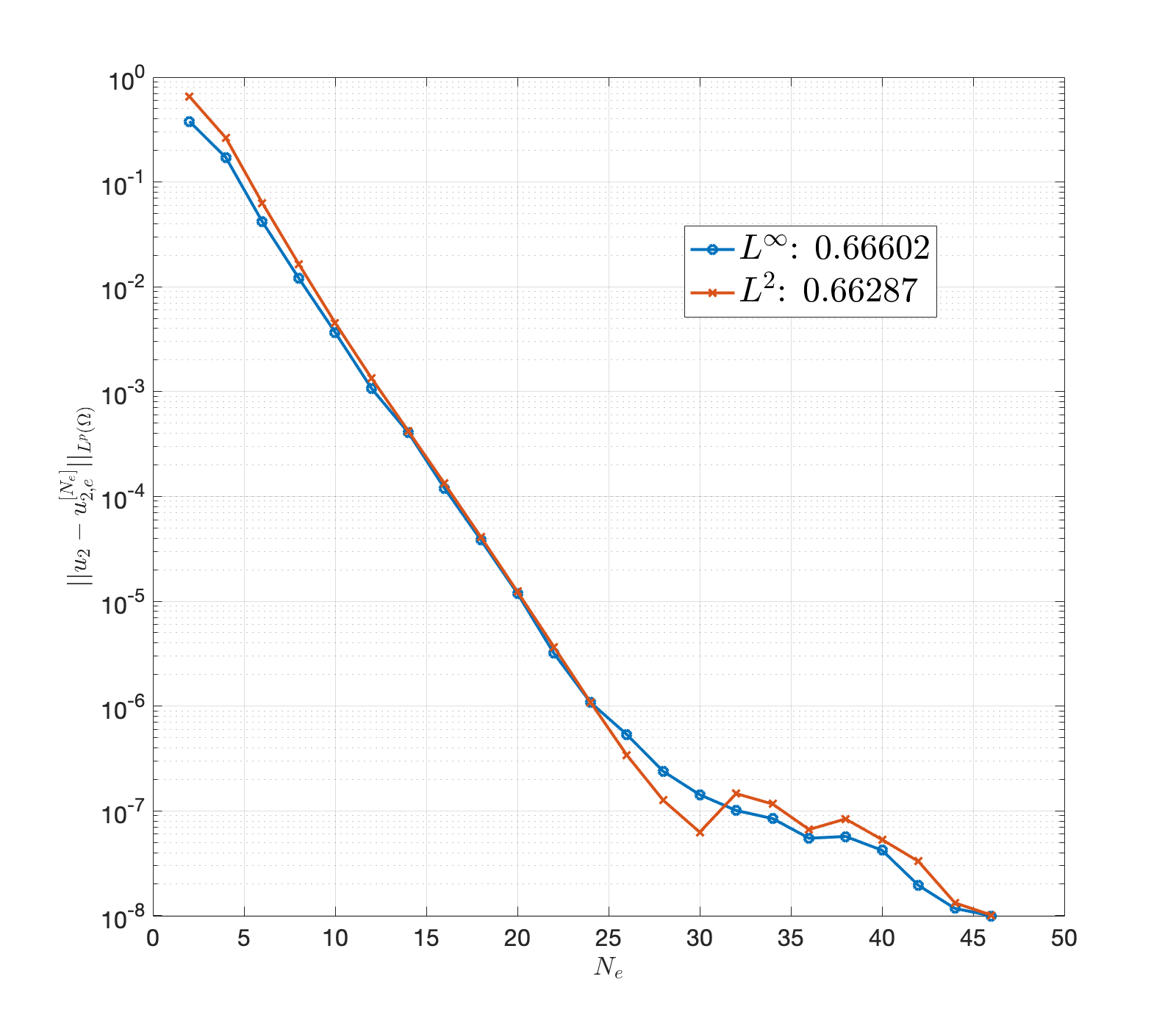}
		\includegraphics[scale=\sclcm]{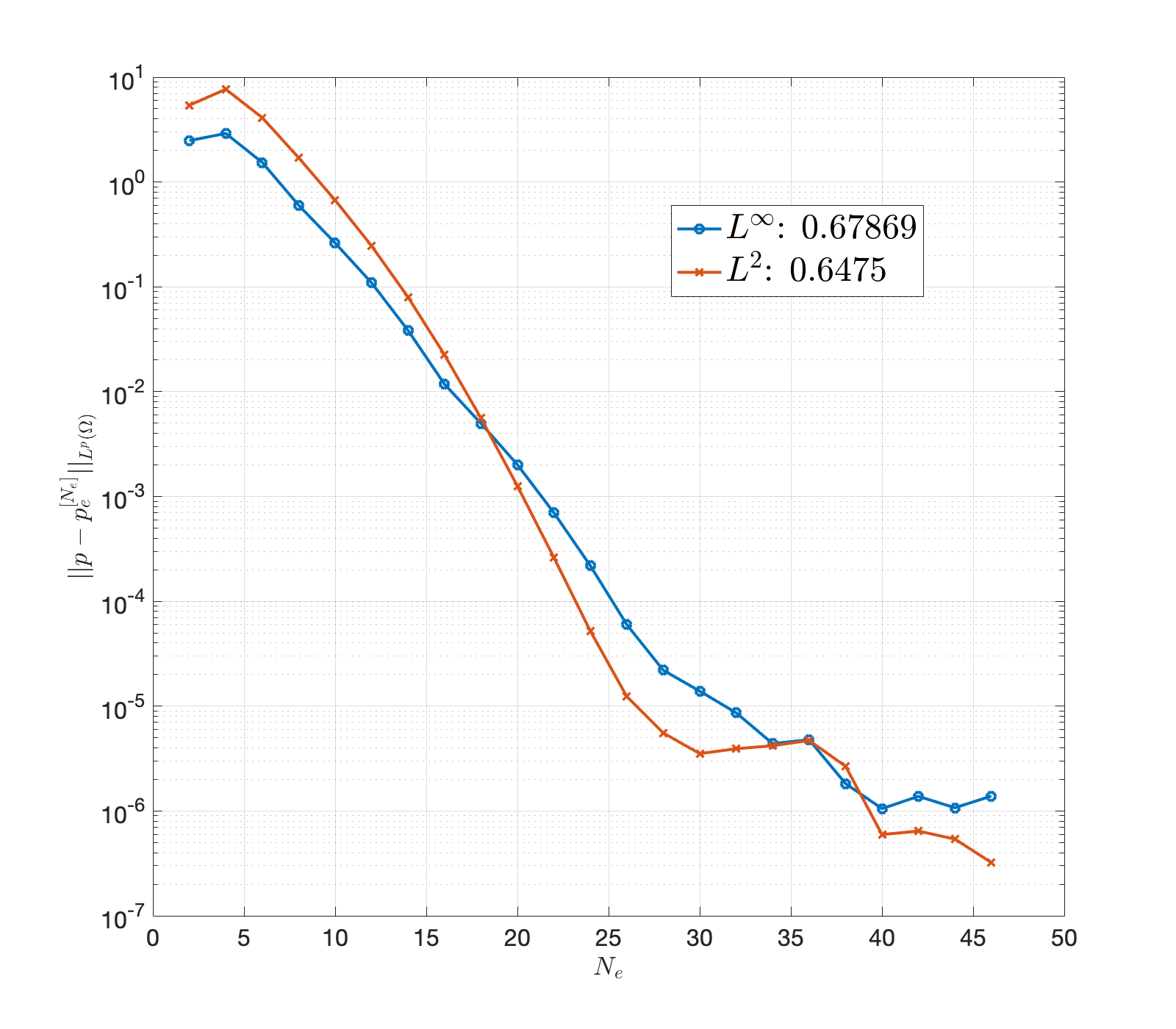}
		\label{NSex1}
	}
	\caption{Convergence plots for the solutions $u_1$, $u_2$ and $p$ to the Navier--Stokes test problems on the domain shown in Figure \ref{rectdom}. All the unknowns appear to converge exponentially at $O(a^{N_\text{e}})$; the corresponding best-fit values of $a$ are shown for the different norms.
	}\label{NSconv}
\end{figure}

We again use the BDF-4 scheme to discretize the time-derivative. However, the nonlinear advection term in the momentum equation makes a fully implicit solver hard to implement. Instead, we treat this team explicitly; specifically, we have, following \cite{hundsdorfer2007imex},
\eqn{
	&& \left(\mathbb{I} - \frac{12(\Delta t)}{25} \Delta\right)\tb{u}^{n+1} + \left(\frac{12(\Delta t)}{25}\right)\nabla p^{n+1} = \frac{48\tb{u}^{n} - 36\tb{u}^{n-1} + 16\tb{u}^{n-2} - 3\tb{u}^{n-3}}{25} + \nonumber\\
	&& \left(\frac{12(\Delta t)}{25}\right) \left(\tb{f}^{n+1} - 4(\tb{u}^{n} \cdot \nabla)\tb{u} ^{n} + 6(\tb{u}^{n-1} \cdot \nabla)\tb{u}^{n-1}   - 4(\tb{u}^{n-2} \cdot \nabla)\tb{u}^{n-2}   + (\tb{u}^{n-3} \cdot \nabla)\tb{u}^{n-3} \right). \nonumber\\	\label{NSbdf4}
}

Together with the incompressibility and boundary constraints 
\eqn{
	\left\{ \begin{matrix}
		\nabla \cdot \tb{u}^{n+1} = 0, &  \text{for } \tb{x} \in  \Omega,  \\
		\tb{u}^{n+1} = \tb{g}^{n+1}, &   \text{for } \tb{x} \in \partial \Omega, \\ 
	\end{matrix} \right. 
	\label{NScontraints}
}
each iteration simply requires solving (\ref{stokes}) with $\mathcal{L} = \mathbb{I} - \frac{12(\Delta t)}{25} \Delta$; this operator is invertible on $\Pi$ so we do not need to account for a non-trivial null space. The backward Euler scheme is employed for the first few time-steps to initialize the iterative procedure.

We test our technique on the domain shown in Figure \ref{rectdom} by varying the forcing and boundary conditions. Firstly, we use the exact solution
\eqn{
	\tb{u}(t,x,y) =  \begin{pmatrix}
		e^{\sin(x)}\cos(y) \\
		-e^{\sin(x)}\sin(y)\cos(x)
	\end{pmatrix}e^t, \qquad p(t,x,y) = e^{2\cos(x)}e^t \label{NSEx0}
}
to find the corresponding $\tb{f}$ and $\tb{g}$, and use our method to recover the solutions. Next, we use the same forcing $\tb{f}$ as given by (\ref{NSEx0}), but change the boundary conditions to no-slip everywhere. Finally, we solve the problem with zero external forcing and prescribed flow conditions at the inlet-outlet walls, as in example (\ref{stokesEx1}) and Figure \ref{flowBC}. The initial condition $\tb{u}_0$ used for all three tests is obtained from (\ref{NSEx0}).

The parameter values $N$ and $n_\text{b}$ are the same as in the Stokes tests in Subsection \ref{stokessub}. The problems are solved till $T = 1$, with time-step size $\Delta t = 10^{-3}$ for both the backward Euler and BDF-4 schemes. The exponential convergence is demonstrated in Figure \ref{NSconv}; as seen earlier, the problem with a known exact solution converges rapidly while the other two are markedly slower.

\begin{figure}[tbph]
	\centering
	\subfigure[Problem (\ref{NSEx3})]
	{\includegraphics[scale=\sclcm]{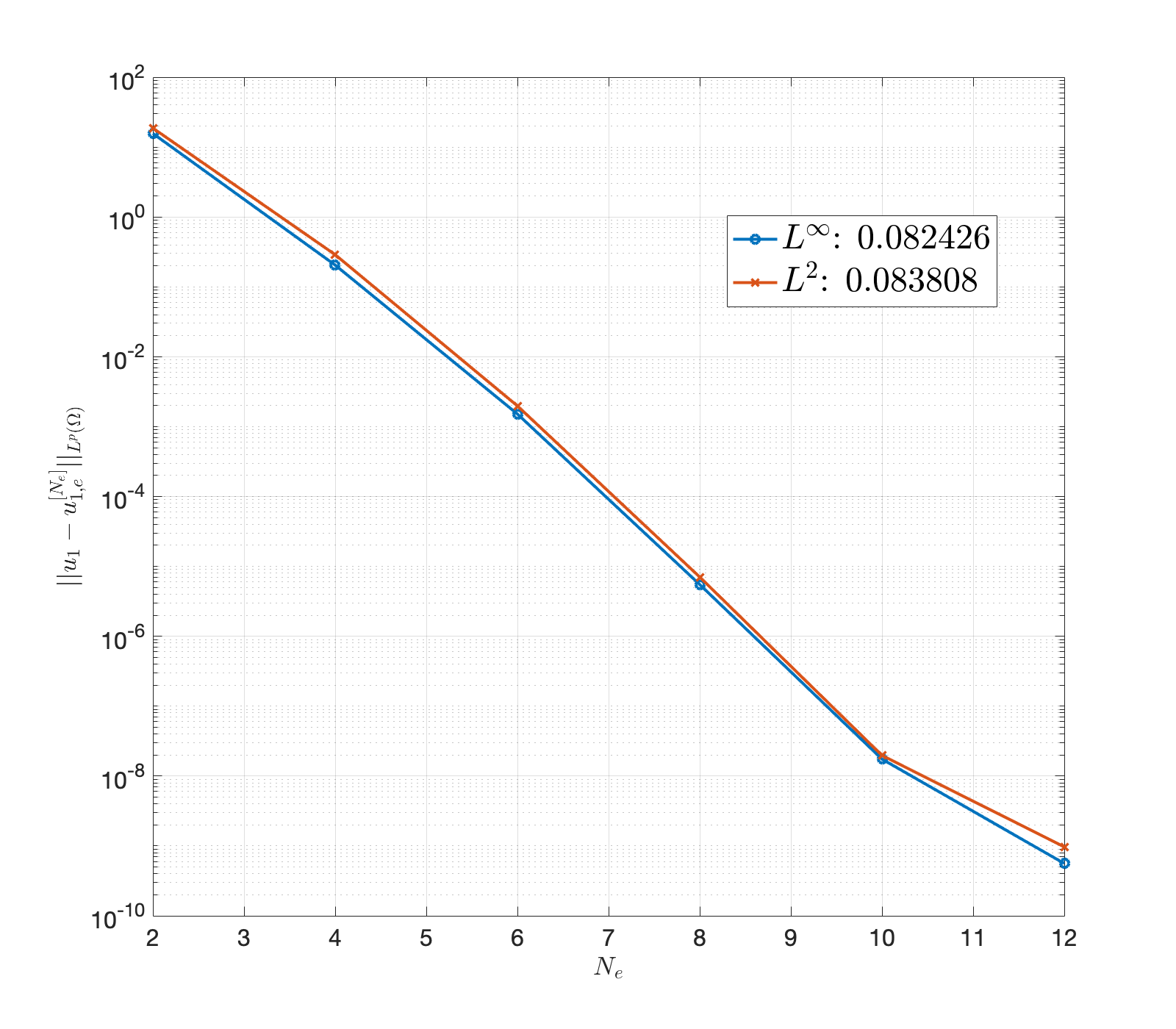}
		\includegraphics[scale=\sclcm]{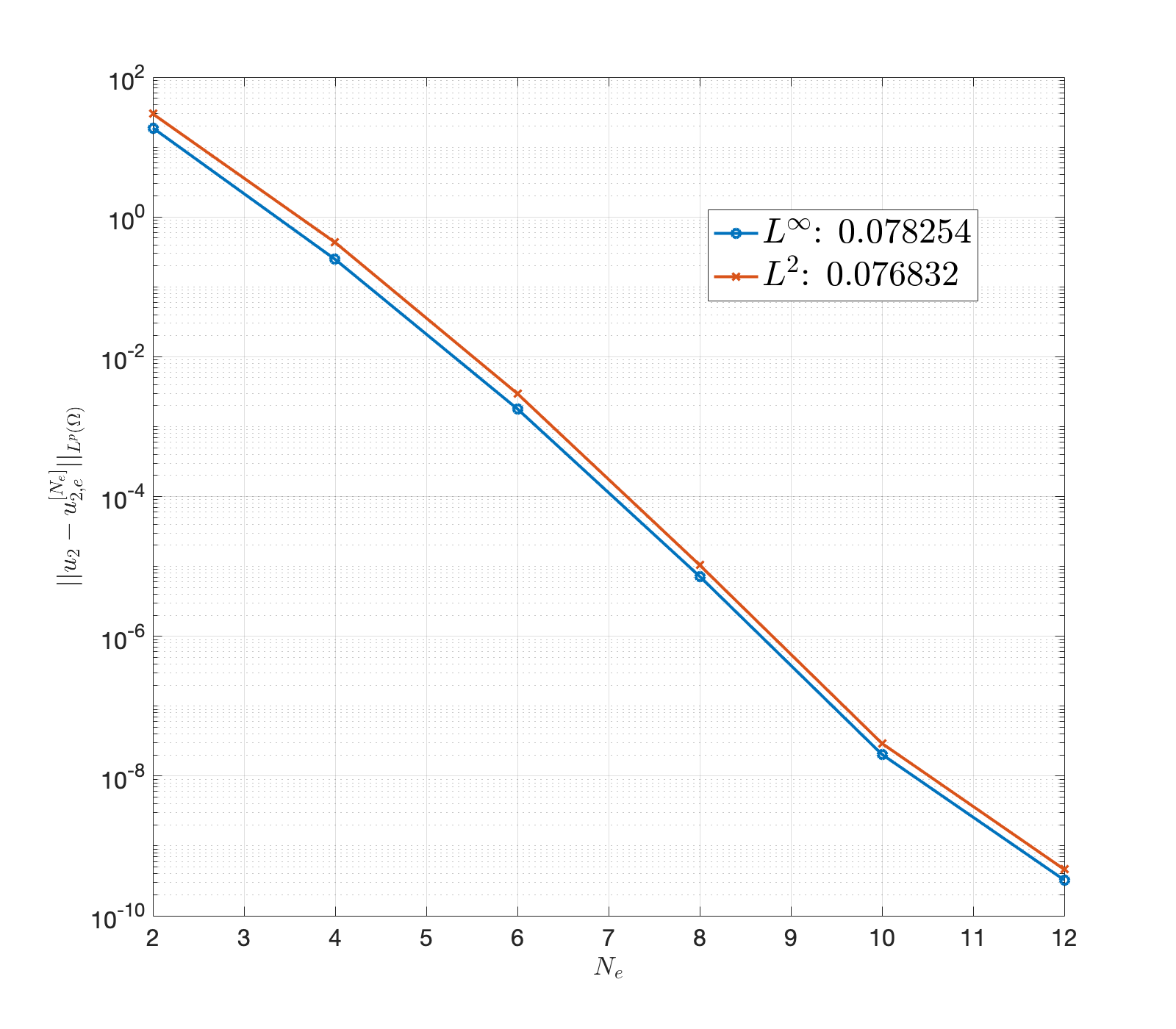}
		\includegraphics[scale=\sclcm]{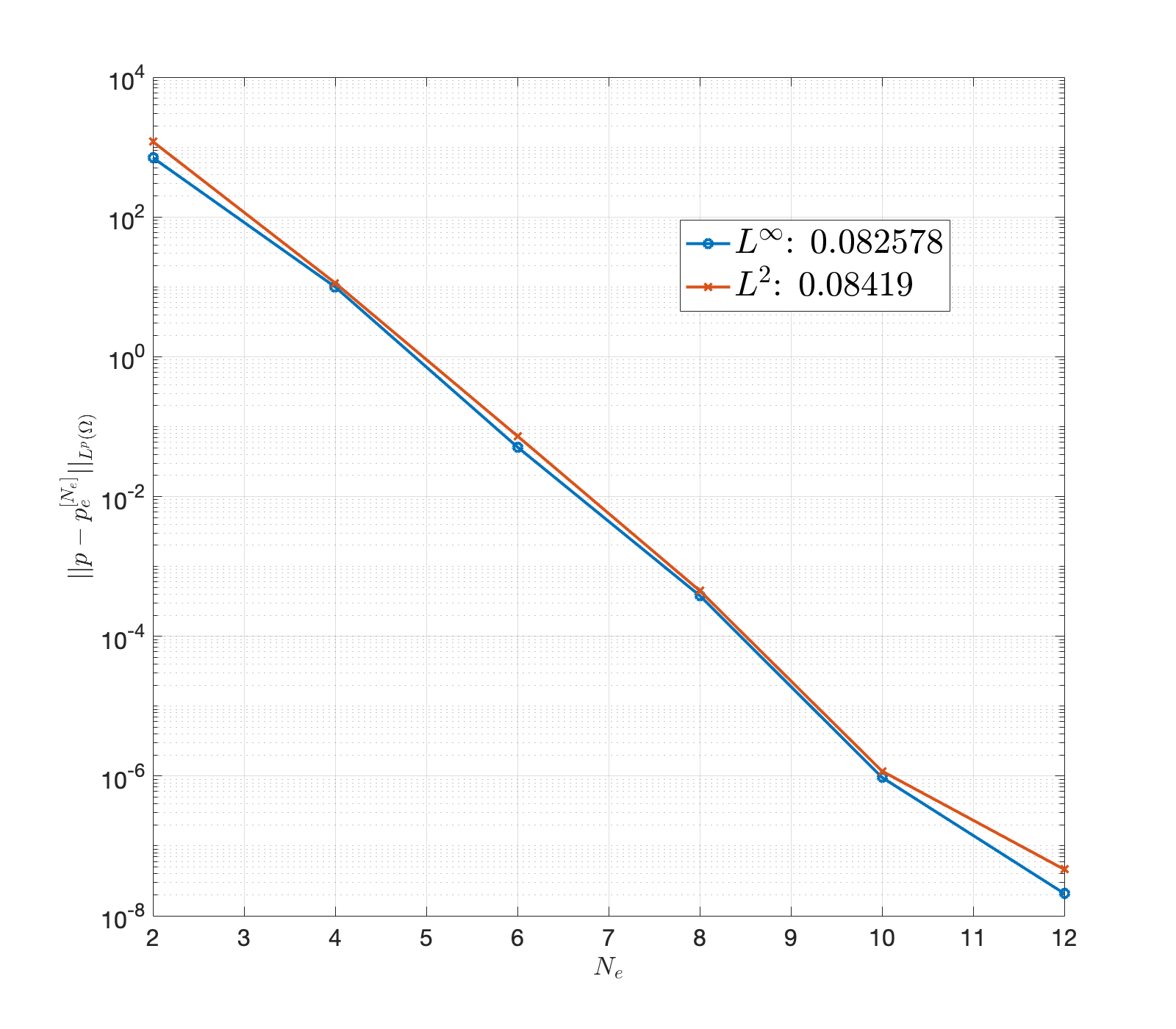}
		\label{NSex3}
	}
	\subfigure[Problem (\ref{NSEx4})]
	{\includegraphics[scale=\sclcm]{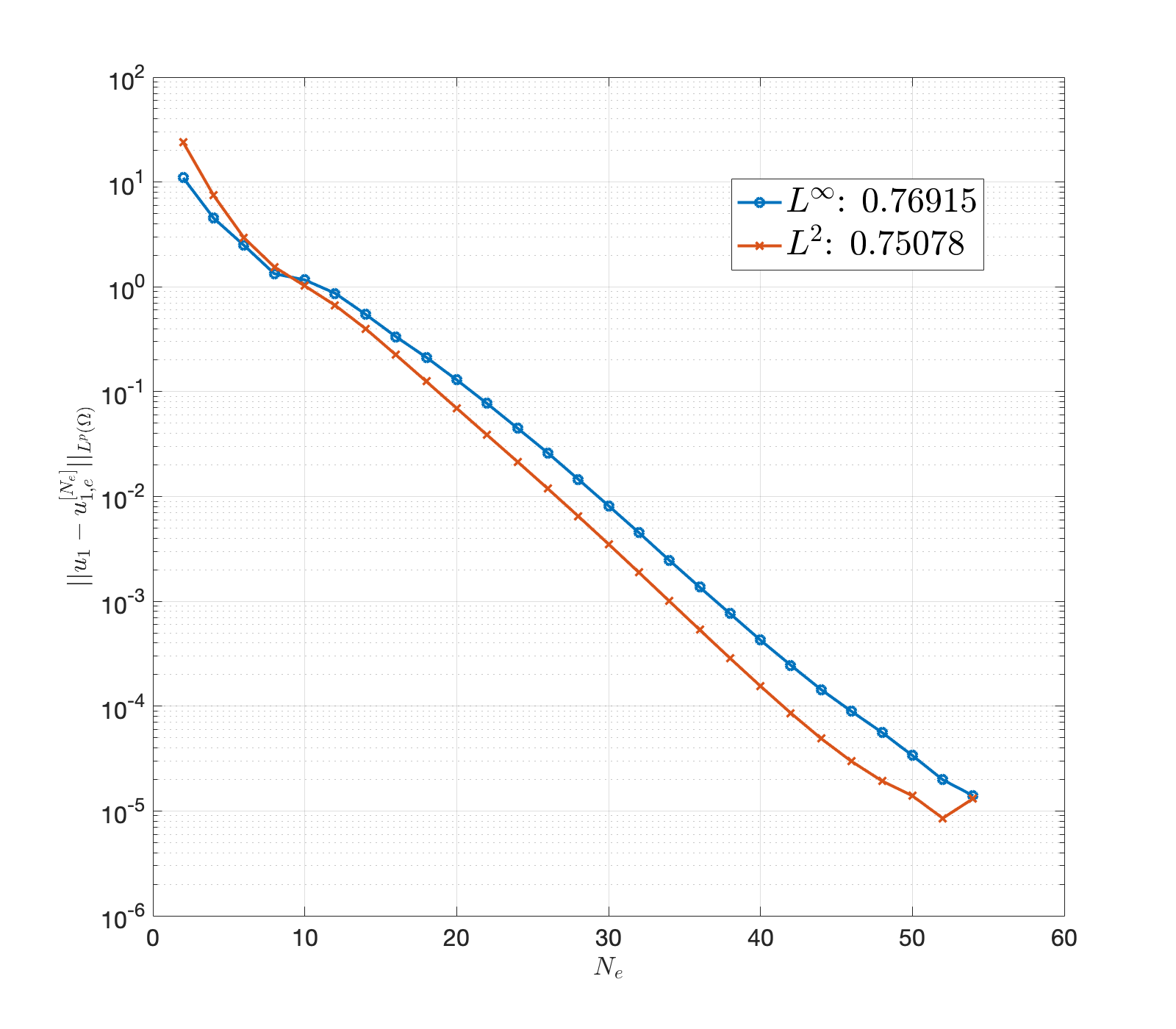}
		\includegraphics[scale=\sclcm]{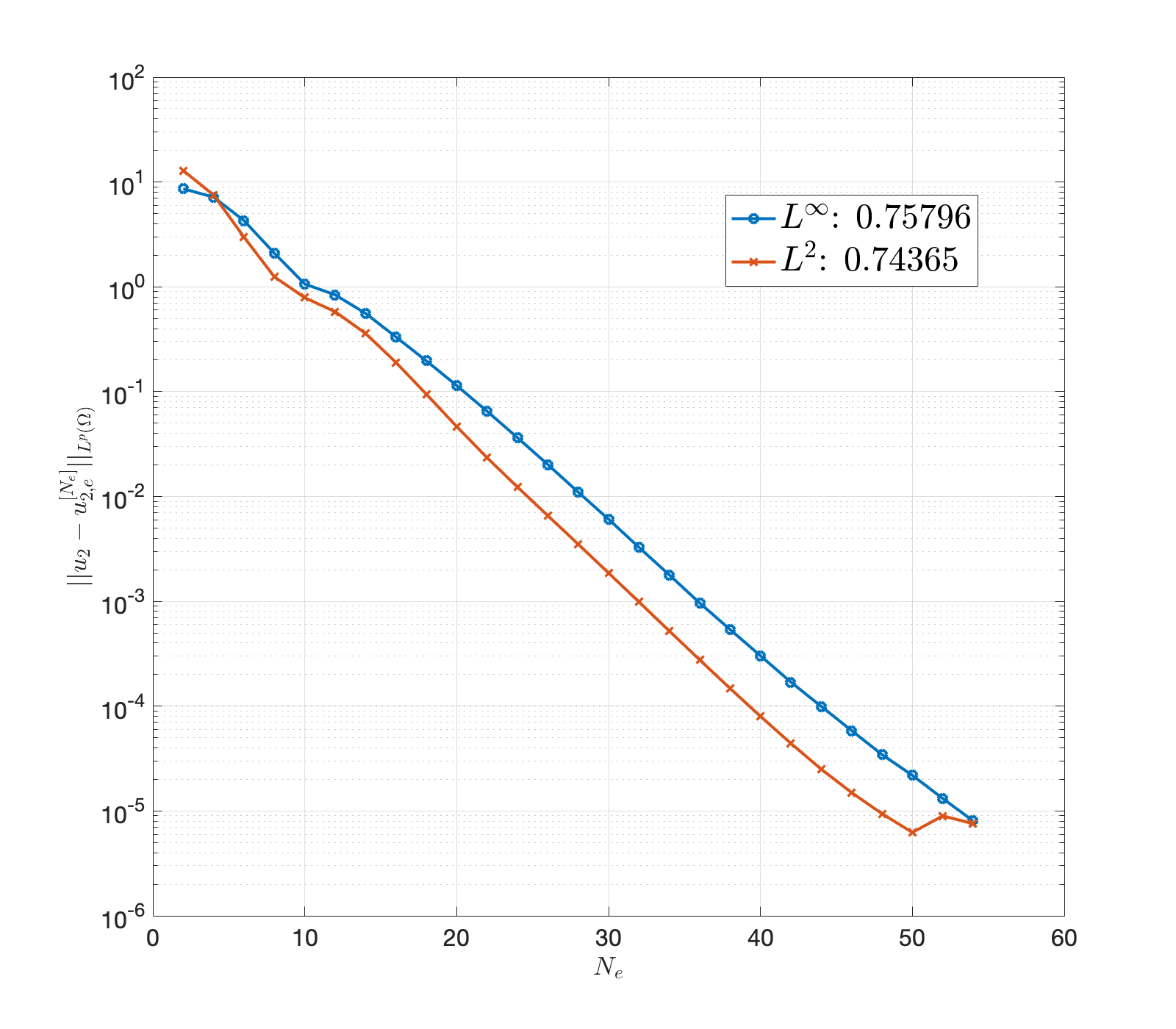}
		\includegraphics[scale=\sclcm]{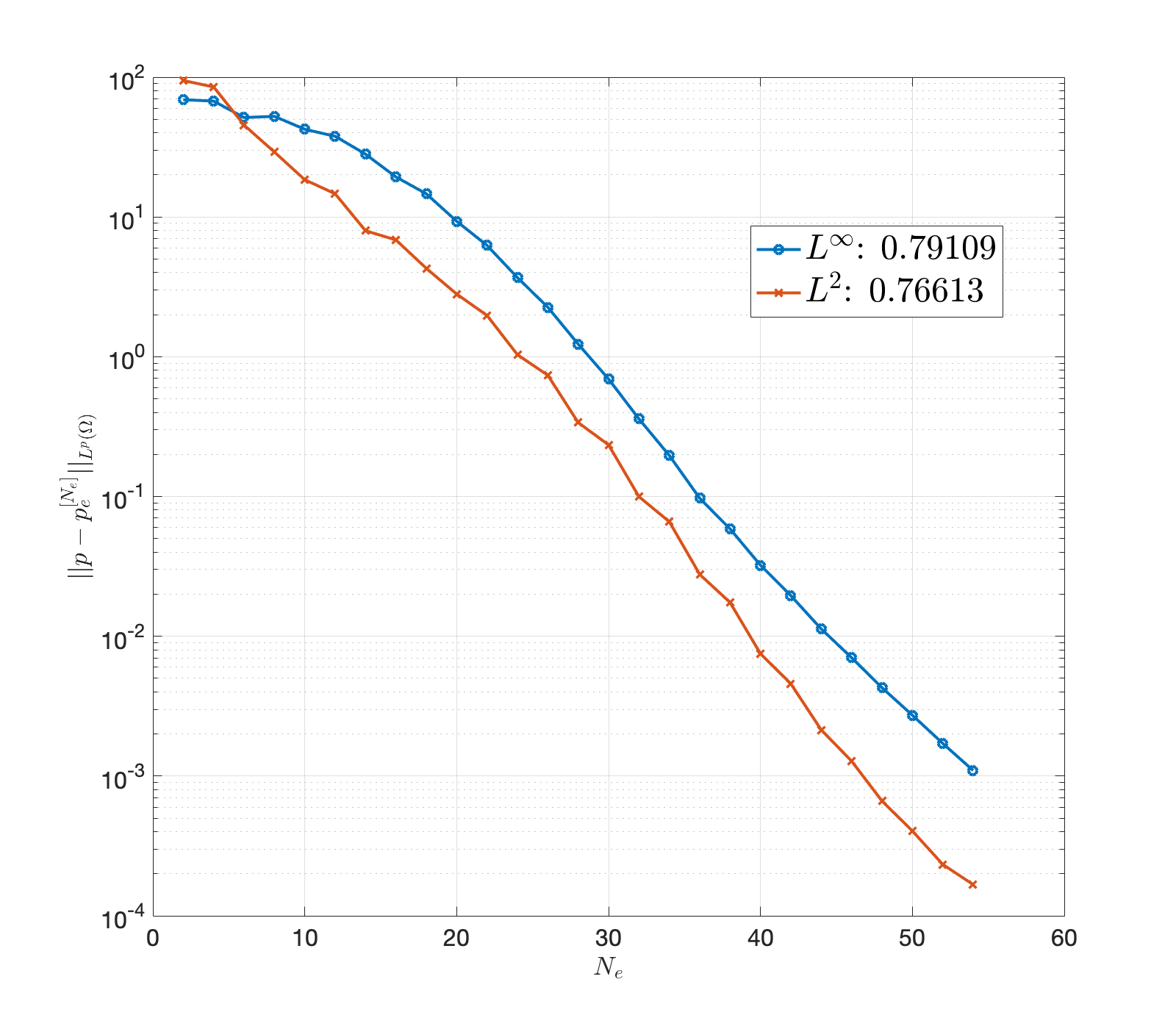}
		\label{NSex4}
	}
	\caption{Convergence plots for the solutions $u_1$, $u_2$ and $p$ to the Navier--Stokes test problems on the domain studied in (\ref{stokesEx3}). All the unknowns appear to converge exponentially at $O(a^{N_\text{e}})$ with the corresponding best-fit values of $a$ shown for the different norms.}
	\label{NSconvChannel}
\end{figure}

In a similar vein, we can use the alternative Stokes formulation in Subsection \ref{AltStokes} to solve (\ref{NSbdf4}) and (\ref{NScontraints}) over the same domain as in (\ref{stokesEx3}). We perform one test with a known exact solution and another with prescribed flow-rate and no-slip conditions; the errors in the latter case are calculated by comparing successively refined solutions. In both instances, we use the forward Euler method for initializing the BDF iterations. For the first test, we use the exact solution
\eqn{
	\tb{u}(t,x,y) =  \begin{pmatrix}
		4ye^{\sin(x)}(y^2-4) \\
		-e^{\sin(x)}\cos(x)(y^2-4)^2
	\end{pmatrix}e^{t}, \qquad p(t,x,y) = e^{2\cos(x)}\cos(y)e^{t}, \label{NSEx3}
} 
and solve it up to $T = 1$ using $\Delta t = 10^{-7}$ for forward Euler and $\Delta t = 10^{-5}$ for the BDF scheme, while the grid parameters are set at $N_x = 64$, $N_y = 72$. The solutions can be seen to converge exponentially to the true solution in Figure \ref{NSex3}.

For the second test, we consider
\eqn{
	\left\{ \begin{matrix}
		\tb{u}_t + (\tb{u} \cdot \nabla)\tb{u} = -\nabla p + \Delta \tb{u} + (\alpha,0)
		, &  \text{for } \tb{x} \in \Omega, t > 0,  \\
		\nabla \cdot \tb{u} = 0, &  \text{for } \tb{x} \in  \Omega, t > 0,  \\
		\tb{u} = \tb{0}, &   \text{for } \tb{x} \in \partial \Omega, t > 0, \\ 
		\int_{-2} ^2 u_1 \ dy = q, & \text{at } x = 0, \\
	\end{matrix} \right. 
	\label{NSEx4}
}
with the same initial condition as (\ref{NSEx3}). The problem is solved for $q = 1$ up to $T = 1$ with $\Delta t = 4 \times 10^{-5}$ throughout and $N_x = 256$ and $N_y = 190$. The error plots, displayed in Figure \ref{NSex4}, again exhibit exponential convergence, albeit at a much slower rate compared to Figure \ref{NSex3}.

\section{Viscoelastic Fluid Models}\label{ViscoFluids}
In this section, we apply our technique to solve a simple viscoelastic fluid flow problem. Such models typically consist of a pair of coupled systems describing the viscous and polymeric parts of the fluid. Numerical methods for such systems combine a Newtonian fluid solver with a scheme for the polymeric constitutive equations. Our effective techniques for viscous fluids, described in Section \ref{FluidProbs}, then only need to be combined with an appropriate treatment for the polymeric components to yield efficient methods for this class of problems.

The constitutive equations for viscoelastic fluids describe the evolution of the polymeric stresses $\bs{\tau}_\text{p}(t,\tb{x})$. Given the Weissenberg number $\text{Wi}$ (with units of time) and polymeric viscosity $\nu_\text{p}$, the Oldroyd--B model prescribes
\eqn{
	\pti \bs{\tau}_\text{p} + \left(\tb{u}\cdot\nabla\right)\bs{\tau}_\text{p} = \left(\nabla \tb{u}\right)\bs{\tau}_\text{p} + \bs{\tau}_\text{p}\left(\nabla \tb{u}\right)^T + \frac{1}{\text{Wi}}\left(\nu_\text{p} \dot{\bs{\gamma}} - \bs{\tau}_\text{p}\right), \label{TauEvol}
}
where $\dot{\bs{\gamma}} = \nabla \tb{u} + \left(\nabla \tb{u}\right)^T$ is the rate of strain tensor and $\left(\nabla \tb{u}\right)_{ij} = \partial_{x_j}u_i$. In the low Reynolds number regime, the fluid can be assumed to adjust instantaneously to changes in the polymeric stresses, in which case it can be modeled by the incompressible Stokes equation
\eqn{
	\left\{ \begin{matrix}
		-\nu_\text{s} \Delta \tb{u} + \nabla p = \nabla \cdot \bs{\tau}_\text{p},\\
		\nabla \cdot \tb{u} = 0,
	\end{matrix} \right. 
	\label{PolyStokes}
}
where $\nu_\text{s}$ is the solvent viscosity. Together with boundary and initial conditions, this provides a closed model of a viscoelastic fluid.

This system is conventionally re-formulated in terms of the non-dimensional conformation tensor $\bs{\sigma}(t,\bf{x})$, which is defined by
\eqn{
	\bs{\sigma} = \frac{\text{Wi}}{\nu_\text{p}}\bs{\tau}_\text{p} + \tb{I}, \label{PolyStConf}
}
where $\tb{I}$ is the identity matrix. Choosing typical length and velocity scales $X$ and $U$ respectively and defining
\eqn{
	\tilde{\tb{x}} = \frac{\tb{x}}{X}, \quad \tb{v} = \frac{\tb{u}}{U} \label{NonDim}
}
allows us to re-write (\ref{TauEvol}) as 
\eqn{
	\partial_{\tilde t} \bs{\sigma}  + \left(\tb{v}\cdot \tilde\nabla \right) \bs{\sigma} =  \left(\tilde \nabla \tb{v}\right) \bs{\sigma} + \bs{\sigma} \left(\tilde \nabla \tb{v}\right)^T + \frac{1}{\lambda}\left(\tb{I} - \bs{\sigma}\right), \label{ConfEvol}
}
and (\ref{PolyStokes}) as
\eqn{
	\left\{ \begin{matrix} 
		-\tilde \Delta \tb{v} + \tilde \nabla P = \xi \tilde \nabla \cdot \bs{\sigma}, \\
		\tilde \nabla \cdot \tb{v} = 0. \\
	\end{matrix}
	\right. \label{StokesEvol}
}

Here, $\lambda = \left(\frac{U}{X}\right)\text{Wi}$ and $\xi = \lambda^{-1}\left(\frac{\nu_\text{p}}{\nu_\text{s}}\right)$ are dimensionless quantities, and $\tilde t = t\left(\frac{U}{X}\right)$ and $P = p\left(\frac{X}{U\nu_\text{s}}\right)$ are appropriate rescalings of time and pressure respectively.

By definition, the conformation tensor is required to be positive semi-definite for all time. Evolving (\ref{ConfEvol}) numerically, however, leads invariably to the introduction of round-off and discretization errors that destroy this structure. To maintain this property, several re-formulations of the tensor evolution equation have been proposed. The more popular among these are the log-conformation representation of \cite{fattal2005time,hulsen2005flow} and the square-root approach pioneered in \cite{balci2011symmetric}. The former is specifically designed to avoid the exponential build-up due to the deformative terms in (\ref{ConfEvol}), hence eliminating a mechanism that may trigger the high Weissenberg number problem (HWNP). While our technique is well-suited to both, we have opted for the square-root formulation in the results shown here due to the computational advantages it offers, and also because this choice permits a direct comparison with the results in \cite{stein2019convergent}.

The square-root approach requires that, instead of (\ref{ConfEvol}), one evolve the square-root $\tb{b} = \sqrt {\bs{\sigma}}$ of the conformation tensor, leading to 
\eqn{
	\partial_{\tilde t} \tb{b} + \left(\tb{v}\cdot \tilde \nabla \right) \tb{b} = \tb{b}\left(\tilde \nabla \tb{v}\right)^T + \tb{a}\tb{b} + \frac{1}{2\lambda}\left(\left(\tb{b}^T\right)^{-1} - \tb{b}\right), \label{BEvol}
}
where $\tb{a}$ is any anti-symmetric matrix. By choosing $\tb{a}$ carefully, it can in addition be guaranteed that (\ref{BEvol}) evolves $\tb{b}$ so it remains a symmetric matrix, assuming it starts off as one, leading to significant savings in computational time and effort \cite{balci2011symmetric}.

We first consider the issues arising from the temporal discretization of (\ref{BEvol}). Despite the abundance of terms depending linearly on $\tb{b}$, the presence of velocity terms that change with time makes it prohibitive to treat most of the terms implicitly. We opt for the extrapolated-BDF-3 method, which is a three-step explicit scheme with the same memory and computational requirements as the Adams--Bashforth-3 method, with the added advantage of being monotonicity-preserving, assuming that it is initialized using the forward Euler method \cite{hundsdorfer2003monotonicity}. Among other properties, this confers total variation bounded (TVB) and positivity properties to our technique, and helps avoid introducing spurious oscillations that may manifest as undershoots and overshoots. We define 
\eqn{
	\tb{F}(\tb{b},\tb{v}) = -\left(\tb{v}\cdot \tilde \nabla \right) \tb{b} + \tb{b}\left(\tilde \nabla \tb{v}\right)^T + \tb{a}\tb{b} + \frac{1}{2\lambda}\left(\left(\tb{b}^T\right)^{-1} - \tb{b}\right), \label{BevF}
}
so the time-stepping routine reads
\eqn{
	\tb{b}^{n+1} = \frac{18\tb{b}^n - 9\tb{b}^{n-1} + 2\tb{b}^{n-2}}{11} + \left(\frac{6 \Delta t}{11}\right)\left(3\tb{F}\left(\tb{b}^{n},\tb{v}^{n} \right)  - 3\tb{F} \left( \tb{b}^{n-1},\tb{v}^{n-1} \right) + \tb{F}\left(\tb{b}^{n-2},\tb{v}^{n-2} \right) \right). \label{BEbdf}
}

This scheme can also be seen to be motivated by the IMEX approach employed in (\ref{NSbdf4}), except that all the terms are treated explicitly. We prefer the three-step scheme over the four-step version for the results presented here because it poses less severe time-step restrictions with no discernible loss of accuracy.

Evaluating (\ref{BEbdf}) at each time step requires that we be able to calculate all the terms in (\ref{BevF}). In particular, the computation of the partial derivatives of the components of $\tb{b}$ in the advection term in (\ref{BevF}) demands careful treatment as their values may only be known on an irregular domain $\Omega$. Taking our cue from the ideas underlying the PE method, we express all components of $\tb{b}$ in terms of a finite-dimensional basis on a larger computational domain $\Pi$, in essence computing their smooth extensions, and use the expansions to perform the necessary numerical operations. The coefficients $\bs{\beta} = \left(\beta^{(m)}_j\right)_{m \in \mathcal{M},j \in \mathcal{J}}$ of the expansions are chosen by minimizing the objective function
\eqn{
	H(\bs{\beta}) = \sum_{m \in \mathcal{M}} \norm{\sum_{j \in \mathcal{J}} \beta^{(m)}_j \eta_j - {b}^{(m)}}^2_{L^2(\Omega)}, \label{Hdef}
}
where $\left({b}^{(m)}\right)_{m \in \mathcal{M}} = \tb{b}$, and $\{\eta_j\}_{j \in \mathcal{J}}$ is an appropriate basis on $\Pi$. Barring any additional constraints, this basis can be chosen to be the same as the pressure basis constructed in Subsection \ref{AltStokes} without discounting the constant function. As a result, we do not need to develop any additional tools, and can simply re-use the machinery assembled for solving the associated Stokes system (\ref{StokesEvol}). This additional projection is also helpful in avoiding the build-up of large gradients across boundaries due to the non-physicality of the computed solutions outside $\Omega$, and echoes the stepwise re-extension procedure employed in \cite{stein2019convergent}.

\def \sclim{0.07}

\begin{figure}[tbph]
	\centering
	\subfigure[]
	{\includegraphics[scale=\sclbm]{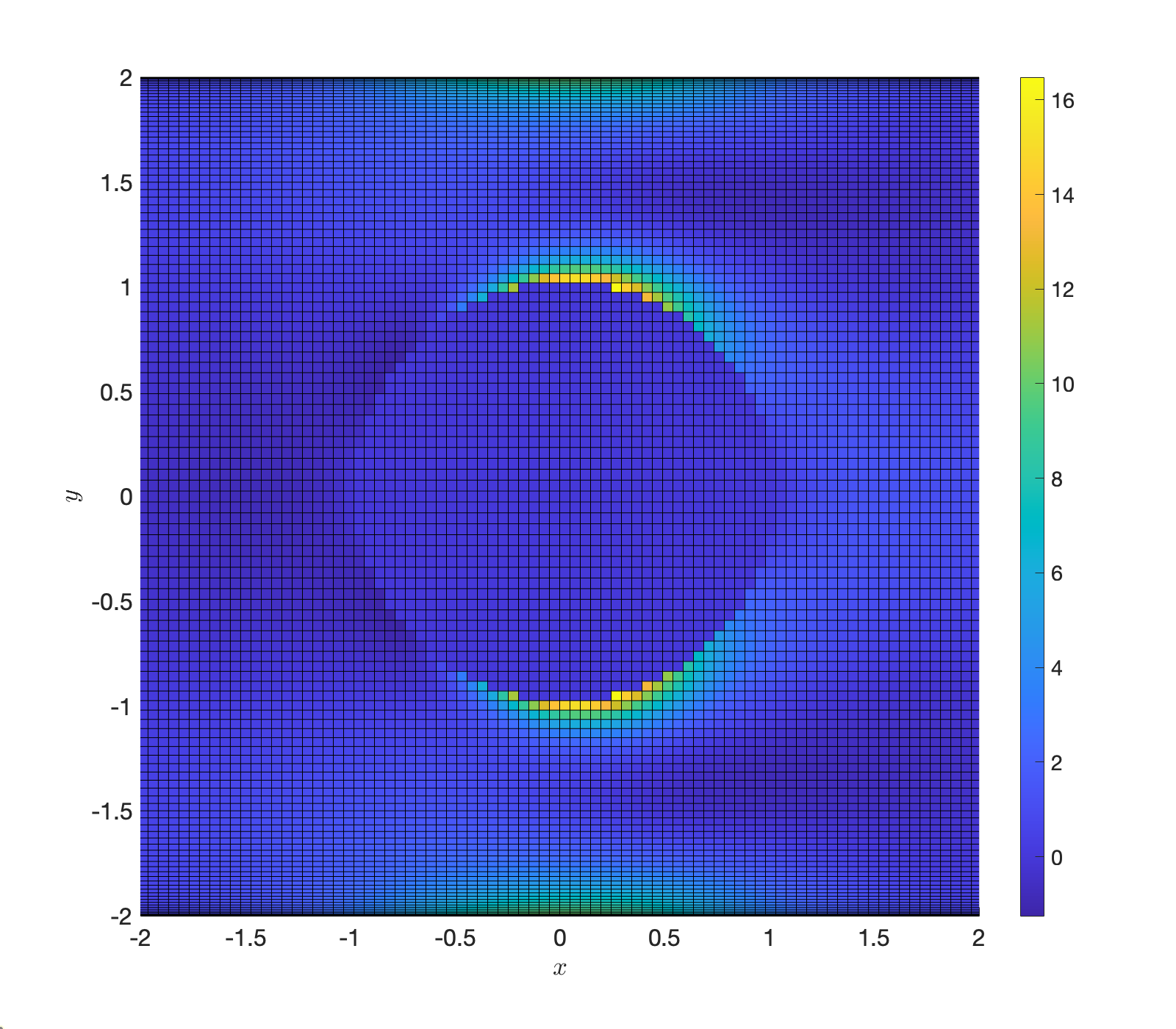}
		\label{IsobTPxx}
	}
	\subfigure[]
	{\includegraphics[scale=\sclbm]{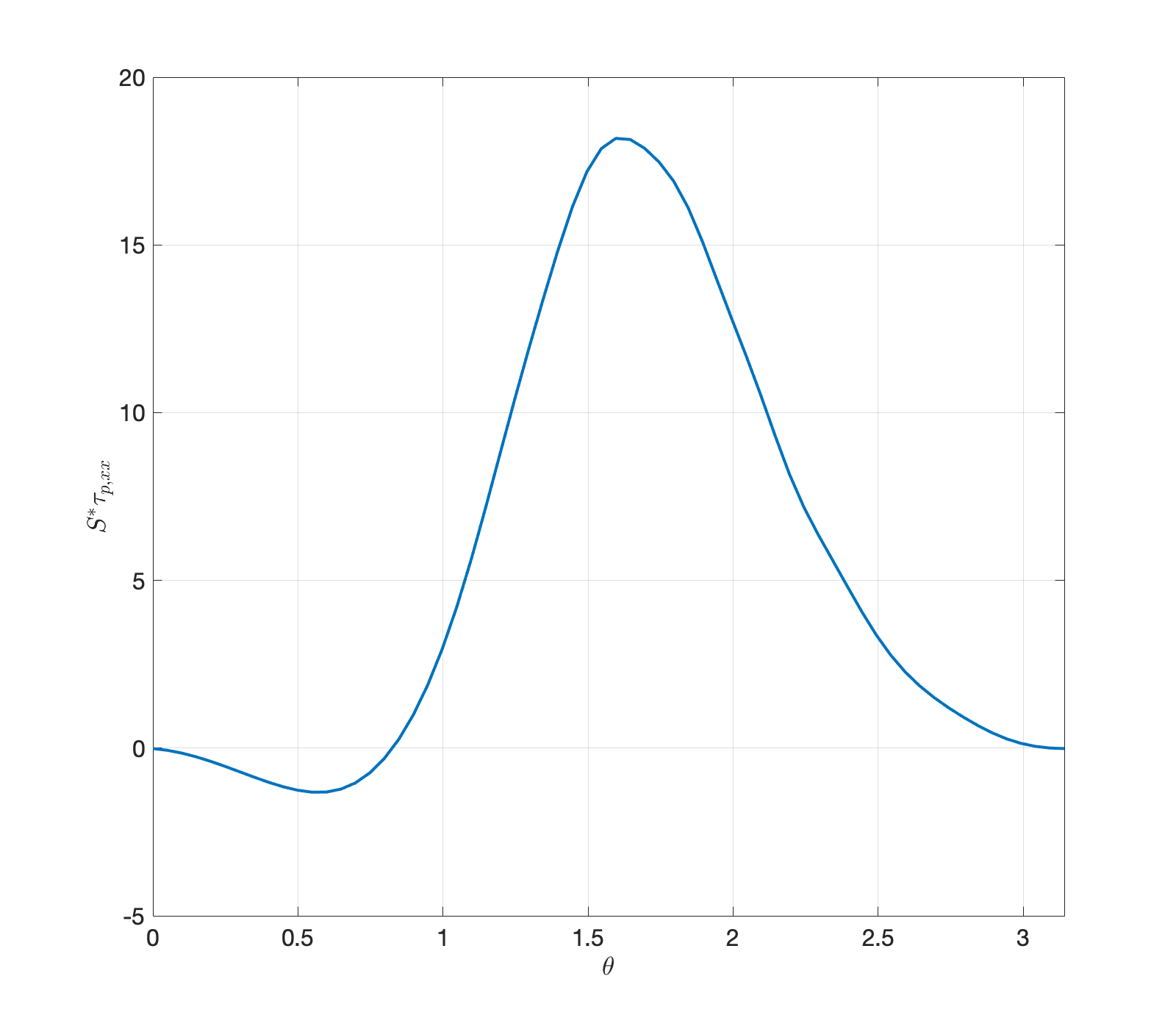}
		\label{IsobTPxxint}
	}
	\subfigure[]
	{\begin{tikzpicture}
		[node distance = 10mm,inner sep = 0pt ]
		\node (n0) at (0,0)  {\includegraphics[scale=\sclbm]{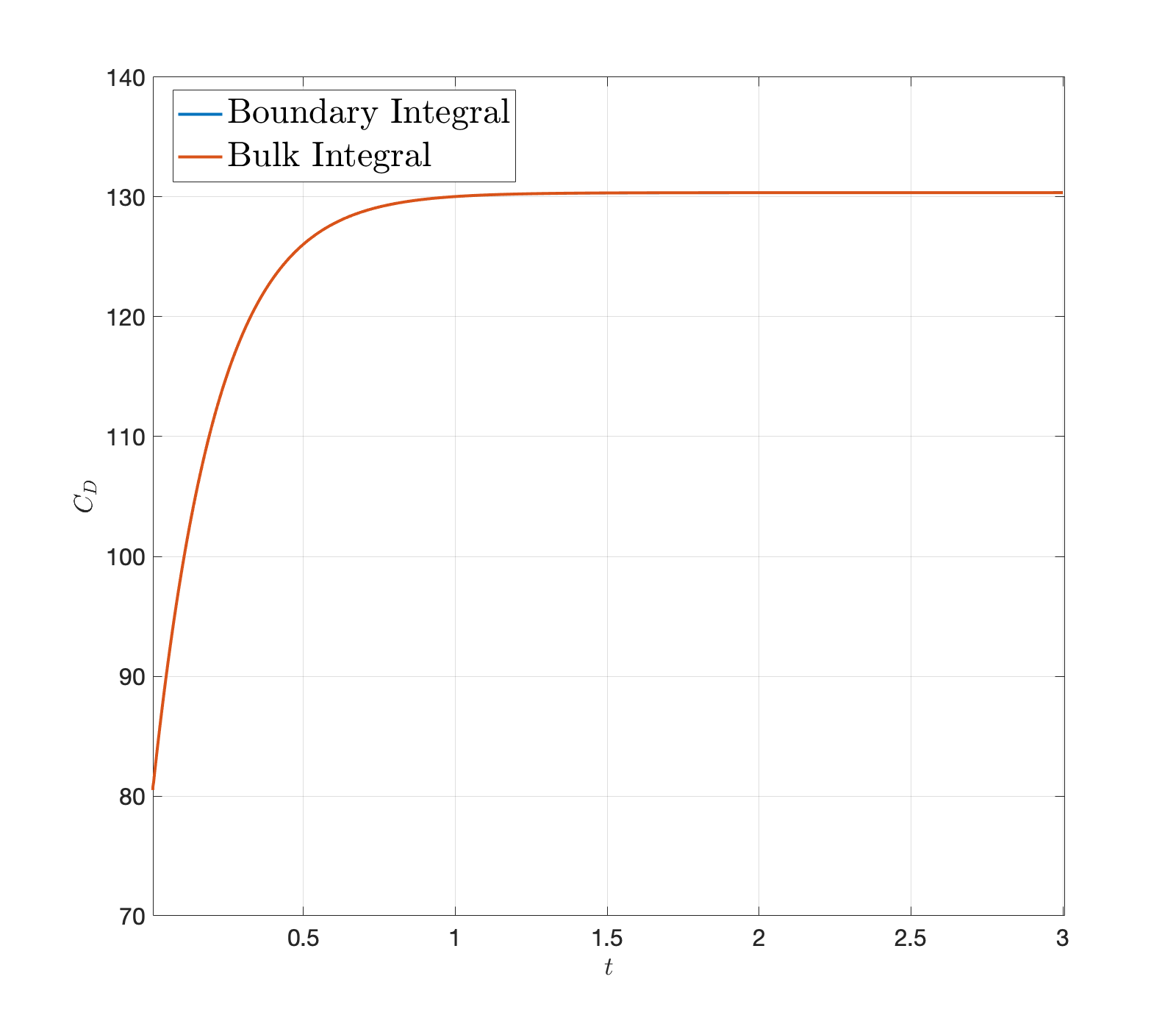}};
		\node (n3) at (1.2,-0.4)
		{\includegraphics[scale=\sclim]{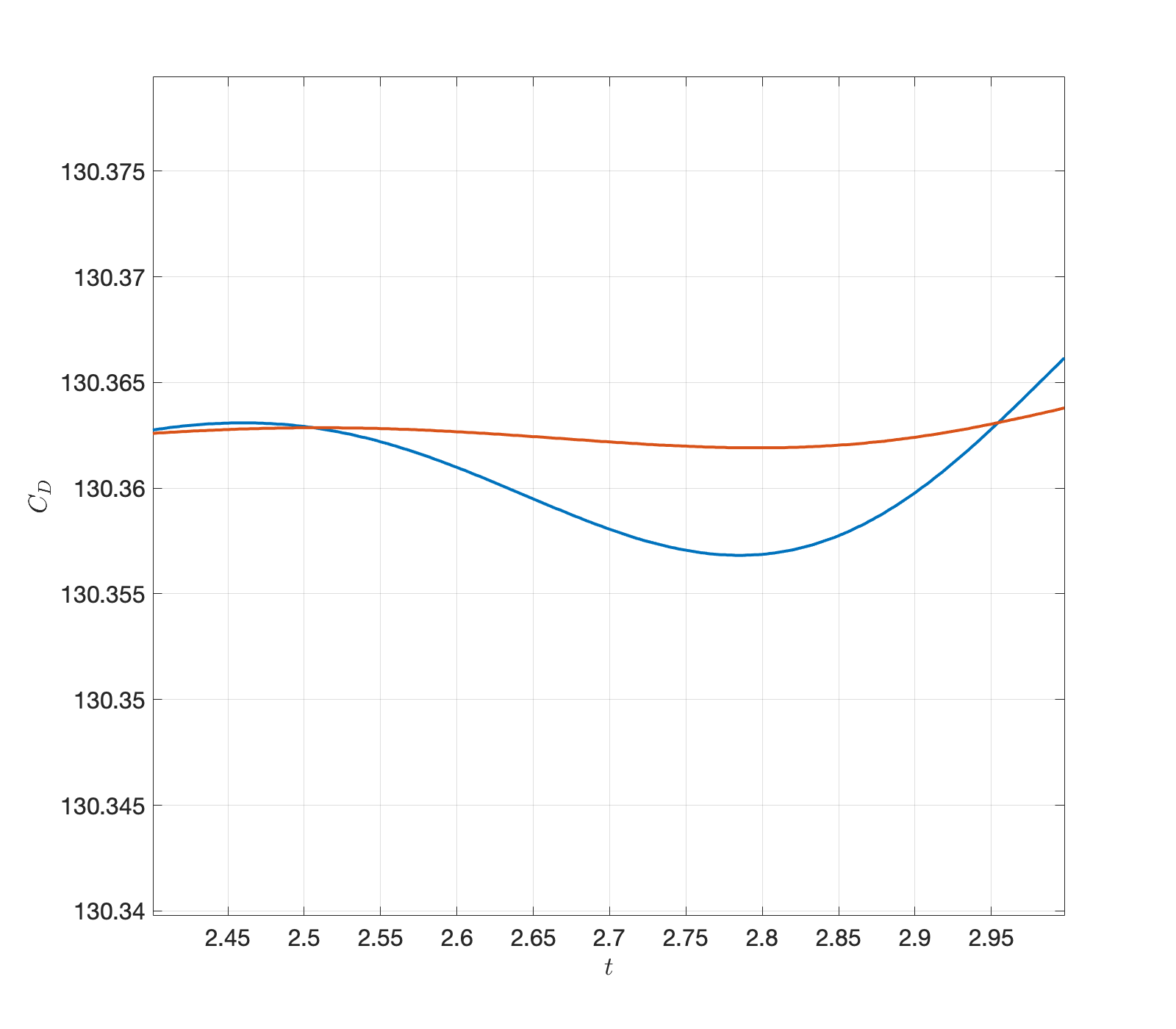}};   
		\draw[densely dashed]   (3.2,2.1) -- (2.8,1.05) 
		(1.85,2.1) -- (-0.25,1.05);              
		\end{tikzpicture}
		\label{IsobDrag}
	}
	\subfigure[]
	{\includegraphics[scale=\sclbm]{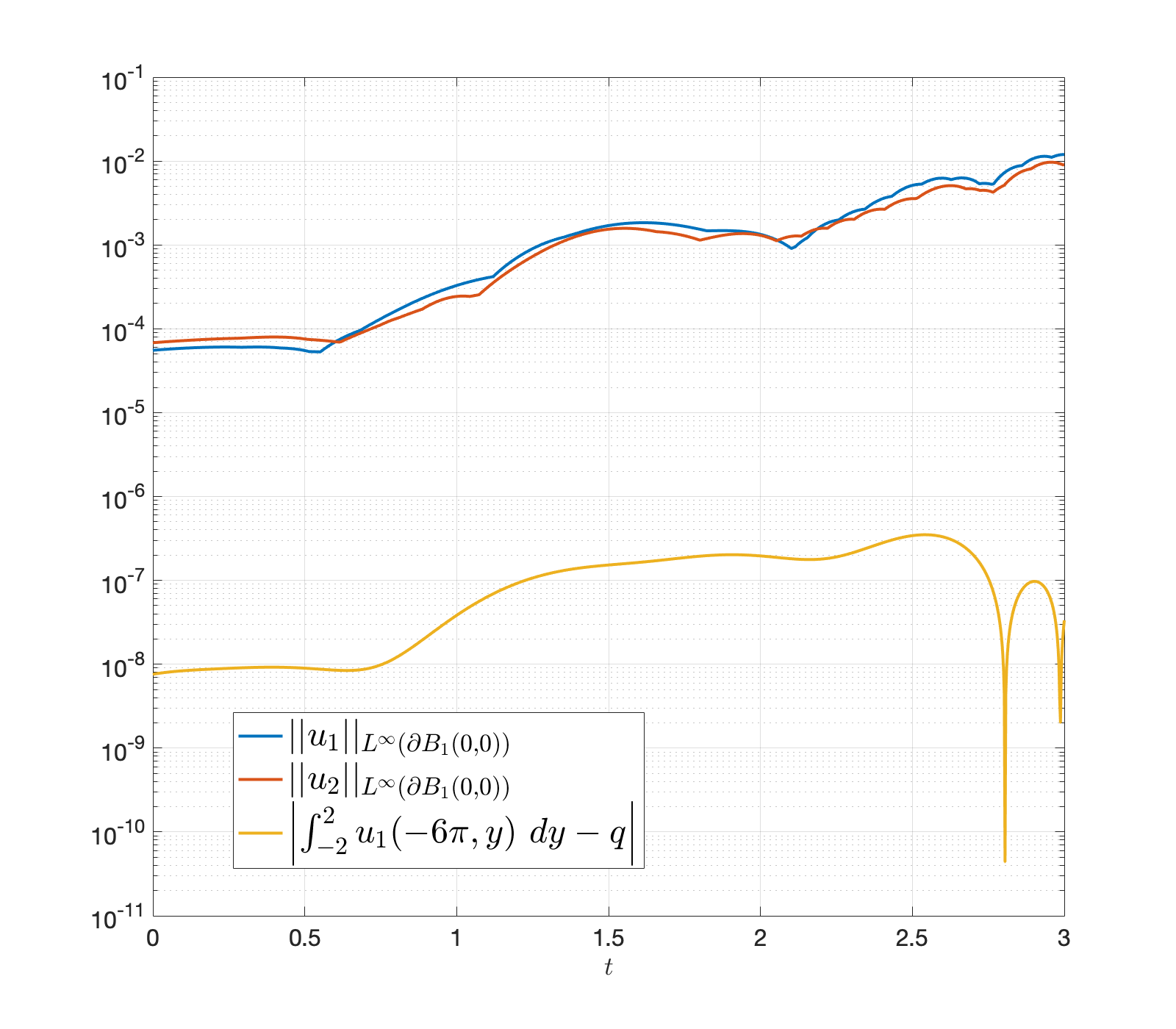}
		\label{IsobEIs}
	}
	\caption{(a) The $xx$ component of the polymeric stress tensor $\bs{\tau}_\text{p}$ for $\text{Wi} = 0.1$, shown at $t = 3$ for $-2 \leq x \leq 2$; the values are clearly maximized at the boundaries. (b) The tensor values from (a) interpolated to the circular boundary, with $\theta = 0$ denoting the point at the center of the channel upstream of the obstacle. (c) The evolution of the drag coefficient computed using both the boundary integral (\ref{DragCoeff}) and bulk integral (\ref{DragCoeff2}) approaches. Both techniques appear to agree broadly, with steady state achieved around $t = 2$. Zooming in, however, it can be seen that the bulk integral values are less oscillatory, hinting at the robustness of this approach. (d) The errors in the imposition of the no-slip and flow-rate conditions. While the latter is enforced fairly accurately, the no-slip constraints at the circular boundary leave a lot to be desired.}\label{Isob}
\end{figure}

We test this methodology by solving the benchmark problem of Stokes--Oldroyd--B flow past a circular obstacle. The channel is chosen to be an elongated version of the one studied in problem (\ref{stokesEx3}). We define $\Omega = ((-6\pi,6\pi) \times (-2,2)) - \overline{B_1(0,0)}$ and take $C = (-6\pi,6\pi) \times (-2,2)$ to be the computational domain. We identify $x = -6\pi$ and $x = 6\pi$, so the channel is again assumed to be periodic in $x$, and apply no-slip boundary conditions at $y = \pm 2$ and at the walls of the circular obstacle. A similar flow-rate condition $\int_{-2}^2 u_1 \ dy = q$ is also enforced at $x = -6\pi$, with an additional constant forcing term $(\alpha,0)$ included to act as a Lagrange multiplier to enforce this flow. As there is no inflow across any portion of the boundary at any time, i.e.,
\eqn{
	\partial \Omega_\text{in}(t) = \{\tb{x} \in \partial \Omega: \tb{u}(t,\tb{x}) \cdot \tb{n}(\tb{x}) < 0\} = \{\}, \label{InflowBoundary}
} 
we do not need to prescribe boundary conditions for the polymeric stresses. At $t = 0$, we set $\bs{\tau}_\text{p} = \tb{0}$ for all $\tb{x} \in \Omega$.

For our experiments, we set $\nu_\text{s} = 0.59$, $\nu_\text{p} = 0.41$, $R = 1$ and $q = 4$. The length and velocity scales in (\ref{NonDim}) are taken to be $X = 1$ (the radius of the disc) and $U = q/4 = 1$ (the average velocity across any vertical cross-section). The corresponding $\lambda$ then matches the definition of the non-dimensional Weissenberg number as used, e.g., in \cite{claus2013viscoelastic,hulsen2005flow}. The problem is solved till steady state is achieved, typically around $t = 20\text{Wi}$, at which point we compute the drag coefficient $C_D$ around the disc in two ways. Using the definition
\eqn{
	C_D = \int_{\partial B_1(0,0)} \hat{x} \cdot \left(\alpha \left(\hat{x} \otimes \hat{x}\right) +  \bs{\tau}_\text{p} - p\tb{I} + \nu_s\left(\nabla \tb{u} + \left(\nabla \tb{u}\right)^T\right)\right) \cdot \tb{n}(\tb{x}) \ ds \label{DragCoeff}
} 
provides one obvious route: interpolate all the terms in the integrand at the boundary and use the trapezoid rule to approximate (\ref{DragCoeff}). Our use of modal expansions reduces the task of interpolation to one of mere evaluation at the boundary nodes, considerably simplifying the task. The alternate method is obtained through the application of the divergence theorem to (\ref{DragCoeff}), yielding
\eqn{
	C_D = \iint_{B_1(0,0)} \di\left(\hat{x} \cdot \left(\alpha \left(\hat{x} \otimes \hat{x}\right) + \bs{\tau}_\text{p} - p\tb{I} + \nu_s\left(\nabla \tb{u} + \left(\nabla \tb{u}\right)^T\right)\right)\right) \ dA. \label{DragCoeff2}
}

The integrand can be read as the $x$-component of the divergence of the momentum current tensor. Conservation of momentum requires that this expression vanish identically on $\Omega$ (hence yielding the $x$-component of the Stokes equations (\ref{PolyStokes})); however, the integral in (\ref{DragCoeff2}) is over the interior of the obstacle, which is not part of $\Omega$, and hence may be non-zero. As our technique computes highly accurate extensions of all the terms in this integral, we can apply this procedure as well to compute the drag coefficient without any significant additional effort.

We used the time-step $\Delta t = 2.5 \times 10^{-3}$ for both the BDF and forward Euler steps, and set $N_x = 768$ and $N_y = 120$ as the grid parameters; this is coarser than the grid used for the Stokes and Navier--Stokes results in the earlier sections and is necessitated by the much longer channel required by this test. The highest frequencies of the extension functions were set at $(N_{\text{e},x},N_{\text{e},y}) = (240,40)$. 

The results for $\text{Wi} = \lambda = 0.1$ are shown in Figure \ref{Isob}. The $xx$ component of the polymeric stress is correctly maximized at the boundary, as noted by \cite{claus2013viscoelastic,stein2019convergent}. In addition, the trace of $\bs{\tau}_{p,xx}$ at the circular boundary is in close agreement with the benchmark results. The evolution of the drag coefficient values, calculated using both the boundary integral (\ref{DragCoeff}) and the bulk integral (\ref{DragCoeff2}), are shown in Figure \ref{IsobDrag} for $0 \leq t \leq 3$. Both approaches largely agree, indicating that steady state sets in at around $t = 2$. Focusing on the evolution further, however, reveals some oscillations, more so in the values obtained from (\ref{DragCoeff}). The more robust method (\ref{DragCoeff2}) yields $C_D = 130.3638$ at $t = 3$, which is in excellent agreement with earlier results \cite{claus2013viscoelastic,hulsen2005flow}.

To understand the origins of this oscillatory behavior, we assess the accuracy with which the no-slip and flow-rate conditions are satisfied. In Figure \ref{IsobEIs}, it can be seen that the latter constraint is met fairly accurately. On the other hand, the traces of both velocity components at the circular boundary increase in magnitude with time and, by the end, are quite large. In general, we have observed that boundary errors are a good indicator of the $L^{\infty}$ errors on the entire domain. The sizes of these errors suggest that our particular choices of discretization parameters provide only passable accuracy and are responsible to a large extent for the fluctuating behavior in Figure \ref{IsobDrag}. Nevertheless, the results obtained for Stokes and Navier--Stokes equations on similar domains give us good reason to be confident that, with minor modifications in our implementation, more accurate solutions can be computed for this problem too. Our current results can be seen as furnishing a proof of concept, in that our treatment of the conformation tensor yields solutions that are in good agreement with those obtained from more computationally intensive techniques.

\section{Conclusion}\label{concl}

In this paper, we have introduced an approach that yields spectrally accurate solutions to various PDEs on complex geometries. The technique fundamentally relies on the insight that functions can be, in principle, extended arbitrarily smoothly to more regular computational domains. Our method then searches for an extension in a finite dimensional space, ensuring the faithfulness of the equation and the boundary conditions by minimizing an objective function. By treating only the linear parts of the PDE implicitly and handling any nonlinear components explicitly, the minimization reduces to a simple linear least squares problem. As the dimension of the projection space is increased, the growth in its descriptive power enables highly accurate solutions to be computed for fairly challenging problems on a wide array of domains.

Our approach permits two formulations: the first explicitly extends the forcing terms while the second one seeks to do the same to the solution directly. While the two are equivalent if the extension functions are eigenfunctions of the differential operator, in a more general setting, they both have particular strengths. The forcing extension framework allows for, in principle, non-differentiable functions to be used as extension functions, e.g., piecewise polynomials. The solution extension approach, in contrast, requires the extension functions to possess some minimal regularity but is significantly easier to set up for incompressible fluid problems, and makes even daunting problems such as fluid flow on curved domains tractable.

On top of these advantages, our method is demonstrably stable under iterative application, which is a crucial requirement for solving time dependent problems. This is a consequence of the fact that our extension methodology relies only on the values on the physical domain and does not carry forward the non-physical values to the following steps. This is a significant improvement on the SFE approach which implicitly makes use of the non-physical values while applying the regularity constraints. Moreover, that technique requires these conditions to be imposed at boundary nodes, with the higher order constraints necessarily requiring normal vectors to be associated with each node. In cases where the boundary is rough and unsymmetrical, choosing these vectors is far from obvious. The PE technique, on the other hand, is completely free of these shortcomings, with its reliance on the trapezoid rule for boundary integrals lending increased accuracy. In addition, the linear systems are much easier and faster to set up than for the SFE technique as we make use of the full extension functions instead of truncating them to the extension domain and applying the inverse differential operators.

Our technique also shows promise for viscoelastic problems. A rudimentary application of our approach allows us to compute the velocity and stress profiles for a benchmark problem with adequate accuracy, with the drag coefficient values agreeing very well with the results in the literature. A more careful treatment is required to ensure greater accuracy and to permit the computation of viscoelastic flow at higher Weissenberg numbers. The most severe bottleneck in our approach is the highly taxing memory requirement due to the large number of nodes in the long discretized channel. Overcoming this impediment requires taking advantage of the structure of the domain: as most of it is composed essentially of a long rectangular channel, it is reasonable to partition it into multiple components and apply our technique locally. This would, in turn, entail replacing the Fourier basis in the $x$-direction by a more flexible polynomial basis and handling communication across the element boundaries. We anticipate further savings in memory and computational requirements by specifically tailoring the QR decomposition to our problem. Finally, we are confident that allying our ideas with the log-conformation formulation of the tensor evolution will lead to stable and accurate solutions to viscoelastic flow at higher Weissenberg numbers.

\section{Acknowledgements}
This work is supported by the Department of Energy (DOE) Office of Advanced Scientific Computing Research (ASCR) through the ASCR Distinguished Computational Mathematics Postdoc Project under ASCR Project 71268, and by the National Science Foundation awards DMS 1664645, OAC 1450327, OAC 1652541, OAC 1931516, and CBET 1944156. The authors thank Aleksandar Donev, Robert D. Guy, Becca Thomases, M. Gregory Forest, and J. Thomas Beale for their insightful suggestions and comments. 

\bibliography{refs}

\begin{thebibliography}{10}
\expandafter\ifx\csname url\endcsname\relax
  \def\url#1{\texttt{#1}}\fi
\expandafter\ifx\csname urlprefix\endcsname\relax\def\urlprefix{URL }\fi
\expandafter\ifx\csname href\endcsname\relax
  \def\href#1#2{#2} \def\path#1{#1}\fi

\bibitem{boyd2001chebyshev}
J.~P. Boyd, Chebyshev and Fourier spectral methods, Courier Corporation, 2001.

\bibitem{gottlieb1977numerical}
D.~Gottlieb, S.~A. Orszag, Numerical analysis of spectral methods: theory and
  applications, Vol.~26, Siam, 1977.

\bibitem{bernardi1997spectral}
C.~Bernardi, Y.~Maday, Spectral methods, Handbook of numerical analysis 5
  (1997) 209--485.

\bibitem{canuto2007spectral}
C.~Canuto, M.~Y. Hussaini, A.~Quarteroni, T.~A. Zang, Spectral methods:
  fundamentals in single domains, Springer Science \& Business Media, 2007.

\bibitem{canuto2012spectral}
C.~Canuto, M.~Y. Hussaini, A.~Quarteroni, A.~Thomas~Jr, et~al., {Spectral
  methods in fluid dynamics}, Springer Science \& Business Media, 2012.

\bibitem{karniadakis2013spectral}
G.~Karniadakis, S.~Sherwin, Spectral/hp element methods for computational fluid
  dynamics, Oxford University Press, 2013.

\bibitem{patera1984spectral}
A.~T. Patera, A spectral element method for fluid dynamics: laminar flow in a
  channel expansion, Journal of computational Physics 54~(3) (1984) 468--488.

\bibitem{peskin1972flow}
C.~S. Peskin, Flow patterns around heart valves: a numerical method, Journal of
  computational physics 10~(2) (1972) 252--271.

\bibitem{peskin1977numerical}
C.~S. Peskin, Numerical analysis of blood flow in the heart, Journal of
  computational physics 25~(3) (1977) 220--252.

\bibitem{leveque1994immersed}
R.~J. Leveque, Z.~Li, {The immersed interface method for elliptic equations
  with discontinuous coefficients and singular sources}, SIAM Journal on
  Numerical Analysis 31~(4) (1994) 1019--1044.

\bibitem{li2006immersed}
Z.~Li, K.~Ito, The immersed interface method: numerical solutions of PDEs
  involving interfaces and irregular domains, SIAM, 2006.

\bibitem{stein2016immersed}
D.~B. Stein, R.~D. Guy, B.~Thomases, {Immersed boundary smooth extension: a
  high-order method for solving PDE on arbitrary smooth domains using Fourier
  spectral methods}, Journal of Computational Physics 304 (2016) 252--274.

\bibitem{stein2017immersed}
D.~B. Stein, R.~D. Guy, B.~Thomases, {Immersed Boundary Smooth Extension
  (IBSE): A high-order method for solving incompressible flows in arbitrary
  smooth domains}, Journal of Computational Physics 335 (2017) 155--178.

\bibitem{spagnolie2015complex}
S.~E. Spagnolie, Complex fluids in biological systems, Biological and Medical
  Physics, Biomedical Engineering.

\bibitem{bruno2010high}
O.~P. Bruno, M.~Lyon, {High-order unconditionally stable FC-AD solvers for
  general smooth domains I. Basic elements}, Journal of Computational Physics
  229~(6) (2010) 2009--2033.

\bibitem{lyon2010high}
M.~Lyon, O.~P. Bruno, {High-order unconditionally stable FC-AD solvers for
  general smooth domains II. Elliptic, parabolic and hyperbolic PDEs;
  theoretical considerations}, Journal of Computational Physics 229~(9) (2010)
  3358--3381.

\bibitem{agress2019smooth}
D.~Agress, P.~Guidotti, D.~Yan, {The Smooth Selection Embedding Method with
  Chebyshev Polynomials}, arXiv preprint arXiv:1902.03713.

\bibitem{agress2021smooth}
D.~J. Agress, P.~Q. Guidotti, The smooth extension embedding method, SIAM
  Journal on Scientific Computing 43~(1) (2021) A446--A471.

\bibitem{klinteberg2019fast}
L.~af~Klinteberg, T.~Askham, M.~C. Kropinski, {A fast integral equation method
  for the two-dimensional Navier-Stokes equations}, Journal of Computational
  Physics 409 (2020) 109353.

\bibitem{fryklund2019integral}
F.~Fryklund, M.~C.~A. Kropinski, A.-K. Tornberg, An integral equation based
  numerical method for the forced heat equation on complex domains, arXiv
  preprint arXiv:1907.08537.

\bibitem{fryklund2018partition}
F.~Fryklund, E.~Lehto, A.-K. Tornberg, Partition of unity extension of
  functions on complex domains, Journal of Computational Physics 375 (2018)
  57--79.

\bibitem{qadeer2021smooth}
S.~Qadeer, B.~E. Griffith, The smooth forcing extension method: A high-order
  technique for solving elliptic equations on complex domains, Journal of
  Computational Physics (2021) 110390.

\bibitem{boyd2002comparison}
J.~P. Boyd, {A comparison of numerical algorithms for Fourier extension of the
  first, second, and third kinds}, Journal of Computational Physics 178~(1)
  (2002) 118--160.

\bibitem{henle2010hydrodynamics}
M.~L. Henle, A.~J. Levine, Hydrodynamics in curved membranes: The effect of
  geometry on particulate mobility, Physical Review E 81~(1) (2010) 011905.

\bibitem{saffman1975brownian}
P.~Saffman, M.~Delbr{\"u}ck, Brownian motion in biological membranes,
  Proceedings of the National Academy of Sciences 72~(8) (1975) 3111--3113.

\bibitem{hundsdorfer2007imex}
W.~Hundsdorfer, S.~J. Ruuth, {IMEX} extensions of linear multistep methods with
  general monotonicity and boundedness properties, Journal of Computational
  Physics 225~(2) (2007) 2016--2042.

\bibitem{fattal2005time}
R.~Fattal, R.~Kupferman, Time-dependent simulation of viscoelastic flows at
  high weissenberg number using the log-conformation representation, Journal of
  Non-Newtonian Fluid Mechanics 126~(1) (2005) 23--37.

\bibitem{hulsen2005flow}
M.~A. Hulsen, R.~Fattal, R.~Kupferman, Flow of viscoelastic fluids past a
  cylinder at high weissenberg number: stabilized simulations using matrix
  logarithms, Journal of Non-Newtonian Fluid Mechanics 127~(1) (2005) 27--39.

\bibitem{balci2011symmetric}
N.~Balci, B.~Thomases, M.~Renardy, C.~R. Doering, Symmetric factorization of
  the conformation tensor in viscoelastic fluid models, Journal of
  Non-Newtonian Fluid Mechanics 166~(11) (2011) 546--553.

\bibitem{stein2019convergent}
D.~B. Stein, R.~D. Guy, B.~Thomases, {Convergent solutions of Stokes Oldroyd-B
  boundary value problems using the Immersed Boundary Smooth Extension (IBSE)
  method}, Journal of Non-Newtonian Fluid Mechanics 268 (2019) 56--65.

\bibitem{hundsdorfer2003monotonicity}
W.~Hundsdorfer, S.~J. Ruuth, R.~J. Spiteri, Monotonicity-preserving linear
  multistep methods, SIAM Journal on Numerical Analysis 41~(2) (2003) 605--623.

\bibitem{claus2013viscoelastic}
S.~Claus, T.~N. Phillips, Viscoelastic flow around a confined cylinder using
  spectral/hp element methods, Journal of Non-Newtonian Fluid Mechanics 200
  (2013) 131--146.

\end{thebibliography}

\end{document}